\documentclass[letterpaper, 10 pt, conference]{ieeeconf}  

\IEEEoverridecommandlockouts                              

\usepackage{array}

\usepackage{cite}
\usepackage{amsmath,amsthm,amssymb,amsfonts}

\usepackage{xpatch}
\makeatletter
\xpatchcmd{\@thm}{\thm@headpunct{.}}{\thm@headpunct{}}{}{}
\makeatother

\usepackage{bm}
\usepackage[T1]{fontenc}

\usepackage{algorithm,algorithmic}
\usepackage{graphicx}
\usepackage{textcomp}
\usepackage{xcolor}

\usepackage[short]{optidef}
\usepackage[colorlinks=false, pdfborderstyle={/S/0/0 4}
, linkbordercolor={1 0 0}, citebordercolor={0 1 0}, urlbordercolor={0 0 1}]{hyperref}

\usepackage{aliascnt}

\newtheorem{theorem}{Theorem}
\newtheorem{assumption}{Assumption}

\newtheorem{lemma}{Lemma}

\newtheorem{remark}{Remark}

\DeclareMathOperator{\col}{col}

\DeclareMathOperator{\diag}{diag}

\def\bsg{{\mathbf{g}}}

\def\bsq{{\mathbf{q}}}

\def\bsv{{\mathbf{v}}}

\def\bsx{{\mathbf{x}}}
\def\bsy{{\mathbf{y}}}
\def\bsz{{\mathbf{z}}}

\def\bsE{{\mathbf{E}}}
\def\bsF{{\mathbf{F}}}

\def\bsH{{\mathbf{H}}}

\def\bsL{{\mathbf{L}}}

\def\bsM{{\mathbf{M}}}

\usepackage{fancyhdr}
\newtheorem{lemmaalt}{Lemma}[lemma]
\newenvironment{lemmap}[1]{
  \renewcommand\thelemmaalt{#1}
  \lemmaalt
}{\endlemmaalt}

\allowdisplaybreaks 

\usepackage{colortbl}
\definecolor{LightGray}{gray}{0.9}
\definecolor{LightCyan}{rgb}{0.4863    0.6157    0.5922}
\definecolor{LightGreen}{rgb}{0.8824    0.8000    0.6941}
\definecolor{LightOrange}{rgb}{0.8314    0.7294    0.6784}

\usepackage{multirow}
\usepackage{array}
\newcolumntype{M}[1]{>{\centering\arraybackslash}m{#1}}
\newcolumntype{N}{@{}m{0pt}@{}}

\usepackage{fancyhdr}
\title{\LARGE \bf
Compressed Zeroth-Order Algorithm for Stochastic Distributed Nonconvex Optimization}

\author{Haonan Wang$^1$, Xinlei Yi$^{1,2}$, Yiguang Hong$^{1,2}$ 
\thanks{This work was supported by the National Natural Science Foundation of China 62088101.}
\thanks{$^1$Department of Control Science and Engineering, College of Electronics and Information Engineering, Tongji University, Shanghai, 201804, China.
        {\tt\small\{hnwang, xinleiyi, yghong\}@tongji.edu.cn}.}%
\thanks{$^2$Shanghai Research Institute for Intelligent Autonomous Systems, Tongji University, Shanghai 201210, China.}%
}
\begin{document}


\maketitle
\thispagestyle{empty}
\pagestyle{empty}

\title{\LARGE \bf
Compressed Zeroth-Order Algorithm for Stochastic Distributed Nonconvex Optimization}
\begin{abstract}
	This paper studies the stochastic distributed nonconvex optimization problem over a network of agents,  where agents only access stochastic zeroth-order information about their local cost functions and collaboratively optimize the global objective over bandwidth-limited communication networks. 
	To mitigate communication overhead and handle the unavailability of explicit gradient information, we propose a communication compressed zeroth-order stochastic distributed (CZSD) algorithm.
	 By integrating a generalized contractive compressor and a stochastic two-point zeroth-order oracle, CZSD achieves convergence rates comparable to its exact communication counterpart while reducing both communication overhead and sampling complexity. 
	  Specifically,  to the best of our knowledge, CZSD is the first compressed zeroth-order algorithm achieving linear speedup, with convergence rates of $\mathcal{O}(\sqrt{p}/\sqrt{nT})$ and $\mathcal{O}(p/(nT))$ under general nonconvex settings and  the Polyak--{\L}ojasiewicz condition, respectively.
	   Numerical experiments validate the algorithm's effectiveness and communication efficiency.
\end{abstract}

\section{INTRODUCTION}

Distributed nonconvex optimization has recently attracted considerable attention due to its wide applications, including machine learning \cite{chang2020distributed} and multi-agent systems \cite{guo2016case}.
Typically, distributed optimization involves multiple agents collaboratively optimizing a global objective function through information exchange and local computation, 
thus facing two major practical challenges: substantial communication overhead among agents, and difficulty in explicitly obtaining local gradient information.

To reduce communication overhead, several compression-based distributed nonconvex algorithms have been proposed. 
A key strategy involves the development of various compression techniques, including unbiased compressors \cite{alistarh2017qsgd}, standard uniform quantizers \cite{xu_quantized_2024}, contractive compressors \cite{zhao2022beer}, and sparsification methods \cite{lin2018deep}.
Another line of research aims to unify these compressors. Specifically, \cite{beznosikov2023biased} provided a comprehensive review of biased compressors; \cite{Yi_CommunicationCompression_2023} introduced three classes of general compressors including globally and locally bounded absolute error compressors; and \cite{Liao_RobustCompressed_2024} proposed a unified compressor incorporating both relative and absolute compression errors. 
Additionally, several studies focus on accelerating convergence rates of compressed algorithms, demonstrating convergence performance comparable to exact communication. 
For instance, \cite{pmlr-v119-lu20a} and \cite{pmlr-v97-yu19d} achieved linear speedup for stochastic gradient descent (SGD) algorithms
through the error compensation and momentum methods, respectively. Recent papers
\cite{zhao2022beer,Liao_RobustCompressed_2024,Yi_CommunicationCompression_2023} employed error compensation/feedback and dynamic scaling techniques to obtain the well-known $\mathcal{O}(1/T)$ convergence rate for first-order algorithms,
where $T$ is the iteration number. 
These convergence guarantees, comparable to exact communication counterparts, exhibit greater potential for more flexible algorithmic extensions and broader applicability of compression methods.

In scenarios where agents cannot explicitly obtain gradient information due to the complexity of the objective functions, integrating zeroth-order methods with compressed distributed algorithms becomes a natural choice. However, this direct integration poses significant challenges due to coupled errors from gradient estimation and compression.
Specifically, compressed zeroth-order algorithms either exhibit slower convergence compared to standard zeroth-order methods or impose stringent conditions. 
For convex optimization algorithms, \cite{ding_distributed_consensus_2017} established convergence without explicit convergence rates. 
Under a relatively strong Lipschitz condition, \cite{Hua_CZOO_2025} achieved an $\mathcal{O}(\sqrt{p/T})$ convergence rate by employing a variance-reduction technique, where $p$ denotes the dimension of the optimization variables.
 The method in \cite{ding_quantized_gradient_2017} attained an $\mathcal{O}(1/\sqrt{T})$ convergence rate comparable to exact communication counterparts, but it requires the compressor precision to increase linearly with the iteration number $T$ and assumes the Lipschitz condition.
	Additionally, \cite{Cao_X_Decentralized_2020} and \cite{singh_decentralized_2024} studied constrained distributed optimization problems, obtaining convergence rates of $\mathcal{O}((p^{\frac{1}{4}}n^{\frac{1}{2}}+p^{\frac{7}{4}}n^{-\frac{1}{2}})/\sqrt{T})$ and $\mathcal{O}(p^2/\sqrt{T})$, respectively, where $n$ represents the number of agents.
Nonconvex optimization presents even greater challenges. Centralized and distributed compressed zeroth-order algorithms proposed by \cite{Kaya_CommunicationEfficient_2023} and \cite{xu2024compressed} achieved convergence rates of $\mathcal{O}(\sqrt{p}/\sqrt{T})$ and $\mathcal{O}(p/\sqrt{T})$, respectively. Although the algorithm in \cite{xu_quantized_2024} achieved linear convergence, it required sampling $\mathcal{O}(pn)$ points per iteration to construct the zeroth-order gradient, leading to significant sampling overhead and necessitating a standard uniform quantizer with a saturation range of $\mathcal{O}(\sqrt{pn})$.

From the discussion above, we pose the following key question:
\emph{Can we design a compressed zeroth-order algorithm for distributed nonconvex optimization that achieves a convergence rate comparable to standard zeroth-order methods 
while relaxing requirements on compressor (e.g., allowing for more general compression techniques) and gradient estimation (e.g., reducing the number of sampled points)?}
This paper provides affirmative answers.

{\bf Main Contributions:} 
This paper proposes a Communication \underline{C}ompressed \underline{Z}eroth-Order \underline{S}tochastic \underline{D}istributed (CZSD) algorithm and analyzes its convergence rate. Specifically, our main contributions are as follows:

 { (i)}  
The algorithm adopts a stochastic two-point zeroth-order oracle, requiring only $2n$ function evaluations per iteration, reducing the sampling overhead. To the best of our knowledge, this is the first compressed zeroth-order algorithm for stochastic distributed nonconvex optimization.

 { (ii)}  
The algorithm utilizes a generalized contractive compressor, encompassing a range of both unbiased and biased compressors. Furthermore, each agent transmits only a single variable per iteration, reducing communication overhead.

{ (iii)}  
The  algorithm achieves an $\mathcal{O}(\sqrt{p}/\sqrt{nT})$ convergence rate for finding stationary points. Furthermore, when the global cost function satisfies the Polyak--{\L}ojasiewicz (P--{\L}) condition, the convergence rate improves to $\mathcal{O}(p/(nT))$ for finding the global optimum, as established in Theorem~\ref{zerosg:thm-sg-smT} and Theorem~\ref{zerosg:thm-sg-diminishingt}, respectively. 
To the best of our knowledge, this is the first compressed zeroth-order algorithm achieving linear speedup for stochastic (and deterministic) distributed nonconvex optimization under both the general setting and the P--{\L} condition, attaining performance comparable to exact communication methods.

 { (iv)}  
The algorithm achieves linear convergence in Theorem~\ref{zerosg:thm-random-pd-fixed} when the relaxed growth condition and state-dependent variance condition are strengthened.





 {\bf Notations:} 
 Let $[n]=\{1,\dots,n\}$ for any positive integer $n$. 
$\|\cdot\|$ denotes the vector $2$-norm.
$\mathbb{B}^p$ and $\mathbb{S}^p$ respectively are the unit ball and unit sphere centered at the origin in $\mathbb{R}^p$.
$\mathrm{Unif(\cdot)}$ denotes the uniform distribution over a given set.
For a differentiable function $f$, $\nabla f$ represents its gradient. 
${\bf 1}_n$ is the $n$-dimensional all-one vector, and ${\bf I}_n$ is the identity matrix of dimension $n$. 
$\col(x_1,\dots,x_k)$ denotes the column vector obtained by concatenating vectors $x_i\in\mathbb{R}^{p_i}$, $i\in[k]$. 
Given a vector $[t_1,\dots,t_n]^\top\in\mathbb{R}^n$, define $\diag([t_1,\dots,t_n])$ as the diagonal matrix with diagonal entries $t_i$. 
$A\otimes B$ denotes the Kronecker product of matrices $A$ and $B$. 
$\rho(\cdot)$ stands for the spectral radius of a matrix, and $\rho_2(\cdot)$ denotes its minimum positive eigenvalue. 
For a positive semi-definite matrix $A$, we define $\|x\|_A=\sqrt{\langle x, Ax\rangle}$.

 {\bf Organization:}
The rest of this paper is organized as follows.
Section~\ref{Problem Formulation} formulates the problem.
Section~\ref{Algorithm and Analysis}  introduces the CZSD algorithm and presents preliminary convergence analysis.
Section~\ref{Main results} states  the main results.
Section~\ref{nonconvex:sec-simulation} provides a numerical example. 
Appendix contains all detailed proofs.


\section{Problem Formulation}\label{Problem Formulation}

Consider $n$ agents, each equipped with 
 a local data set $\Xi_i$ following distribution $\mathcal{D}_i$ and a local cost function $f_i: \mathbb{R}^{p}\mapsto \mathbb{R}$ (not necessarily convex). 
 They collaborate to achieve consensus and optimal model by solving the following optimization problem:
\begin{align}\label{zerosg:eqn:xopt}				
 \min_{x\in \mathbb{R}^p} f(x) \triangleq \frac{1}{n}\sum_{i=1}^nf_i(x) \triangleq \frac{1}{n}\sum_{i=1}^n\mathbb{E}_{\xi_i}[F_i(x,\xi_i)],
\end{align}
where $x\in \mathbb{R}^p$ is the model parameter, $\xi_i$ denotes the local
data that follows the local distribution $\mathcal{D}_i$, and $F_i(x,\xi_i): \mathbb{R}^{p}\times \Xi_i \mapsto \mathbb{R}$ is a stochastic local cost function. 


To solve the problem \eqref{zerosg:eqn:xopt}, there are two key aspects: exchanging information with other agents and collecting information about the local cost function. 
We adopt the following mechanisms for these two aspects, along with additional assumptions.

\subsection{Communication Network and Compression}

This section discusses the use of compressed communication methods for information exchange between agents over a communication network.

To achieve consenus and optimal model, agents must communicate during the optimization process. 
We model the communication network as an undirected and connected graph \( \mathcal{G} = (\mathcal{V}, \mathcal{E}, A) \), where \( \mathcal{V} = [n] \) is the set of agents, \( \mathcal{E} \subseteq \mathcal{V} \times \mathcal{V} \) is the edge set, and \( A =[A]_{ij}\) is the weighted adjacency matrix. \( [A]_{ij} > 0 \) if \( (i,j) \in \mathcal{E} \), and \( [A]_{ij} = 0 \) otherwise. The neighbor set of agent \( i \) is denoted by \( \mathcal{N}_i = \{j \in \mathcal{V} : (i,j) \in \mathcal{E}\} \). The degree matrix \( D \) is defined as $D = \text{diag}(d_1, d_2, \dots, d_n)$, where \( d_i \) denotes the \( i \)-th row sum of \( A \). The Laplacian matrix of the graph is \( L = D - A \).


To improve communication efficiency, we consider the scenario that the communication between agents is compressed. Specifically, we consider a class of compressors with bounded relative compression error satisfying the following assumption \cite{Yi_CommunicationCompression_2023}.

\begin{assumption}[Generalized contractive 	compressor]\label{nonconvex:ass:compression}
	The compressor $\mathcal{C}:\mathbb{R}^p\mapsto\mathbb{R}^p$ satisfies
	\begin{align}\label{nonconvex:ass:compression_equ_scaling}
		\mathbb{E}_{\mathcal{C}}\Big[\Big\|\frac{\mathcal{C}(x)}{r}-x\Big\|^2\Big]\le (1-\delta)\|x\|^2,~\forall x\in\mathbb{R}^p,
	\end{align}
	for some constants $\delta\in(0,1]$ and $r>0$. Here $\mathbb{E}_{\mathcal{C}}[\cdot]$ denotes the expectation over the internal randomness of the stochastic compression operator $\mathcal{C}$.
\end{assumption}
This futher implies 
\begin{align}\label{nonconvex:ass:compression_equ}
	\mathbb{E}_{\mathcal{C}}[\|\mathcal{C}(x)-x\|^2]
		&\le \delta_0\|x\|^2,~\forall x\in\mathbb{R}^p,
\end{align}
where $\delta_0=2r^2(1-\delta)+2(1-r)^2$.


\subsection{Stochastic Zeroth-Order Gradient Estimator}

	 This section concerns the use of the stochastic zeroth-order (ZO) information to estimate the gradients of local cost functions.
In the scenario where explicit gradient expressions are unavailable, we assume each agent can query the
local function value at finite points. Based on two-point sampling, we introduce the stochastic ZO gradient estimator.

Let $f(x)=\mathbb{E}_\xi[F(x, \xi)]$, where $F(x, \xi): \mathbb{R}^p \times \Xi \mapsto \mathbb{R}$ is a  stochastic differentiable function. 
The authors of \cite{Yi_Zerothorder_2022} proposed the following random gradient estimator, which estimates the gradient of $f(x)$ by sampling the function values of the stochastic function $F(x, \xi)$:
\begin{align}
	&g^z (x,\mu,\xi,\zeta)=\frac{p(F(x+ \mu \zeta,\xi)-F(x,\xi))}{ \mu}\zeta, \label{dbco:gradient:model2-st:1}
\end{align}
where $\mu>0$ is an exploration parameter and $\zeta\in\mathbb{S}^p$ is a uniformly distributed random vector indicating the distance and direction of the two sampling points, respectively. Then, $\mathbb{E}_{\zeta\in\mathbb{S}^p}[\mathbb{E}_{\xi}[ g^z(x,\mu,\xi,\zeta) ]]=\nabla \hat{f}(x,\mu)$, where $\hat{f}(x,\mu)\triangleq\mathbb{E}_{\zeta\in\mathbb{B}^p}[f(x+\mu \zeta)]$. Namely, $g^z(x,\mu,\xi,\zeta)$ is a unbiased estimation of $\nabla \hat{f}(x,\mu)$ instead of $\nabla f(x)$. 
Furthermore, if $F(\cdot,\xi)$ is $\ell$-smooth, the gradient estimator $g^z(x,\mu,\xi,\zeta)$ has a non-zero variance as follows:
\begin{align}
	&\mathbb{E}_{\zeta\in\mathbb{S}^p}[\|g^z(x,\mu,\xi,\zeta)\|^2]\le 2p\|\nabla_x F(x,\xi)\|^2+\frac{1}{2}p^2\mu^{2} \ell^2.\label{zerosg:lemma:uniformsmoothing-equ6-1}
	\end{align}

	

\subsection{Other Assumptions}
The following standard assumptions for the problem \eqref{zerosg:eqn:xopt} are made throughout the paper.





\begin{assumption}[Smoothness]\label{zerosg:ass:zeroth-smooth}
	For almost all $\xi_i$, the stochastic local cost function $F_i(\cdot,\xi_i)$ is $\ell$-smooth, i.e.,  there exists a constant $\ell>0$ such that for any $i\in[n]$ and $x,y\in\mathbb{R}^p$,
	\begin{align}\label{nonconvex:smooth}
		\|\nabla F_i(x,\xi_i)-\nabla F_i(y,\xi_i)\|\le \ell\|x-y\|.
	\end{align}
\end{assumption}

\begin{assumption}[Relaxed growth condition]\label{zerosg:ass:zeroth-variance}
	The variance of the stochastic gradient grows proportionally with the gradient norm, i.e., there exist constants $\eta_1, \sigma_1 \ge 0$ such that for any $i \in [n]$ and $x \in \mathbb{R}^p$,
	\begin{align}
		&\mathbb{E}_{\xi_i}[\|\nabla_xF_i(x,\xi_i)-\nabla f_i(x)\|^2]\le\eta^2_1\|\nabla f_i(x)\|^2+\sigma^2_1.
	\end{align}
\end{assumption}
\begin{remark}
	It should be highlighted that Assumption~\ref{zerosg:ass:zeroth-variance} is more general than the bounded variance assumption ($\eta_1 = 0$), which is widely  used in stochastic optimization \cite{	pmlr-v119-lu20a,pmlr-v97-yu19d,Hua_CZOO_2025 }.
\end{remark}

\begin{assumption}[State-dependent variance]\label{zerosg:ass:fig}
	Each local gradient $\nabla f_i(x)$ has state-dependent variance, i.e.,
	there exist constants $\eta_2, \sigma_2\ge0$ such that for any $i\in[n],~x\in\mathbb{R}^p$,
	\begin{align}
		\|\nabla f_i(x)-\nabla f(x)\|^2\le\eta_2^2\|\nabla f(x)\|^2+\sigma^2_2.
	\end{align}
\end{assumption}
\begin{remark}
	
	It should be highlighted that Assumption~\ref{zerosg:ass:fig} is more general than both the bounded variance assumption and the Lipschitz condition.
	It reduces to the bounded variance case when $\eta_2 = 0$, 
	while the Lipschitz condition 
	is much more restrictive.
	 Both assumptions are commonly used in ZO optimization, such as in \cite{ Cao_X_Decentralized_2020, singh_decentralized_2024, ding_quantized_gradient_2017, ding_distributed_consensus_2017}.
\end{remark}
\begin{assumption}[Polyak--{\L}ojasiewicz condition]\label{nonconvex:ass:fil} 
	The global cost function $f(x)$ satisfies the Polyak--{\L}ojasiewicz (P--{\L}) condition, i.e., there exists a constant $\nu>0$ such that for any $x\in\mathbb{R}^p$,
\begin{align}
 \frac{1}{2}\|\nabla f(x)\|^2\ge \nu( f(x)-f^*), 
\label{nonconvex:equ:plc}
\end{align}
where $f^*\triangleq\inf_{x\in\mathbb{R}^p}f(x)>-\infty$.
\end{assumption}

\begin{remark}
	It should be highlighted that no convexity assumptions are made, in contrast to the convex conditions adopted in \cite{Cao_X_Decentralized_2020, singh_decentralized_2024, ding_quantized_gradient_2017, ding_distributed_consensus_2017,Hua_CZOO_2025}.
	Instead, we utilize the P--{\L} condition, a nonconvex analog of strong convexity that ensures global optimality at stationary points and accelerates convergence,
	 applicable broadly to nonconvex problems such as neural network training \cite{liu2022loss}.
\end{remark}

\section{Algorithm and Analysis} \label{Algorithm and Analysis}

	In this section, we propose a Communication \underline{C}ompressed \underline{Z}eroth-Order \underline{S}tochastic \underline{D}istributed (CZSD) algorithm, and analyze its convergence properties based on Lyapunov stability theory.
	First, we introduce the CZSD algorithm and its motivation in \ref{Algorithm Description}.
	Then we provide five supporting lemmas in Section~\ref{Supporting Lemma}, each corresponding to a component term to construct a Lyapunov function. 
	Finally, we discuss the property of this Lyapunov function, which in turn establishes the convergence of the CZSD algorithm in Section~\ref{Lyapunov Function}.	

	\subsection{Algorithm Description} \label{Algorithm Description}

	We first introduce a ZO algorithm, then extend it with communication compression technique.

	\cite{Yi_Zerothorder_2022} proposed the following ZO algorithm:
	\begin{subequations}\label{nonconvex:kia-algo-dc-compress}
		\begin{align}
			&x_{i,k+1} = x_{i,k}-\alpha_k\Big(\beta_k\sum\nolimits_{j=1}^{n}L_{ij}x_{j,k}+\gamma_k v_{i,k}+g^z_{i,k}\Big), \label{nonconvex:kia-algo-dc-x-compress}\\
			&v_{i,k+1} =v_{i,k}+ \alpha_k\gamma_k\sum\nolimits_{j=1}^{n}L_{ij}x_{j,k},  \label{nonconvex:kia-algo-dc-v-compress}\\
			&g^z_{i,k}=\frac{p(F_i(x_{i,k}+ \mu _{i,k} \zeta_{i,k},\xi_{i,k})-F_i(x_{i,k},\xi_{i,k}))}{ \mu _{i,k}}\zeta_{i,k}, \label{dbco:gradient:model2-st}
		\end{align}
	\end{subequations}
	 where $x_{i,k}\in\mathbb{R}^p$ represents agent $i$'s estimation of the solution to the problem \eqref{zerosg:eqn:xopt} at the $k$-th iteration,
	 and $v_{i,k}\in\mathbb{R}^p$ is its dual variable; 
 	$\alpha_k$, $\beta_k$ and $\gamma_k$ are positive algorithm parameters at the $k$-th iteration;
	 $\mu_{i,k}$ is the exploration parameter, and $\zeta_{i,k} \in \mathbb{S}^p$ is a uniformly distributed random vector sampled for agent $i$ at the $k$-th iteration.

To reduce communication cost, we introduce a compressed variable $\hat{x}_{i,k}$ to replace ${x}_{i,k}$, and employ two auxiliary variables $y_{i,k}, z_{i,k} \in \mathbb{R}^p$. Specifically, we define   
\begin{align}
	\hat{x}_{i,k}=y_{i,k}+\mathcal{C}(x_{i,k}-y_{i,k}). \label{nonconvex:kia-algo-dc-compact-xhat}
\end{align}
Here, $y_{i,k}$ helps  reduce compression error. From \eqref{nonconvex:ass:compression_equ} and \eqref{nonconvex:kia-algo-dc-compact-xhat}, we obtain:
\begin{align}\label{nonconvex:xminusxhat}
\mathbb{E}_{\mathcal{C}}[\|x_{i,k} - \hat{x}_{i,k}\|^2] 
&= \mathbb{E}_{\mathcal{C}}[\|x_{i,k} - y_{i,k} - \mathcal{C}(x_{i,k} - y_{i,k})\|^2] \nonumber \\
&\le \delta_0 \mathbb{E}_{\mathcal{C}}[\|x_{i,k} - y_{i,k}\|^2].
\end{align}
Thus, the compression error $\|x_{i,k} - \hat{x}_{i,k}\|^2 \to 0$  as $\|x_{i,k} - y_{i,k}\|^2 \to 0$.
Another auxiliary variable $z_{i,k}$ is used to compute $\sum_{j=1}^{n} L_{ij} y_{i,k}$, allowing for the transmission of only $\mathcal{C}(x_{i,k} - y_{i,k})$, thereby reducing the overall communication cost.
Without this approach, transmitting $\hat{x}_{i,k}$ would increase the communication cost, necessitating the transmission of both $\mathcal{C}(x_{i,k} - y_{i,k})$ and $y_{i,k}$.



	\begin{algorithm}[!tb]
		\caption{Communication \underline{C}ompressed \underline{Z}eroth-Order \underline{S}tochastic \underline{D}istributed (CZSD) Algorithm}
		\label{nonconvex:algorithm-pdgd}
		\begin{algorithmic}[1]
			\STATE \textbf{Input}: positive sequences $\{\alpha_k\}$, $\{\beta_k\}$, $\{\gamma_k\}$ and $\{ \mu _{i,k}\}$; positive parameter $\omega$.
			\STATE \textbf{Initialize}: $ x_{i,0}\in\mathbb{R}^p$, $v_{i,0}=y_{i,0}=z_{i,0}={\bf 0}_p$, and $q_{i,0}=\mathcal{C}(x_{i,0}),~\forall i\in[n]$.
			\FOR{$k=0,1,\dots$}
			\FOR{$i=1,\dots,n$ \textbf{in parallel}} 
			\STATE  \textbf{Communication}: Send $q_{i,k}$ to $\mathcal{N}_i$ and receive $q_{j,k}$ from $j\in\mathcal{N}_i$.
			\STATE \textbf{Stochastic zeroth-order gradient}:\\
			(i) Sample $\xi_{i,k} \sim \mathcal{D}_i$;\\
			(ii) Sample $\zeta_{i,k} \sim \mathrm{Unif}(\mathbb{S}^p)$;\\
			(iii) Sample $F_i(x_{i,k},\xi_{i,k})$, $F_i(x_{i,k}+\mu _{i,k}\zeta_{i,k},\xi_{i,k})$;\\
			 (iv) Compute $g^z_{i,k}$ using \eqref{dbco:gradient:model2-st}.
			\STATE \textbf{Update auxiliary variables}:
			\begin{subequations}\label{zerosg:algorithm-random-pd}
				\vspace{-0.45em} 
				\begin{align}
					y_{i,k+1}&=y_{i,k}+\omega q_{i,k}, \label{nonconvex:kia-algo-dc-a}\\
					z_{i,k+1}&=z_{i,k}+\omega\sum\nolimits_{j=1}^{n}L_{ij}q_{j,k}. \label{nonconvex:kia-algo-dc-b}
				\end{align}\vspace{-1.15em} 
			\STATE \textbf{Update primal and dual variables}:
			\vspace{-0.45em} 				
			\begin{align}
					x_{i,k+1} &= x_{i,k}-\alpha_k\beta_k\Big(z_{i,k}+\sum\nolimits_{j=1}^{n}L_{ij}q_{j,k}\Big)\nonumber\\
					&\qquad-\alpha_k(\gamma_k v_{i,k}+g^z_{i,k} 	), \label{nonconvex:kia-algo-dc-x}\\
					v_{i,k+1} &=v_{i,k}+ \alpha_k\gamma_k\Big(z_{i,k}+\sum\nolimits_{j=1}^{n}L_{ij}q_{j,k}\Big).  \label{nonconvex:kia-algo-dc-v}
				\end{align}\vspace{-1.15em} 
			\STATE \textbf{Compression}:				
			\vspace{-0.5em} 
				\begin{align}
					q_{i,k+1}&=\mathcal{C}(x_{i,k+1}-y_{i,k+1}). \label{nonconvex:kia-algo-dc-q}
				\end{align}\vspace{-1.15em} 
			\end{subequations}
			\ENDFOR
			\ENDFOR
			\STATE  \textbf{Output}: $\{x_{i,k}\}$.
		\end{algorithmic}
	\end{algorithm}

	Based on \eqref{nonconvex:kia-algo-dc-compress} and \eqref{nonconvex:kia-algo-dc-compact-xhat}, we propose the CZSD algorithm (Algorithm~\ref{nonconvex:algorithm-pdgd}).
	 While Algorithm~\ref{nonconvex:algorithm-pdgd} integrates the ZO primal--dual algorithm \cite{Yi_Zerothorder_2022} and the generalized contractive compressor \cite{Yi_CommunicationCompression_2023},
	  its convergence analysis  is mainly challenged by the coupling of gradient and compression errors and the use of the more general compressor.
	Additional difficulties arise from relaxed Assumptions~\ref{zerosg:ass:zeroth-variance},~\ref{zerosg:ass:fig} and nonconvexity. 
	  We address these challenges by designing a Lyapunov function.

	In the subsequent analysis section, we introduce the following notations: $\bsx_k=\col(x_{1,k}, \dots,x_{n,k})$, $\bar{x}_k=\frac{1}{n}({\bm 1}_n^\top\otimes{\bf I}_p)\bsx_k$, $\bar{\bsx}_k={\bm 1}_n\otimes\bar{x}_k$,   
	$\tilde{f}(\bsx_k)=\sum_{i=1}^{n}f_i(x_{i,k})$, $\bsg_k=\nabla\tilde{f}(\bsx_k)$,  $\bsg^0_k=\nabla\tilde{f}(\bar{\bsx}_k)$, $\bar{\bsg}_k^0=\bsH\bsg^0_{k}={\bm 1}_n\otimes\nabla f(\bar{x}_k)$, 
	$\mu_{k}=\max_{i\in[n]}\{\mu_{i,k}\}$, $\bsg_k^z=\col(g^z_{1,k},\dots,g^z_{n,k})$, $g^\mu_{i,k}=\nabla \hat{f}_{i}(x_{i,k},\mu_{i,k})$, $\bsg^\mu_k=\col(g^\mu_{1,k},\dots,g^\mu_{n,k})$, $\bar{\bsg}^\mu_k=\bsH\bsg^\mu_k$,
	$\bsv_k=\col(v_{1,k},\dots,v_{n,k})$, $\bsy_k=\col(y_{1,k},\dots,y_{n,k})$, $\hat{\bsx}_k=\col(\hat{x}_{1,k}, \dots,\hat{x}_{n,k})$, $\bsH=\frac{1}{n}{\bm 1}_n{\bm 1}_n^\top\otimes{\bf I}_p$, and  $\bsL=L\otimes {\bf I}_p$. 
	$\mathcal{B}_k$ represents the $\sigma$-algebra generated by the independent random variables 
	$\xi_{1,k}, \dots, \xi_{n,k}, \zeta_{1,k}, \dots, \zeta_{n,k}$. 
	$\mathcal{C}_k$ represents the $\sigma$-algebra generated by the randomness of the compressor in the $k$-th iteration. 
Let $\mathcal{A}_k = \mathcal{B}_k \cup \mathcal{C}_k$, and define 
$\mathcal{F}_k = \bigcup_{t=0}^k \mathcal{A}_t$.
	Due to the independence of $\xi_{i,k}$, $\zeta_{i,k}$, and the randomness of the compressor, 
	it follows that $\mathcal{B}_k$ and $\mathcal{C}_k$ are independent, $x_{i,k}$ and $v_{i,k}$ for $i \in [n]$ depend on $\mathcal{F}_{k-1}$ and are independent of $\mathcal{A}_t$ for all $t \geq k$.


	
\subsection{Component Lemmas} \label{Supporting Lemma}

In this section, we introduce five component terms $e_{1,k}\sim e_{5,k}$ that constitute the Lyapunov function
:
\begin{align*}
	&{\rm (Consensus ~error)}	&e_{1,k}&\triangleq\frac{1}{2}\|\bsx_{k} \|^2_{\bsE},\\
	&{\rm (Dual ~term)}		&e_{2,k}&\triangleq\frac{1}{2}\Big\|\bsv_{k}+\frac{1}{\gamma_{k}}\bsg_{k}^0\Big\|^2_{\frac{\beta_k+\gamma_k}{\gamma_k}\bsF},\\
	&{\rm (Cross ~term)}		&e_{3,k}&\triangleq\bsx_{k}^\top\bsE\bsF\Big(\bsv_{k}+\frac{1}{\gamma_{k}}\bsg_{k}^0\Big),\\
	&{\rm (Optimality )}		&e_{4,k}&\triangleq\tilde{f}(\bar{\bsx}_{k})-nf^*,\\
	&{\rm (Compression ~error)}&e_{5,k}&\triangleq\|\bsx_{k}-\bsy_{k}\|^2,
\end{align*}
where $\bsE=E\otimes {\bf I}_p$, $\bsF=F_M\otimes {\bf I}_p$, $E={\bf I}_n-\frac{1}{n}{\bm 1}_n{\bm 1}^{\top}_n$, 	
	\begin{align*}
		&F_M=\left[\begin{array}{ll}q&Q
		\end{array}\right]\left[\begin{array}{ll}\lambda_{n+1}^{-1}&0\\
			0&\Lambda_1^{-1}
		\end{array}\right]\left[\begin{array}{l}q^\top\\
			Q^\top
		\end{array}\right],
	\end{align*}
$\Lambda_1=\diag([\lambda_2,\dots,\lambda_n])$ and the orthogonal matrix $[q \ Q]\in \mathbb{R}^{n \times n}$ come form the orthogonal decomposition of the Laplacian matrix $L$ with
	$0<\lambda_2\leq\dots\leq\lambda_n$ being the nonzero eigenvalues of $L$, $q=\frac{1}{\sqrt{n}}\mathbf{1}_n$, $Q \in \mathbb{R}^{n\times (n-1)}$, and $\lambda_{n+1}$ is an arbitrary constant in $[\lambda_2,~\lambda_n]$.
	\begin{remark}
	$e_{1,k}=\frac{1}{2}\|\bsx_{k} \|^2_{\bsE}=\frac{1}{2}\sum\nolimits_{i=1}^{n}\|x_{i,k}-\bar{x}_k\|^2$ represents the consensus error among agents.
	In standard distributed primal--dual algorithms,  
	 the convergence criteria are $\bsv_k \to \frac{1}{\gamma_k}\bsg_k$ and $\bsx_k \to \bar{\bsx}_k$ \cite{Jakovetic2019Unification}, making  $e_{2,k}$ a measure of dual variable convergence.
	$e_{3,k}$ acts as the cross term between $e_{1,k}$ and $e_{2,k}$, while
	$e_{4,k} = \tilde{f}(\bar{\bsx}_k) - nf^* = n(f(\bar{x}_k) - f^*)$ directly reflects the optimality.
	Finally, $e_{5,k}$, as shown in \eqref{nonconvex:xminusxhat}, measures the compression error $\|\bsx_k - \hat{\bsx}_k\|^2$.

\end{remark}	\vspace{-0.75em}
\begin{remark}
	Our goal is to construct a Lyapunov function that captures both consensus error and optimality. However, a direct formulation fails to ensure convergence. 	To address this, we incorporate additional terms $e_{2,k}, e_{3,k}$ and $e_{5,k}$ to enforce a non-increasing structure, 
	with $e_{5,k}$ specifically being used to handle the compression error.
\end{remark}
\vspace{-0.15em}


The following Lemmas~\ref{Lemma Consensus}--\ref{Lemma Compression} characterize the variation of $e_{i,k}$ for $i = 1, \ldots, 5$.  
These derivations involve technical challenges distinct from prior work, primarily due to the coexistence of the compressed variable $\hat{x}_{i,k}$ and the stochastic zeroth-order gradient $g^z_{i,k}$ in all terms.  
To overcome this, we develop a key technique that separates the randomness of compression and gradient estimation at each iteration by introducing the $\sigma$-algebras $\mathcal{C}_k$ and $\mathcal{B}_k$, enabling their decoupling in the analysis.

\begin{lemma}[Consensus error]\label{Lemma Consensus}
		Suppose Assumptions~\ref{zerosg:ass:zeroth-smooth} holds. 
		Let $\{\bsx_k\}$ be the sequence generated by Algorithm~\ref{nonconvex:algorithm-pdgd}.
		Then,
		\begin{align}
			&\mathbb{E}_{\mathcal{A}_k}[e_{1,k+1}]\le e_{1,k}-\|\bsx_k\|^2_{\frac{\alpha_k\beta_k}{2}\bsL-\frac{\alpha_k}{2}\bsE- \alpha_k(1+5\alpha_k) \ell ^2\bsE}
			\nonumber\\
			&\quad
					+ \mathbb{E}_{\mathcal{C}_k}\big[\|\hat{\bsx}_k\|^2_{\frac{3}{2}\alpha_k^2\beta_k^2\bsL^2}\big]
				+ n \ell ^2\alpha_k(1+5\alpha_k)\mu^2_k \nonumber\\
				&\quad + \frac{\alpha_k}{2}(\beta_k+2\gamma_k)\rho(L)\mathbb{E}_{\mathcal{C}_k}[\|\bsx_k-\hat{\bsx}_k\|^2] \nonumber\\
				&\quad -\alpha_k\gamma_k\mathbb{E}_{\mathcal{C}_k}[\hat{\bsx}^\top_k]\bsE\Big(\bsv_k+\frac{1}{\gamma_k}\bsg_k^0\Big)  + 2\alpha_k^2\mathbb{E}_{\mathcal{B}_k}[\|\bsg_k^z\|^2]\nonumber\\
				&\quad +\Big\|\bsv_k+\frac{1}{\gamma_k}\bsg_k^0\Big\|^2_{\frac{6\alpha_k^2\gamma_k^2\rho(L)+\alpha_k\gamma_k}{4}\bsF}.\label{nonconvex:v1k}
		\end{align}
	\end{lemma}
	\begin{proof}
		See Appendix~\ref{zerosg:proof-thm-random-pd-sm}  for the proof.
		\end{proof}


	\begin{lemma}[Dual term]\label{Lemma Optimality 1}
		Suppose Assumptions~\ref{zerosg:ass:zeroth-smooth} holds, $\{\gamma_k\}$ is non-decreasing, and $\beta_k/\gamma_k=\epsilon_1$. 
		Let $\{\bsx_k\}$ be the sequence generated by Algorithm~\ref{nonconvex:algorithm-pdgd}.
		Then,
		\begin{align}
			&\mathbb{E}_{\mathcal{A}_k}[ e_{2,k+1}]\le e_{2,k}  
					+\frac{1}{2}\varepsilon_1(\Delta_k+\Delta_k^2)\mathbb{E}_{\mathcal{A}_k}[\|\bsg_{k+1}^0\|^2]
			\nonumber\\
			&\quad
			 +(1+\Delta_k)\alpha_k(\beta_k+\gamma_k)\mathbb{E}_{\mathcal{C}_k}[\hat{\bsx}^\top_k]\bsE\Big(\bsv_k+\frac{1}{\gamma_k}\bsg_k^0\Big)\nonumber\\
			&\quad+\mathbb{E}_{\mathcal{C}_k}\big[\|\hat{\bsx}_k\|^2_{(1+\Delta_k)\frac{1}{2}\alpha_k^2(\beta_k\gamma_k+\gamma_k^2)\bsL + (1+\Delta_k)\frac{1}{2}\alpha_k^2(\beta_k+\gamma_k)^2\bsE}\big] \nonumber\\
			&\quad +\Big\|\bsv_k+\frac{1}{\gamma_k}\bsg_{k}^0\Big\|^2_{\big(\alpha_k\frac{\gamma_k}{4}+\Delta_k(\frac{\beta_k+\gamma_k}{2\gamma_k}+\frac{\alpha_k\gamma_k}{4})\big)\bsF}\nonumber\\
			&\quad +(1+\Delta_k)(\epsilon_4\alpha_k+\epsilon_5\alpha_k^2) \ell ^2\mathbb{E}_{\mathcal{B}_k}[\|\bar{\bsg}_k^z\|^2], 
			\label{zerosg:v2k}
		\end{align}
	where $\Delta_k=\frac{1}{\gamma_{k}}-\frac{1}{\gamma_{k+1}}$, and $\epsilon_4, \epsilon_5, \varepsilon_1$ are given in Appendix~\ref{zerosg:proof-thm-random-pd-sm}.
	\end{lemma}
	\begin{proof}
		See Appendix~\ref{zerosg:proof-thm-random-pd-sm}  for the proof.
		\end{proof}

	\begin{lemma}[Cross term]\label{Lemma Cross}
		Suppose Assumptions~\ref{zerosg:ass:zeroth-smooth} holds and $\{\gamma_k\}$ is non-decreasing. 
		Let $\{\bsx_k\}$ be the sequence generated by Algorithm~\ref{nonconvex:algorithm-pdgd}.
		Then,
		\begin{align}
			&\mathbb{E}_{\mathcal{A}_k}[ e_{3,k+1}]
				\le e_{3,k} +\|\bsx_k\|^2_{\alpha_k(\frac{\gamma_k+2}{4}+\frac{1}{2} \ell ^2)\bsE+3\alpha_k^2 \ell ^2\bsE}
			\nonumber\\
			&\quad
			-(1+\Delta_k)\alpha_k\beta_k\mathbb{E}_{\mathcal{C}_k}[\hat{\bsx}_k^\top]\bsE\Big(\bsv_k+\frac{1}{\gamma_k}\bsg_{k}^0\Big)\nonumber\\
			&\quad +\mathbb{E}_{\mathcal{C}_k}\Big[\|\hat{\bsx}_k\|^2_{\alpha_k\gamma_k\bsE
			 +\alpha_k^2\big((\frac{1}{2}\beta_k^2+\gamma_k^2)\bsE-\beta_k\gamma_k\bsL\big)+\frac{1}{2}\Delta_k\alpha_k\beta_k\bsE}\Big] \nonumber\\
			&\quad -\Big\|\bsv_k+\frac{1}{\gamma_k}\bsg_{k}^0\Big\|^2_{\big(\alpha_k(\frac{3}{4}\gamma_k-\rho_2^{-1}(L))-\alpha_k^2\gamma_k^2\rho(L)-\frac{1}{2}\rho(L)\Delta_k\alpha_k\beta_k\big)\bsF}\nonumber\\
			&\quad +(\alpha_k\epsilon_6+\alpha_k^2\epsilon_7) \ell ^2\mathbb{E}_{\mathcal{B}_k}\big[\|\bar{\bsg}_k^z\|^2\big]  +n \ell ^2\alpha_k(\frac{1}{2}+3\alpha_k)\mu_k^2 \nonumber\\
			&\quad + \alpha^2_k\mathbb{E}_{\mathcal{B}_k}[\|\bsg_k^z\|^2] +\frac{\alpha_k}{\gamma_k}\rho_2^{-1}(L)\|\bar{\bsg}_k^\mu\|^2 \nonumber\\
			&\quad+\frac{1}{2}\Delta_k\mathbb{E}_{\mathcal{A}_k}[2e_{1,k+1}+\rho_2^{-2}(L)\|\bsg_{k+1}^0\|^2],
		\label{zerosg:v3k}
			\end{align}
			where $\epsilon_6$ and $\epsilon_7$ are given in Appendix~\ref{zerosg:proof-thm-random-pd-sm}.
	\end{lemma}
	\begin{proof}
		See Appendix~\ref{zerosg:proof-thm-random-pd-sm}  for the proof.
		\end{proof}

	\begin{lemma}[Optimality]\label{Lemma Optimality 2}
		Suppose Assumptions~\ref{zerosg:ass:zeroth-smooth} holds and $\{\gamma_k\}$ is non-decreasing. 
		Let $\{\bsx_k\}$ be the sequence generated by Algorithm~\ref{nonconvex:algorithm-pdgd}.
		Then,
		\begin{align}
			&\mathbb{E}_{\mathcal{A}_k}[e_{4,k+1}]
			\le e_{4,k} - \frac{1}{4}\alpha_k\|\bar{\bsg}^\mu_{k}\|^2 + \|\bsx_k\|^2_{\alpha_k \ell ^2\bsE}\nonumber\\
			&\quad + n \ell ^2\alpha_k\mu^2_k-\frac{1}{4}\alpha_k\|\bar{\bsg}_{k}^0\|^2
			+\frac{1}{2}\alpha^2_k \ell \mathbb{E}_{\mathcal{B}_k}[\|\bar{\bsg}^z_{k}\|^2].\label{zerosg:v4k}
		\end{align}
	\end{lemma}
	\begin{proof}
		See Appendix~\ref{zerosg:proof-thm-random-pd-sm}  for the proof.
		\end{proof}

	\begin{lemma}[Compression error] \label{Lemma Compression}
		Suppose Assumptions~\ref{zerosg:ass:zeroth-smooth} holds, $\{\gamma_k\}$ is non-decreasing, and $\omega\le1/r$.
		Let $\{\bsx_k\}$ be the sequence generated by Algorithm~\ref{nonconvex:algorithm-pdgd}.
		Then,
		\begin{align}
			&\mathbb{E}_{\mathcal{A}_k}[e_{5,k+1}]
				\le e_{5,k}-c_2\|\bsx_{k}-\bsy_{k}\|^2
			\nonumber\\
			&\quad +4(1+c_1^{-1})\rho^2(L)\alpha_k^2\beta_k^2 \mathbb{E}_{\mathcal{C}_k}[\|\bsx_k-\hat{\bsx}_k\|_\bsE^2] \nonumber\\
			&\quad + \|\bsx_k\|^2_{4(1+c_1^{-1})\alpha^2_k(\beta_k^2\rho^2(L)+4 \ell ^2)\bsE}\nonumber\\
			&\quad + 16(1+c_1^{-1})n \ell ^2\alpha_k^2\mu_k^2 + 8(1+c_1^{-1})\alpha_k^2\mathbb{E}_{\mathcal{B}_k}[\|\bsg_k^z\|^2]\nonumber\\
			&\quad + \Big\|\bsv_k+\frac{1}{\gamma_k}\bsg_k^0\Big\|^2_{4(1+c_1^{-1})\rho(L)\alpha_k^2\gamma_k^2\bsF},
			\label{nonconvex:xminush}
		\end{align}
		where $c_1$ and $c_2$ are given in Appendix~\ref{zerosg:proof-thm-random-pd-sm}.
	\end{lemma}
	\begin{proof}
		See Appendix~\ref{zerosg:proof-thm-random-pd-sm}  for the proof.
		\end{proof}

\subsection{Lyapunov Function} \label{Lyapunov Function}
In this section, we propose a Lyapunov function $\mathcal{L}_{1,k}\triangleq\sum\nolimits_{i=1}^{5}e_{i,k}$
by directly combining the five error terms and then analyze its variation, which demonstrates the convergence of the CZSD algorithm.

\begin{lemma}[Lyapunov function]\label{zerosg:lemma:sg:L1}
	Suppose  Assumptions~\ref{nonconvex:ass:compression}--\ref{zerosg:ass:fig} hold, $\{\gamma_k\}$ is non-decreasing, $\beta_k/\gamma_k=\epsilon_1$, $\alpha_k\gamma_k=\epsilon_2$, $\epsilon_1>\kappa_1$, $\epsilon_2>0$, and $\gamma_k\ge\varepsilon_{0}$.
	Let $\{\bsx_k\}$ be the sequence generated by Algorithm~\ref{nonconvex:algorithm-pdgd}.
	Then,
	\begin{align}
	&\mathbb{E}_{\mathcal{A}_k}[\mathcal{L}_{1,k+1}]  \nonumber\\
	&~~ \le \mathcal{L}_{1,k} -\|\bsx_k\|^2_{(2 a_1-\varepsilon_5\Delta_k-b_{1,k})\bsE}
	-\Big\|\bsv_k+\frac{1}{\gamma_k}\bsg_{k}^0\Big\|^2_{b_{2,k}\bsF}\nonumber\\
	&~~ - \alpha_k\Big(\frac{1}{4}-(1+\eta_2^2)(b_{3,k}+8p(1+\eta_1^2)b_{4,k})\alpha_k\Big)\|\bar{\bsg}^0_{k}\|^2 \nonumber\\
	&~~ +2pn\sigma^2_1b_{4,k}\alpha_k^2 + n\sigma^2_2(b_{3,k}+4p(1+\eta_1^2)b_{4,k})\alpha_k^2 \nonumber\\
	&~~ +b_{5,k}\alpha_k\mu_k^2	- (c_2-\delta_0\varepsilon_2-\delta_0b_{6,k}\Delta_k)\|\bsx_{k}-\bsy_{k}\|^2, 
		\label{zerosg:sgproof-vkLya:L1}
	\end{align}
	where $\kappa_1, \varepsilon_0, \varepsilon_2, \varepsilon_5,  a_1, b_{1,k} \sim b_{6,k}$ are given in Appendix~\ref{zerosg:proof-thm-random-pd-sm}.
\end{lemma}
\begin{proof}
	Combining Lemma~\ref{Lemma Consensus}--\ref{Lemma Compression} yields the desired result. 
	For the detailed proof, see Appendix~\ref{zerosg:proof-thm-random-pd-sm}.
\end{proof}


\begin{remark}

		Since $\mathcal{L}_{1,k}$ incorporates both consensus error and optimality terms, its convergence directly implies the convergence of the algorithm. 
By properly selecting the algorithm parameters, we ensure that all coefficients of the squared norm terms in \eqref{zerosg:sgproof-vkLya:L1} are negative, allowing us to establish the non-increasing property of $\mathcal{L}_{1,k}$ via Lemma~\ref{zerosg:lemma:sg:L1}. 
Moreover, it is clear from the definition that $\mathcal{L}_{1,k}$ is non-negative.
Together, these two conditions guarantee the convergence of $\mathcal{L}_{1,k}$, and hence of the algorithm.
Finally, the convergence rate is obtained by summarizing the inequality \eqref{zerosg:sgproof-vkLya:L1} and exploiting the recursive structure.

\end{remark}


\section{Main results}\label{Main results}

In this section, we present the convergence results of Algorithm~\ref{nonconvex:algorithm-pdgd}.
Since the error caused by the stochastic ZO gradient does not tend to zero as shown in \eqref{zerosg:lemma:uniformsmoothing-equ6-1},
fixed (depending on the number of iterations $T$) or time-varying algorithm parameters 
  must be used. 
Under different parameter settings and assumptions, 
  we have the following results.

	If \(T\) is known in advance, linear speedup can be achieved for finding a stationary point of a nonconvex problem. 

	\begin{theorem}[Stationary point]\label{zerosg:thm-sg-smT}
	Suppose Assumptions~\ref{nonconvex:ass:compression}--\ref{zerosg:ass:fig} hold. Consider the sequence $\{\bsx_k\}$ generated by Algorithm~\ref{nonconvex:algorithm-pdgd} under the following setting:
	\begin{align}\label{zerosg:step:eta2-sm}
		&\alpha_k=\frac{\sqrt{n}}{\sqrt{pT}}, ~\beta _k=\epsilon _1\gamma _k, ~\gamma _k=\epsilon_2/\alpha _k, ~\mu _{i,k}\le\frac{\kappa_\mu\sqrt{p\alpha _k}}{\sqrt{n+p}}, \nonumber\\
		& \qquad~\omega\le1/r, ~T\ge\max\{\frac{n(\tilde{\kappa}_0(\epsilon_1,\epsilon_2))^2}{p\epsilon_2^2},~\frac{n^3}{p}\}, \nonumber\\
		&\qquad \quad ~~\epsilon _1>\kappa _1, ~\epsilon _2\in ( 0,\kappa _2( \epsilon _1 ) ), ~\kappa_\mu>0,
	\end{align}	
	where $\kappa_2(\epsilon_1)$ and $\tilde{\kappa}_0(\epsilon_1,\epsilon_2)$ are given in Appendix~\ref{zerosg:proof-thm-sg-smT}.
	Then,
	\begin{subequations}
	\begin{align}
	&~~\frac{1}{T}\sum\nolimits_{k=0}^{T-1}\mathbb{E}[\|\nabla f(\bar{x}_k)\|^2]
	=\mathcal{O}\Big(\frac{\sqrt{p}}{\sqrt{nT}}\Big)+\mathcal{O}\Big(\frac{n}{T}\Big),\label{zerosg:coro-sg-sm-equ3}\\
	&~~\mathbb{E}[f(\bar{x}_{T})]-f^*=\mathcal{O}(1),\label{zerosg:coro-sg-sm-equ4}\\
	&~~\mathbb{E}\Big[\frac{1}{n}\sum\nolimits_{i=1}^{n}\|x_{i,T}-\bar{x}_T\|^2\Big]
	=\mathcal{O}\Big(\frac{n}{T}\Big).\label{zerosg:coro-sg-sm-equ3.1}
	\end{align}
	\end{subequations}
	\end{theorem}
	\begin{proof}
	The detailed proof proof is given in Appendix~\ref{zerosg:proof-thm-sg-smT}.
	\end{proof}

	\begin{remark}
		It should be highlighted that the constants omitted in the first dominant term 	on the right-hand side of \eqref{zerosg:coro-sg-sm-equ3} are independent of any parameters related to the communication network.
		As a result, Algorithm~\ref{nonconvex:algorithm-pdgd} achieves linear speedup with a convergence rate of $\mathcal{O}(\sqrt{p}/\sqrt{nT})$, enabling accelerated convergence as more agents are added.
		To the best of our knowledge, this is the fastest proved rate for finding stationary points in distributed compressed zeroth-order algorithms, outperforming the proved rates in \cite{xu2024compressed,Cao_X_Decentralized_2020,singh_decentralized_2024,ding_quantized_gradient_2017,ding_distributed_consensus_2017,Hua_CZOO_2025}.
The compression constants slow down convergence by affecting consensus and indirectly optimality, but this effect vanishes as $\delta$ and $r$ approach $1$. Details, reflected in the coefficients $a_1$, $\tilde{a}_4\sim\tilde{a}_6$, and $\tilde{a}_5^{\prime}$, are given in the Appendix~\ref{zerosg:proof-thm-sg-smT}.
	\end{remark}
	If the global cost function satisfies the P--{\L} condition, using a time-varying step size leads to faster convergence to the global optimum. 

	\begin{theorem}[Global optimum]\label{zerosg:thm-sg-diminishingt}
Suppose Assumptions~\ref{nonconvex:ass:compression}--\ref{nonconvex:ass:fil} hold. 
Consider the sequence $\{\bsx_k\}$ generated by Algorithm~\ref{nonconvex:algorithm-pdgd} 
under the following setting:
\begin{align}\label{zerosg:step:eta1t1}
	&~~~\gamma_k=\epsilon_3(k+m), ~\beta_k=\epsilon_1\gamma_k,
	~ \alpha_k=\epsilon_2/\gamma_k,\nonumber\\
	&\qquad \quad~~~\mu_{i,k}\le \frac{\kappa_\mu\sqrt{p\alpha_k}}{\sqrt{n+p}}, ~\omega\le1/r, \nonumber\\
	&\epsilon_1>\kappa_1, ~\epsilon_2\in(0,\kappa_2(\epsilon_1)), ~\epsilon_3\in[3\hat{\kappa}_0\nu\epsilon_2/16,3\nu\epsilon_2/16),\nonumber\\
	&\quad ~m\ge\hat{\kappa}_3(\epsilon_1,\epsilon_2,\epsilon_3), ~\hat{\kappa}_0\in(0,1), ~\kappa_\mu>0, 
	\end{align}
where 
$\hat{\kappa}_3(\epsilon_1,\epsilon_2,\epsilon_3)$ is given in Appendix~\ref{zerosg:proof-thm-sg-diminishingt}.
Then,
	\begin{subequations}
	\begin{align}
	&\mathbb{E}\Big[\frac{1}{n}\sum\nolimits_{i=1}^{n}\|x_{i,T}-\bar{x}_T\|^2\Big]
	=\mathcal{O}\Big(\frac{p}{T^{2}}\Big), \label{zerosg:thm-sg-diminishing-equ2.1bounded}\\
	&\mathbb{E}[f(\bar{x}_{T})-f^*]
	=\mathcal{O}\Big(\frac{p}{nT}\Big)+\mathcal{O}\Big(\frac{p}{T^{2}}\Big) \nonumber\\
		&\qquad + \mathcal{O}\Big(\frac{p^{\theta}}{T^{\theta}}\Big)+\mathcal{O}\Big(\frac{p}{T^{3}}\Big),~\forall T\in\mathbb{N}_+,
	\label{zerosg:thm-sg-diminishing-equ2bounded}
	\end{align}
	\end{subequations}
	where $\theta=\frac{3\nu\epsilon_2}{8\epsilon_3}\in(2,\frac{2}{\hat{\kappa}_0}]$.
	\end{theorem}
	\begin{proof}
		The detailed proof proof is given in Appendix~\ref{zerosg:proof-thm-sg-diminishingt}. 
	\end{proof}

\begin{remark}
	It should be highlighted that the constants omitted in the first dominant term on the right-hand side of \eqref{zerosg:thm-sg-diminishing-equ2bounded} are independent of any parameters related to the communication network, thus Algorithm~\ref{nonconvex:algorithm-pdgd} achieves linear speedup with a convergence rate of $\mathcal{O}(p/(nT))$.
	Compared to \cite{xu2024compressed}, our algorithm has a faster proved convergence rate and 
	reduces communication cost by transmitting only one compressed variable instead of two.
The impact of the compression constants is reflected in the coefficients $a_1\sim a_6$ and $a_5^{\prime}$, as detailed in Appendix~\ref{zerosg:proof-thm-sg-diminishingt}.
\end{remark}

Further strengthening Assumptions~\ref{zerosg:ass:zeroth-variance} and~\ref{zerosg:ass:fig} with $\sigma_1=\sigma_2=0$ yields linear convergence,
where $\sigma_1=0$ holds when the deterministic zeroth-order information is  available, and $\sigma_2=0$ holds, for example, when all $\xi_i$ share the same distribution and the local objectives $f_i$ have the same functional form.

\begin{theorem}[Linear convergence]\label{zerosg:thm-random-pd-fixed}
Suppose Assumptions~\ref{nonconvex:ass:compression}--\ref{nonconvex:ass:fil} hold, and $\sigma_1=\sigma_2=0$. Consider the sequence $\{\bsx_k\}$ generated by Algorithm~\ref{nonconvex:algorithm-pdgd} under the following setting:
\begin{align}\label{zerosg:step:eta2-fixed}
&\gamma_k=\gamma, ~\beta _k=\epsilon _1\gamma _k, ~\alpha _k=\epsilon_2/\gamma _k, ~\mu_{i,k}\le\kappa_\mu\tilde{\varepsilon}^{k}, ~\omega\le1/r,\nonumber\\
&\epsilon_1>\kappa_1, \epsilon_2\in(0,\kappa_2(\epsilon_1)), \gamma\ge\tilde{\kappa}_0(\epsilon_1,\epsilon_2),  \tilde{\varepsilon}\in(0,1), \kappa_\mu>0.
\end{align}
Then,  
	\begin{align}\label{zerosg:thm-sg-fixed-equ1-coro1}
	&\mathbb{E}\Big[\frac{1}{n}\sum\nolimits_{i=1}^{n}\|x_{i,k}-\bar{x}_k\|^2+f(\bar{x}_k)-f^*\Big]
	=\mathcal{O}(\varepsilon^{k})
	,
	\end{align}
where $\varepsilon \in (0,1)$ is given in Appendix~\ref{zerosg:proof-thm-random-pd-fixed}.
\end{theorem}
\begin{proof}
	The detailed proof proof is given in Appendix~\ref{zerosg:proof-thm-random-pd-fixed}.
\end{proof}

\begin{remark}
	While \cite{xu_quantized_2024} also achieves linear convergence, their method requires sampling \(\mathcal{O}(pn)\) points  to construct the zeroth-order gradient at each iteration. 	
	This results in $p$ times more samples per iteration compared to our algorithm, leading to greater sampling pressure.
The impact of the compression constants is reflected in the coefficients $\varepsilon_{15}$, as detailed in Appendix~\ref{zerosg:proof-thm-random-pd-fixed}.
\end{remark}

\section{Simulations}\label{nonconvex:sec-simulation}
In this section, we validate theoretical results via numerical experiments on a nonconvex distributed binary classification task \cite{xu2024compressed,Yi_CommunicationCompression_2023}, with comparisons against the state-of-the-art method.
The optimization objective follows \eqref{zerosg:eqn:xopt} with each component function is defined as
\begin{align*}
&F_i(x,\xi_i)=\frac{n}{m_i}\sum_{j=1}^{m_i}\log\Big(1+e^{-t_{ij}x^\top s_{ij}}\Big)+\sum_{l=1}^{p}\frac{\theta\tau[x]_l^2}{1+\tau[x]_l^2},
\end{align*}
where  $m_i=200$ denotes the number of local observations per agent, 
$\xi_i=(s_{ij},t_{ij})$, 
$s_{ij}\in\mathbb{R}^p$ represents Gaussian features regenerated per iteration, labels $t_{ij}\in\{-1,1\}$ are sampled with probability $p(t_{ij}=1)=(1+e^{-t_{ij}x_i^\top s_{ij}})^{-1}$;
Regularization parameters are set to $\theta=0.001$ and $\tau=1$, with $n=20$ agents and dimension $p=50$ following \cite{xu2024compressed,Yi_CommunicationCompression_2023}, 
$[x]_l$ is the $l$-th coordinate of $x\in \mathbb{R}^p$. 
Agent interactions are modeled by a 3D spherical random geometric graph with angular threshold set to $10^\circ$.


We implement an unbiased 2-bit quantizer \cite{Liao_CompGradTrack_2022}  satisfying Assumption~\ref{nonconvex:ass:compression} with $r=1+p/4^k$ and $\delta=1/(1+p/4^k)$:
	\begin{align*}
	\mathcal{C}(x)=\frac{\|x\|_\infty}{2^{k-1}}\text{sgn}(x)\circ
	\bigg\lfloor\frac{2^{k-1}\left|{x}\right|}{\|x\|_\infty}+\varpi\bigg\rfloor,
	\end{align*}
 where $k=2$,
 $\varpi \sim \mathrm{Unif}([0,1]^p)$ is a random vector uniformly distributed in the unit hypercube $[0,1]^p$,
 $\text{sign}(\cdot)$, $\circ$, $\lfloor\cdot\rfloor$, and $\left|{\cdot}\right|$ are the element-wise sign, Hadamard product, floor function, and absolute value, respectively.
 This scheme achieves per-vector communication cost of $(k+1)p + 64$ bits, where the 64-bit floating-point ensures exact precision.

To the best of our knowledge, no existing compressed zeroth-order methods address such stochastic distributed nonconvex optimization.
We, therefore, compare our  Algorithm~\ref{nonconvex:algorithm-pdgd} against 
 zeroth-order stochastic distributed primal--dual algorithm (ZSD-PD). 
Performance metric $P(T)=\min_{k\in [T]}\{\mathbb{E}_{\xi}[\|\nabla_x F(\bar{x}_k,\xi)\|^ {2}]+\frac{1}{n}\sum_{i=1}^{n}\|x_{i,k}-\bar{x}_k\|^ {2}\}
$ is analyzed against both iteration count and inter-agent communication bits.
Hyperparameters used in the experiment are tuned manually and given in TABLE~\ref{tab:para}.
\vspace{-0.4em}
\begin{table}[!ht] 
	\caption{Hyperparameters  for Algorithms}	\label{tab:para} \vspace{-0.7em}
	\centering
	\small
	\begin{tabular}{|M{1.6cm}|M{1cm}|M{1cm}|M{1cm}|M{0.9cm}|M{0.4cm}|}
	\hline
	Algorithm & $\alpha_k$ & $\beta_k$ & $\gamma_k$& $\mu_{ik}$ & $\omega$\\
	\hline
	ZSD-PD &0.1/(k+1) & 3(k+1) & 0.1(k+1) & $0.99^k$ &--\\
	\hline
	 Algorithm~\ref{nonconvex:algorithm-pdgd} &0.1/(k+1) & 3(k+1) & 0.1(k+1) & $0.99^k$ &0.1\\
	\hline
	\end{tabular}
	\end{table}

\begin{figure}[!ht]
	\vspace*{0.3cm}
\centering
  \includegraphics[width=0.48\textwidth]{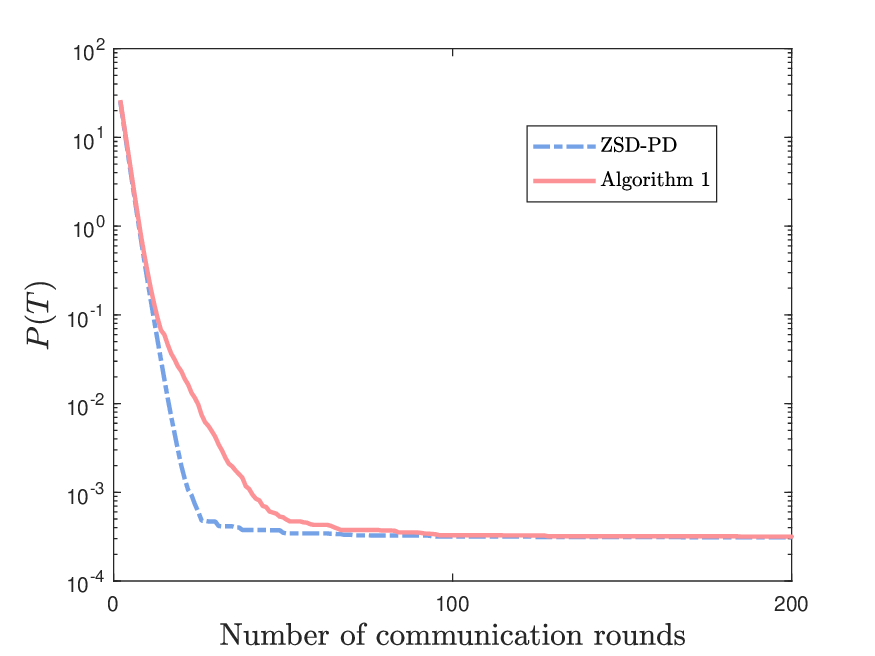}
  \caption{Evolutions of $P(T)$ with respect to the number of iterations.}
  \label{nonconvex:fig:iterations}
\end{figure}

\begin{figure}[!ht]
\centering
  \includegraphics[width=0.48\textwidth]{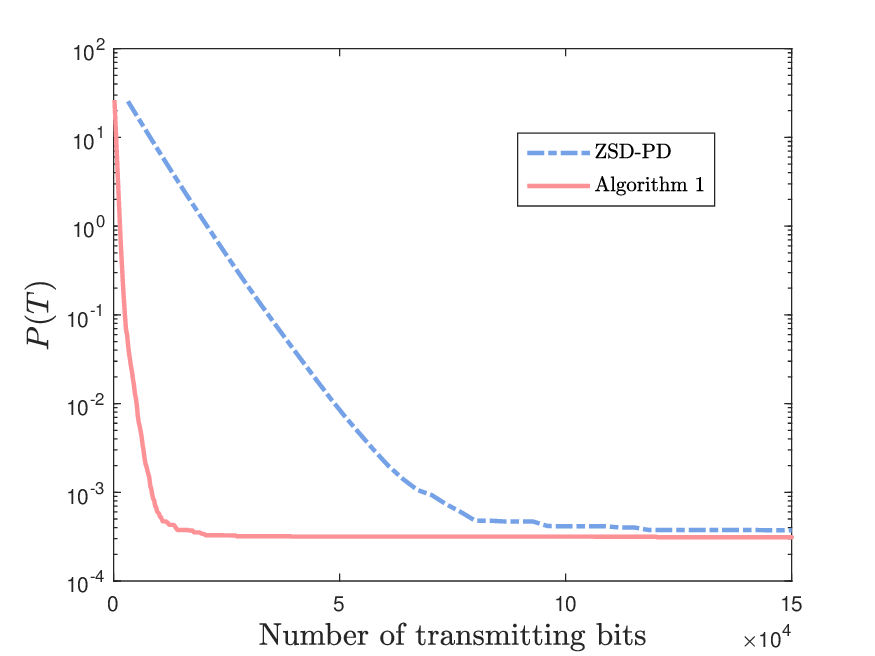}
  \caption{Evolutions of $P(T)$ with respect to the number of transmitted bits.}
  \label{nonconvex:fig:bits}
\end{figure}

Fig.~\ref{nonconvex:fig:iterations} illustrates that Algorithm~\ref{nonconvex:algorithm-pdgd} achieves comparable performance to ZSD-PD.
The slightly slower convergence of Algorithm~\ref{nonconvex:algorithm-pdgd} compared to ZSD-PD (exact transmission) is reasonable, given that our algorithm transmits only one compressed variable per iteration.
Moreover, Fig.~\ref{nonconvex:fig:bits} highlights its communication efficiency: to achieve $P(T)\le10^{-3}$, ZSD-PD requires $7.85\times$ more communication bits than Algorithm~\ref{nonconvex:algorithm-pdgd}.

\section{Conclusions} \label{Conclusions}
In this paper, we proposed the CZSD algorithm to balance the trade-offs between communication efficiency, gradient estimation complexity, and convergence speed.
To the best of our knowledge, we established the first linear speedup for  compressed zeroth-order algorithms in stochastic distributed nonconvex optimization, matching comparable convergence rates to exact communication methods.
Future directions include extending the algorithm to directed graphs, exploring improved zeroth-order methods through variance reduction, and developing more generalized compression schemes.

\section{ACKNOWLEDGMENTS}

The authors would like to thank Dr. Shengjun Zhang and Dr. Lei Xu for sharing their codes.

\appendix

\subsection{Useful Lemmas}\label{zero:app-lemmas}
The following results are used in the proofs.

\subsubsection{Useful Inequality}
	\begin{lemma}\label{nonconvex:lemma:inequality-arithmetic}
		For any $x,y\in\mathbb{R}^d$ and $a,b>0$ satisfying $ab=\frac{1}{4}$, it holds that \begin{align}\label{nonconvex:lemma:inequality-arithmetic-equ}
			x\top y\le a\|x\|^2+b\|y\|^2.
		\end{align}
	\end{lemma}
	This lemma is a direct extension of the Cauchy--Schwarz inequality.

	\begin{lemma}\label{zerosg:lemma:sumgeo}
		Let $a,b\in(0,1)$ be two constants, then
		\begin{align}\label{zerosg:lemma:sumgeo-equ}
		\sum_{\tau=0}^{k}a^\tau b^{k-\tau}\le
		\begin{cases}
		\frac{a^{k+1}}{a-b}, & \mbox{if } a>b, \\
		\frac{b^{k+1}}{b-a}, & \mbox{if } a<b, \\
		\frac{c^{k+1}}{c-b}, & \mbox{if } a=b, 
		\end{cases}
		\end{align}
		where $c$ is any constant in $(a,1)$.
		\end{lemma}

	\begin{lemma} 
		Let $f: \mathbb{R}^p \to \mathbb{R}$ be a differentiable function that is $\ell$-smooth for some constant $\ell > 0$. Then for all $x, y \in \mathbb{R}^p$, the following inequalities hold:
		\begin{align}
			&|f(y) - f(x) - (y - x)^\top \nabla f(x)| \le \frac{\ell}{2} \|y - x\|^2, \label{nonconvex:lemma:lipschitz} \\
			&\frac{1}{2} \|\nabla f(x)\|^2 \le \ell (f(x) - f^*).
		\end{align}
	\end{lemma}

\subsubsection{Properties of useful matrices}
 \begin{lemma}[Lemma 3,\cite{Yi_CommunicationCompression_2023}]
	The matrices $L, E$ are positive semi-definite and $F_M$ is positive definite. Moreover they satisfy the following properties:
	\begin{align}
		&EL=LE=L,\label{nonconvex:KL-L-eq}\\ 
		&0\le\rho_2(L)E\le L\le\rho(L)E,\label{nonconvex:KL-L-eq2}\\
		&F_ML=LF_M=E,\label{nonconvex:lemma-eq3}\\
		&\rho^{-1}(L){\bf I}_n\leq F_M \le\rho_2^{-1}(L){\bf I}_n.\label{nonconvex:lemma-eq5}
	\end{align}
 \end{lemma}

\subsubsection{Properties of Gradient Approximation}
\begin{lemma}[Lemma 6, \cite{Yi_Zerothorder_2022}]\label{zerosg:lemma:grad-st}
	Suppose Assumption~\ref{zerosg:ass:zeroth-smooth} holds. Let $\{\bsx_k\}$ be the sequence generated by Algorithm~\ref{nonconvex:algorithm-pdgd}, then
	\begin{subequations}
	\begin{align}
		\mathbb{E}_{\mathcal{B}_k}[\bsg_k^z]&=\bsg^\mu_k,\label{zerosg:rand-grad-esti1}\\
		\|\bsg_k^0-\bsg^\mu_k\|^2&\le 2 \ell ^2\|\bsx_{k}\|^2_{\bsE}+2n \ell ^2\mu_k^2,\label{zerosg:rand-grad-esti8}\\
		\|\bar{\bsg}_k^0-\bar{\bsg}^\mu_k\|^2
		&\le 2 \ell ^2\|\bsx_{k}\|^2_{\bsE}+2n \ell ^2\mu_k^2,\label{zerosg:rand-grad-esti9}\\
		\mathbb{E}_{\mathcal{B}_k}[\|\bar{\bsg}^z_{k}\|^2]
		&\le  \frac{1}{n}\mathbb{E}_{\mathcal{B}_k}[\|\bsg_k^z\|^2]+\|\bar{\bsg}^\mu_{k}\|^2,\label{zerosg:rand-grad-esti5}\\
		\mathbb{E}_{\mathcal{B}_k}[\|\bsg_k^0-\bsg_k^z\|^2]&\le 4 \ell ^2\|\bsx_{k}\|^2_{\bsE}+4n \ell ^2\mu_k^2
		+2\mathbb{E}_{\mathcal{B}_k}[\|\bsg_k^z\|^2],\label{zerosg:rand-grad-esti6}\\
		\|\bsg^0_{k+1}-\bsg^0_{k}\|^2&\le \alpha^2_k \ell ^2\|\bar{\bsg}^z_{k}\|^2
		\le\alpha^2_k \ell ^2\|\bsg^z_{k}\|^2,\label{zerosg:gg-rand-pd}\\
		\|\bar{\bsg}^0_k\|^2&\le 2n \ell (f(\bar{x}_k)-f^*).\label{zerosg:rand-grad-smooth}
	\end{align}
	\end{subequations}

	If Assumptions~\ref{zerosg:ass:zeroth-variance} and \ref{zerosg:ass:fig} also hold, then
	\begin{subequations}
	\begin{align}
		\mathbb{E}_{\mathcal{B}_k}[\|\bsg_k^z\|^2]
		&\le  16p(1+\eta_1^2)(1+\eta_2^2)(\|\bar{\bsg}_{k}^0\|^2
		+ \ell ^2\|\bsx_{k}\|^2_{\bsE})\nonumber\\
		&+4np\sigma^2_1
		+8np(1+\eta_1^2)\sigma^2_2+\frac{1}{2}np^2 \ell ^2\mu_k^2,\label{zerosg:rand-grad-esti2}\\
		\|\bsg^0_{k+1}\|^2&\le 3(\alpha^2_k \ell ^2\|\bsg^z_{k}\|^2+n\sigma^2_2
		+(1+\eta_2^2)\|\bar{\bsg}_{k}^0\|^2).\label{zerosg:rand-grad-esti4}
	\end{align}
	\end{subequations}
	\end{lemma}



	\subsubsection{Convergence properties of the sequence}
	\begin{lemma}[Lemma 5, \cite{Yi_Zerothorder_2022}]\label{zerosg:serise:lemma:sequence} 
		Let $\{\psi_k\}$, $\{r_{1,k}\}$, and $\{r_{2,k}\}$ be sequences. Suppose there exists $m\in\mathbb{N}_+$ such that
		\begin{subequations}
		\begin{align}
		&\psi_k\ge0,\label{zerosg:serise:lemma:sequence-equ0}\\
		&\psi_{k+1}\le(1-r_{1,k})\psi_k+r_{2,k},\label{zerosg:serise:lemma:sequence-equ1}\\
		&1> r_{1,k}\ge\frac{d_1}{(k+m)^{ \theta_1}},\label{zerosg:serise:lemma:sequence-equ2}\\
		&r_{2,k}\le\frac{d_2}{(k+m)^{ \theta_2}},~\forall k\in\mathbb{N}_0, \label{zerosg:serise:lemma:sequence-equ3}
		\end{align}
		\end{subequations}
		where $d_1>0$, $d_2>0$, $ \theta_1\in[0,1]$, and $ \theta_2> \theta_1$ are constants.

	(i) If $ \theta_1=1$, then
	\begin{align}\label{zerosg:serise:lemma:sequence-equ5}
	\psi_{k}&\le \phi_1(k,m,d_1,d_2, \theta_2,\psi_{0}),~\forall k\in\mathbb{N}_+,
	\end{align}
	where
	\begin{align}\label{zerosg:serise:lemma:sequence-equ5-phi3}
	\phi_1(k,m,d_1,d_2, \theta_2,\psi_{0})&=\frac{m^{d_1}\psi_{0}}{(k+m)^{d_1}}
	+\frac{d_2}{(k+m-1)^{ \theta_2}}\nonumber\\
	&+\Big(\frac{m+1}{m}\Big)^{ \theta_2}d_2s_2(k+m),
	\end{align}
	and
	\begin{align*}
	s_2(k)=
	\begin{cases}
	  \frac{1}{(d_1- \theta_2+1)k^{ \theta_2-1}}, & \mbox{if } d_1- \theta_2>-1, \\
	  \frac{\ln(k-1)}{k^{d_1}}, & \mbox{if } d_1- \theta_2=-1, \\
	  \frac{-m^{d_1- \theta_2+1}}{(d_1- \theta_2+1)k^{d_1}}, & \mbox{if } d_1- \theta_2<-1.
	\end{cases}
	\end{align*}
	
	(ii) If $ \theta_1=0$, then
	\begin{align}\label{zerosg:serise:lemma:sequence-equ6}
	\psi_{k}&\le \phi_2(k,m,d_1,d_2, \theta_2,\psi_{0}),~\forall k\in\mathbb{N}_+,
	\end{align}
	where
	\begin{align}\label{zerosg:serise:lemma:sequence-equ6-phi4}
	&\phi_2(k,m,d_1,d_2, \theta_2,\psi_{0})
	=(1-d_1)^k\psi_{0}+d_2(1-d_1)^{k+m-1}\nonumber\\
	&\qquad\Big([m_2-m]_+s_3(m)+([m_3-m]_+-[m_2-m]_+)s_3(m_3)\Big)\nonumber\\
	&\qquad+\frac{{\bm 1}_{(k+m-1\ge m_3)}2 d_2}{-\ln(1-d_1)(k+m)^{ \theta_2}(1-d_1)},
	\end{align}
	$\quad s_3(k)=\frac{1}{k^{ \theta_2}(1-d_1)^{k}}$, $m_2=\lceil \frac{- \theta_2}{\ln(1-d_1)}\rceil$, and $m_3=\lceil \frac{-2 \theta_2}{\ln(1-d_1)}\rceil$.
	\end{lemma}

	(iii) If $ \theta_1\in(0,1)$, then
	\begin{align}\label{zerosg:serise:lemma:sequence-equ4}
	\psi_{k}\le\phi_3(k,m,d_1,d_2, \theta_1, \theta_2,\psi_{0}),~\forall k\in\mathbb{N}_+,
	\end{align}
	where
	\begin{align}\label{zerosg:serise:lemma:sequence-equ4-phi2}
	&\phi_3(k,m,d_1,d_2, \theta_1, \theta_2,\psi_{0})
	=\frac{1}{s_1(k+m)}\Big(s_1(m)\psi_{0}\\
	&+\frac{[m_4-1-m]_+s_1(m+1)d_2}{m^{ \theta_2}}\Big)\nonumber
	\quad+\frac{d_2}{(k+m-1)^{ \theta_2}}\\
	&+\frac{{\bm 1}_{(k+m-1\ge m_4)}(\frac{m+1}{m})^{ \theta_2}d_2 \theta_2}{d_1 \theta_1(k+m)^{ \theta_2- \theta_1}},
	\end{align}
	$s_1(k)=e^{\frac{d_1}{1- \theta_1}k^{1- \theta_1}}$ and $m_4=\lceil(\frac{ \theta_2}{d_1})^{\frac{1}{1- \theta_1}}\rceil$.

\subsection{Proof of Lemma~\ref{Lemma Consensus}--\ref{zerosg:lemma:sg:L1}}\label{zerosg:proof-thm-random-pd-sm}

	Denote $\mu_{k}=\max_{i\in[n]}\{\mu_{i,k}\}$, $\bsx=\col(x_1,\dots,x_n)$, $\tilde{f}(\bsx)=\sum_{i=1}^{n}f_i(x_i)$,  $\bar{x}_k=\frac{1}{n}({\bm 1}_n^\top\otimes{\bf I}_p)\bsx_k$, $\bar{\bsx}_k={\bm 1}_n\otimes\bar{x}_k$, $\bsg_k=\nabla\tilde{f}(\bsx_k)$, $\bar{\bsg}_k=\bsH\bsg_{k}$, $\bsg^0_k=\nabla\tilde{f}(\bar{\bsx}_k)$, $\bar{\bsg}_k^0=\bsH\bsg^0_{k}={\bm 1}_n\otimes\nabla f(\bar{x}_k)$, $\bsg_k^z=\col(g^z_{1,k},\dots,g^z_{n,k})$, $\bar{g}_k^z=\frac{1}{n}({\bm 1}_n^\top\otimes{\bf I}_p)\bsg_k^z$,  $\bar{\bsg}_k^z={\bm 1}_n\otimes\bar{g}_k^z=\bsH\bsg_k^z$, $\hat{f}_{i}(x,\mu_{i,k})=\mathbb{E}_{\zeta\in\mathbb{B}^p}[f_i(x+\mu_{i,k} \zeta)]$, $g^\mu_{i,k}=\nabla \hat{f}_{i}(x_{i,k},\mu_{i,k})$, $\bsg^\mu_k=\col(g^\mu_{1,k},\dots,g^\mu_{n,k})$, and $\bar{\bsg}^\mu_k=\bsH\bsg^\mu_k$. Moreover, without ambiguity, we denote $\mathcal{C}(\bsx)=\col(\mathcal{C}(x_1),\dots,\mathcal{C}(x_n))$.
$\bsx=\col(x_1,\dots,x_n)$,
	We also denote the following notations:
	\begin{align*}
		&\kappa_1=\max\Big\{\frac{13}{2\rho_2(L)},~\rho_2(L)\Big\},\\
		&\epsilon_4=\frac{(\beta_k+\gamma_k)^2}{\gamma_k^5}\rho_2^{-1}(L),\\
		&\epsilon_5=\frac{\beta_k+\gamma_k}{2\gamma_k^3}\rho_2^{-1}(L)
			+ \frac{1}{2\gamma_k^2},\\	
		&\epsilon_6=\frac{1}{2\gamma_k^2}\rho_2^{-2}(L),\\
		&\epsilon_7=\frac{1}{2\gamma_k^2}+\frac{1}{\gamma_k^2}\rho_2^{-2}(L),\\
		&\varepsilon_{0}=\max\{ 1+\frac{5}{2} \ell ^2, 
		~\big(8+16(1+c_1^{-1})+8p(1+\eta_1^2)(1+\\
			&\quad \eta_2^2)(6+ 16(1+c_1^{-1})+ \ell) + 2 b_{7,k}\big)^\frac{1}{2} \ell ,~12\rho_2^{-1}(L),\\
			&\quad 128p(1+\eta_1^2)(1+\eta_2^2)\epsilon_2\ell\}, \\
		&\varepsilon_1=(1+\epsilon_1)\rho_2^{-1}(L),\\
		&\varepsilon_2=\frac{1}{2}\rho(L)\epsilon_1\epsilon_2+(2+\rho(L))\epsilon_2\\
			& \quad + \big(3\rho^2(L)+4(1+c_1^{-1})\rho^2(L)+2\big)\epsilon_1^2\epsilon_2^2\\
			& \quad +\big(2-\rho_2(L)\big)\epsilon_1\epsilon_2^2+\big(3+\rho(L)
			\big)\epsilon_2^2,\\
		&\varepsilon_3=\frac{2\rho_2(L)\epsilon_1-9}{4}-1,\\
		&\varepsilon_4=\big(3+4(1+c_1^{-1})\big)\rho^2(L)\epsilon_1^2-\rho_2(L)\epsilon_1+\rho(L)\\
			&\quad +2\epsilon_1^2+2\epsilon_1+3+1,\\	
		&\varepsilon_5=3\epsilon_1^2\epsilon_2^2\rho^2(L)+ (\epsilon_1\epsilon_2^2+\epsilon_2^2)\rho(L)-\frac{1}{2}\epsilon_1\epsilon_2\rho_2(L) \\
			&\quad + \epsilon_1^2\epsilon_2^2+2\epsilon_1\epsilon_2^2+\epsilon_2^2+\epsilon_1\epsilon_2+\frac{1}{2},\\	
		&\varepsilon_8=(1+\epsilon_1)\rho_2^{-1}(L)+\rho_2^{-2}(L),\\
		&c_1=\frac{\delta\omega r}{2},\\
		&c_2=c_1+2c_1^2,\\	
		& a_1=\frac{1}{2}(\varepsilon_3\epsilon_2-\varepsilon_4\epsilon_2^2),\\
		&b_{7,k}=2\frac{\rho_2^{-1}(L)}{\epsilon_2}.
			\end{align*}

	To simplify the analysis, we present the compact form of \eqref{nonconvex:kia-algo-dc-a}--\eqref{nonconvex:kia-algo-dc-q}.
	\begin{subequations}\label{nonconvex:alg-compress-compact}
			\begin{align}
				\bsy_{k+1}&=\bsy_{k}+\omega\bsq_k,\label{nonconvex:kia-algo-dc-compact-a}\\
				\bsz_{k+1}&=\bsz_{k}+\omega\bsL\bsq_k,\label{nonconvex:alg-compress-compact-b}\\
				\bsx_{k+1}&=\bsx_k-\alpha_k(\beta_k\bsL\hat{\bsx}_k+\gamma_k\bsv_k+\bsg_k^z),\label{nonconvex:kia-algo-dc-compact-x}\\
				\bsv_{k+1}&=\bsv_k+\alpha_k\gamma_k\bsL\hat{\bsx}_k,\label{nonconvex:kia-algo-dc-compact-v}\\
				\bsq_{k+1}&=\mathcal{C}(\bsx_{k+1}-\bsy_{k+1}).\label{nonconvex:kia-algo-dc-compact-q}
			\end{align}	
		\end{subequations}

		Now it is ready to prove Lemma~\ref{Lemma Consensus}--\ref{zerosg:lemma:sg:L1}.

\subsubsection{Proof of Lemma~\ref{Lemma Consensus}}

	We show the relation between $e_{1,k+1}$ and $e_{1,k}$.
		\begin{align}
			&\mathbb{E}_{\mathcal{B}_k}[e_{1,k+1}]=\mathbb{E}_{\mathcal{B}_k}\Big[\frac{1}{2}\|\bsx_{k+1} \|^2_{\bsE}\Big]\nonumber\\
			&=\mathbb{E}_{\mathcal{B}_k}\Big[\frac{1}{2}\|\bsx_k-\alpha_k(\beta_k\bsL\hat{\bsx}_k+\gamma_k\bsv_k+\bsg_k^{z}) \|^2_{\bsE}\Big]\nonumber\\
			&=\mathbb{E}_{\mathcal{B}_k}\Big[\frac{1}{2}\|\bsx_k\|^2_{\bsE}-\alpha_k\beta_k\bsx^\top_k\bsL\hat{\bsx}_k
			+\|\hat{\bsx}_k\|^2_{\frac{\alpha_k^2\beta_k^2}{2}\bsL^2}\nonumber\\
			&\quad-\alpha_k\gamma_k(\bsx^\top_k-\alpha_k\beta_k\hat{\bsx}_k^\top\bsL)\bsE
			\Big(\bsv_k+\frac{1}{\gamma_k}\bsg_k^{z}\Big)\nonumber\\
			&\quad+\Big\|\bsv_k+\frac{1}{\gamma_k}\bsg_k^{z}\Big\|^2_{\frac{\alpha_k^2\gamma_k^2}{2}\bsE}\Big]\nonumber\\
			&=\frac{1}{2}\|\bsx_k\|^2_{\bsE}-\alpha_k\beta_k\bsx^\top_k\bsL(\bsx_k+\hat{\bsx}_k-\bsx_k)
			+\|\hat{\bsx}_k\|^2_{\frac{\alpha_k^2\beta_k^2}{2}\bsL^2}\nonumber\\
			&\quad-\alpha_k\gamma_k(\bsx^\top_k-\alpha_k\beta_k\hat{\bsx}_k^\top\bsL)\bsE\Big(\bsv_k
			+\frac{1}{\gamma_k}\bsg_k^0\nonumber\\
			&\quad +\frac{1}{\gamma_k}\bsg_k^\mu-\frac{1}{\gamma_k}\bsg_k^0\Big)\nonumber\\
			&\quad+\mathbb{E}_{\mathcal{B}_k}\Big[\Big\|\bsv_k+\frac{1}{\gamma_k}\bsg_k^0
			+\frac{1}{\gamma_k}\bsg_k^{z}-\frac{1}{\gamma_k}\bsg_k^0\Big\|^2_{\frac{\alpha_k^2\gamma_k^2}{2}\bsE}\Big]\nonumber\\
			&\le\frac{1}{2}\|\bsx_k\|^2_{\bsE}-\|\bsx_k\|^2_{\alpha_k\beta_k\bsL}
			+\|\bsx_k\|^2_{\frac{\alpha_k\beta_k}{2}\bsL} \nonumber\\
			&\quad +\|\hat{\bsx}_k-\bsx_k\|^2_{\frac{\alpha_k\beta_k}{2}\bsL} 
			+\|\hat{\bsx}_k\|^2_{\frac{\alpha_k^2\beta_k^2}{2}\bsL^2} \nonumber\\
			&\quad-\alpha_k\gamma_k\bsx^\top_k\bsE\Big(\bsv_k+\frac{1}{\gamma_k}\bsg_k^0\Big) \nonumber\\
			&\quad +\frac{\alpha_k}{2}\|\bsx_k\|^2_{\bsE}
			+\frac{\alpha_k}{2}\|\bsg_k^\mu-\bsg_k^0\|^2\nonumber\\
			&\quad +\|\hat{\bsx}_k\|^2_{\frac{\alpha_k^2\beta_k^2}{2}\bsL^2} 
			 +\frac{\alpha_k^2\gamma_k^2}{2}\Big\|\bsv_k+\frac{1}{\gamma_k}\bsg_k^0\Big\|^2\nonumber\\
			&\quad+\|\hat{\bsx}_k\|^2_{\frac{\alpha_k^2\beta_k^2}{2}\bsL^2}
			+\frac{\alpha_k^2}{2}\|\bsg_k^\mu-\bsg_k^0\|^2\nonumber\\
			&\quad+\alpha_k^2\gamma_k^2\Big\|\bsv_k+\frac{1}{\gamma_k}\bsg_k^0\Big\|^2
			+\alpha_k^2\mathbb{E}_{\mathcal{B}_k}\Big[\|\bsg_k^z-\bsg_k^0\|^2\Big]\nonumber\\
			&=\frac{1}{2}\|\bsx_k\|^2_{\bsE}-\|\bsx_k\|^2_{\frac{\alpha_k\beta_k}{2}\bsL-\frac{\alpha_k}{2}\bsE}
			+\|\hat{\bsx}_k\|^2_{\frac{3\alpha_k^2\beta_k^2}{2}\bsL^2}\nonumber\\
			&\quad+\frac{\alpha_k}{2}(1+\alpha_k)\|\bsg_k^\mu-\bsg_k^0\|^2
			+\|\hat{\bsx}_k-\bsx_k\|^2_{\frac{\alpha_k\beta_k}{2}\bsL}\nonumber\\
			&\quad-\alpha_k\gamma_k(\hat{\bsx}_k+\bsx_k-\hat{\bsx}_k)^\top\bsE
			\Big(\bsv_k+\frac{1}{\gamma_k}\bsg_k^0\Big)\nonumber\\
			&\quad+\frac{3\alpha_k^2\gamma_k^2}{2}\Big\|\bsv_k+\frac{1}{\gamma_k}\bsg_k^0\Big\|^2+\alpha_k^2\mathbb{E}_{\mathcal{B}_k}\Big[\|\bsg_k^z-\bsg_k^0\|^2\Big] \nonumber\\
			&\le\frac{1}{2}\|\bsx_k\|^2_{\bsE}-\|\bsx_k\|^2_{\frac{\alpha_k\beta_k}{2}\bsL-\frac{\alpha_k}{2}\bsE}
			+\|\hat{\bsx}_k\|^2_{\frac{3\alpha_k^2\beta_k^2}{2}\bsL^2}\nonumber\\
			&\quad+\frac{\alpha_k}{2}(1+\alpha_k)\|\bsg_k^\mu-\bsg_k^0\|^2
			+\|\hat{\bsx}_k-\bsx_k\|^2_{\frac{\alpha_k}{2}(\beta_k\bsL+2\rho(L)\gamma_k\bsE)}\nonumber\\
			&\quad-\alpha_k\gamma_k\hat{\bsx}^\top_k\bsE\Big(\bsv_k+\frac{1}{\gamma_k}\bsg_k^0\Big)+\alpha_k^2\mathbb{E}_{\mathcal{B}_k}\Big[\|\bsg_k^z-\bsg_k^0\|^2\Big]\nonumber\\	
			&\quad+\frac{6\alpha_k^2\gamma_k^2+\alpha_k\gamma_k\rho^{-1}(L)}{4}\Big\|\bsv_k+\frac{1}{\gamma_k}\bsg_k^0\Big\|^2
			\nonumber\\
			&\le e_{1,k}-\|\bsx_k\|^2_{\frac{\alpha_k\beta_k}{2}\bsL-\frac{\alpha_k}{2}\bsE
				- \alpha_k(1+5\alpha_k) \ell ^2\bsE} + \|\hat{\bsx}_k\|^2_{\frac{3\alpha_k^2\beta_k^2}{2}\bsL^2}\nonumber\\
			&\quad + n \ell ^2\alpha_k(1+5\alpha_k)\mu^2_k + \frac{\alpha_k}{2}(\beta_k+2\gamma_k)\rho(L)\|\bsx_k-\hat{\bsx}_k\|^2 \nonumber\\
			&\quad -\alpha_k\gamma_k\hat{\bsx}^\top_k\bsE\Big(\bsv_k+\frac{1}{\gamma_k}\bsg_k^0\Big)  + 2\alpha_k^2\mathbb{E}_{\mathcal{B}_k}[\|\bsg_k^z\|^2]\nonumber\\
			&\quad +\Big\|\bsv_k+\frac{1}{\gamma_k}\bsg_k^0\Big\|^2_{\frac{6\alpha_k^2\gamma_k^2\rho(L)+\alpha_k\gamma_k}{4}\bsF},\label{nonconvex:v1k_B_k}
		\end{align}
		where the second, third and fourth equalities hold due to \eqref{nonconvex:kia-algo-dc-compact-x}, \eqref{nonconvex:KL-L-eq} and \eqref{zerosg:rand-grad-esti1}, respectively; the first and second inequalities hold due to \eqref{nonconvex:lemma:inequality-arithmetic-equ} and $\rho(\bsE)=1$; and the last  inequality holds due to \eqref{zerosg:rand-grad-esti8}, \eqref{zerosg:rand-grad-esti6} and \eqref{nonconvex:lemma-eq5}.
		Then, since  $\mathcal{B}_k$ and $\mathcal{C}_k$ are independent and $\mathcal{A}_k = \mathcal{B}_k \cup \mathcal{C}_k$, taking the expectation with respect to $\mathcal{A}_k$ on both sides of \eqref{nonconvex:v1k_B_k} yields \eqref{nonconvex:v1k}.

		\subsubsection{Proof of Lemma~\ref{Lemma Optimality 1}}
We show the relation between $e_{2,k+1}$ and $e_{2,k}$.
		\begin{align}
			&e_{2,k+1}=\frac{1}{2}\Big\|\bsv_{k+1}+\frac{1}{\gamma_{k+1}}\bsg_{k+1}^0\Big\|^2_{\frac{\beta_k+\gamma_k}{\gamma_k}\bsF}\nonumber\\
			&\quad=\frac{1}{2}\Big\|\bsv_{k+1}+\frac{1}{\gamma_{k}}\bsg_{k+1}^0
			+\Big(\frac{1}{\gamma_{k+1}}-\frac{1}{\gamma_{k}}\Big)\bsg_{k+1}^0\Big\|^2_{\frac{\beta_k+\gamma_k}{\gamma_k}\bsF}\nonumber\\
			&\quad\le\frac{1}{2}(1+\Delta_k)\Big\|\bsv_{k+1}+\frac{1}{\gamma_{k}}\bsg_{k+1}^0\Big\|^2_{\frac{\beta_k+\gamma_k}{\gamma_k}\bsF}\nonumber\\
			&\qquad+\frac{1}{2}(\Delta_k+\Delta_k^2)\|\bsg_{k+1}^0\|^2_{\frac{\beta_k+\gamma_k}{\gamma_k}\bsF},\label{zerosg:v2k-1}
		\end{align}
		where the inequality holds due to \eqref{nonconvex:lemma:inequality-arithmetic-equ}.
		
		For the first term on the right-hand side of \eqref{zerosg:v2k-1}, we have
		\begin{align}
			&\frac{1}{2}\Big\|\bsv_{k+1}+\frac{1}{\gamma_k}\bsg_{k+1}^0\Big\|^2_{\frac{\beta_k+\gamma_k}{\gamma_k}\bsF}\nonumber\\
			&=\frac{1}{2}\Big\|\bsv_k+\frac{1}{\gamma_k}\bsg_{k}^0+\alpha_k\gamma_k\bsL\hat{\bsx}_k
			+\frac{1}{\gamma_k}(\bsg_{k+1}^0-\bsg_{k}^0) \Big\|^2_{\frac{\beta_k+\gamma_k}{\gamma_k}\bsF}\nonumber\\
			&=\frac{1}{2}\Big\|\bsv_k+\frac{1}{\gamma_k}\bsg_{k}^0\Big\|^2_{\frac{\beta_k+\gamma_k}{\gamma_k}\bsF}\nonumber\\ 
			&\quad +\alpha_k(\beta_k+\gamma_k)\hat{\bsx}^\top_k\bsE\Big(\bsv_k+\frac{1}{\gamma_k}\bsg_k^0\Big)\nonumber\\
			&\quad+\|\hat{\bsx}_k\|^2_{\frac{\alpha_k^2\gamma_k}{2}(\beta_k+\gamma_k)\bsL}
			+\frac{1}{2\gamma_k^2}\|\bsg_{k+1}^0-\bsg_{k}^0\|^2_{\frac{\beta_k+\gamma_k}{\gamma_k}\bsF}\nonumber\\
			&\quad +\frac{\beta_k+\gamma_k}{\gamma_k^2}\Big(\bsv_k+\frac{1}{\gamma_k}\bsg_{k}^0
			\Big)^\top\bsF(\bsg_{k+1}^0-\bsg_{k}^0)\nonumber\\
			&\quad +\alpha_k\frac{\beta_k+\gamma_k}{\gamma_k}\hat{\bsx}_k^\top\bsE(\bsg_{k+1}^0-\bsg_{k}^0)\nonumber\\
			&\le\frac{1}{2}\Big\|\bsv_k+\frac{1}{\gamma_k}\bsg_{k}^0\Big\|^2_{\frac{\beta_k+\gamma_k}{\gamma_k}\bsF}\nonumber\\
			 &\quad +\alpha_k(\beta_k+\gamma_k)\hat{\bsx}^\top_k\bsE\Big(\bsv_k+\frac{1}{\gamma_k}\bsg_k^0\Big)\nonumber\\
			&\quad+\|\hat{\bsx}_k\|^2_{\frac{\alpha_k^2\gamma_k}{2}(\beta_k+\gamma_k)\bsL}
			+\|\bsg_{k+1}^0-\bsg_{k}^0\|^2_{\frac{\beta_k+\gamma_k}{2\gamma_k^3}\bsF}\nonumber\\
			&\quad+\Big\|\bsv_k+\frac{1}{\gamma_k}\bsg_{k}^0\Big\|^2_{\frac{\alpha_k\gamma_k}{4}\bsF}
			+\|\bsg_{k+1}^0-\bsg_{k}^0\|^2_{\frac{(\beta_k+\gamma_k)^2}{\alpha_k\gamma_k^5}\bsF}\nonumber\\
			&\quad+\|\hat{\bsx}_k\|^2_{\frac{\alpha_k^2(\beta_k+\gamma_k)^2}{2}\bsE}
			+\frac{1}{2\gamma_k^2}\|\bsg_{k+1}^0-\bsg_{k}^0\|^2\nonumber\\
			&\le \frac{1}{2}\Big\|\bsv_k+\frac{1}{\gamma_k}\bsg_{k}^0\Big\|^2_{\frac{\beta_k+\gamma_k}{\gamma_k}\bsF} 
			 +\Big\|\bsv_k+\frac{1}{\gamma_k}\bsg_{k}^0\Big\|^2_{\frac{\alpha_k\gamma_k}{4}\bsF} \nonumber\\ 
			&\quad +\alpha_k(\beta_k+\gamma_k)\hat{\bsx}^\top_k\bsE\Big(\bsv_k+\frac{1}{\gamma_k}\bsg_k^0\Big)\nonumber\\
			&\quad+\|\hat{\bsx}_k\|^2_{\frac{\alpha_k^2}{2}(\beta_k\gamma_k+\gamma_k^2)\bsL+\frac{\alpha_k^2}{2}(\beta_k+\gamma_k)^2\bsE} \nonumber\\
			&\quad +\big((\frac{\beta_k+\gamma_k}{2\gamma_k^3}+\frac{(\beta_k+\gamma_k)^2}{\alpha_k\gamma_k^5})\rho_2^{-1}(L)
			  + \frac{1}{2\gamma_k^2}\big)\nonumber\\
			&\quad \times\|\bsg_{k+1}^0-\bsg_{k}^0\|^2 \nonumber\\
			&\le\frac{1}{2}\Big\|\bsv_k+\frac{1}{\gamma_k}\bsg_{k}^0\Big\|^2_{\frac{\beta_k+\gamma_k}{\gamma_k}\bsF}
			+\Big\|\bsv_k+\frac{1}{\gamma_k}\bsg_{k}^0\Big\|^2_{\frac{\alpha_k\gamma_k}{4}\bsF}\nonumber\\ 
			&\quad +\alpha_k(\beta_k+\gamma_k)\hat{\bsx}^\top_k\bsE\Big(\bsv_k+\frac{1}{\gamma_k}\bsg_k^0\Big)\nonumber\\
			&\quad+\|\hat{\bsx}_k\|^2_{\frac{\alpha_k^2}{2}(\beta_k\gamma_k+\gamma_k^2)\bsL+\frac{\alpha_k^2}{2}(\beta_k+\gamma_k)^2\bsE} \nonumber\\
			&\quad +(\epsilon_4\alpha_k+\epsilon_5\alpha_k^2) \ell ^2\|\bar{\bsg}_k^z\|^2
			,\label{nonconvex:v2k}
		\end{align}
		where the first equality holds due to \eqref{nonconvex:kia-algo-dc-compact-v}; the second equality holds due to \eqref{nonconvex:KL-L-eq} and \eqref{nonconvex:lemma-eq3}; the first inequality holds due to \eqref{nonconvex:lemma:inequality-arithmetic-equ} and $\rho(\bsE)=1$; the second inequality holds due to \eqref{nonconvex:lemma-eq5}; and the last inequality holds due to \eqref{zerosg:gg-rand-pd}.

		For the second term on the right-hand side of \eqref{zerosg:v2k-1}, we have
		\begin{align}\label{zerosg:v2k-4}
			\|\bsg_{k+1}^0\|^2_{\frac{\beta_k+\gamma_k}{\gamma_k}\bsF}
			\le\varepsilon_1\|\bsg_{k+1}^0\|^2.
		\end{align}


		Then, from \eqref{zerosg:v2k-1}--\eqref{zerosg:v2k-4}, we have
		\begin{align}
			&e_{2,k+1}\le e_{2,k}\nonumber\\
			&\quad +(1+\Delta_k)\alpha_k(\beta_k+\gamma_k)\hat{\bsx}^\top_k\bsE\Big(\bsv_k+\frac{1}{\gamma_k}\bsg_k^0\Big)\nonumber\\
			&\quad+\|\hat{\bsx}_k\|^2_{(1+\Delta_k)\frac{1}{2}\alpha_k^2(\beta_k\gamma_k+\gamma_k^2)\bsL + (1+\Delta_k)\frac{1}{2}\alpha_k^2(\beta_k+\gamma_k)^2\bsE} \nonumber\\
			&\quad +\Big\|\bsv_k+\frac{1}{\gamma_k}\bsg_{k}^0\Big\|^2_{\big(\alpha_k\frac{\gamma_k}{4}+\Delta_k(\frac{\beta_k+\gamma_k}{2\gamma_k}+\frac{\alpha_k\gamma_k}{4})\big)\bsF}\nonumber\\
			&\quad +(1+\Delta_k)(\epsilon_4\alpha_k+\epsilon_5\alpha_k^2) \ell ^2\|\bar{\bsg}_k^z\|^2 \nonumber\\
			&\quad +\frac{1}{2}\varepsilon_1(\Delta_k+\Delta_k^2)\|\bsg_{k+1}^0\|^2.
			\label{zerosg:v2k_B_k}
		\end{align}
		Since $\mathcal{B}_k$ and $\mathcal{C}_k$ are independent and $\mathcal{A}_k = \mathcal{B}_k \cup \mathcal{C}_k$, taking the expectation with respect to $\mathcal{A}_k$ on both sides of \eqref{zerosg:v2k_B_k} yields \eqref{zerosg:v2k}.

		\subsubsection{Proof of Lemma~\ref{Lemma Cross}}
		We show the relation between $e_{3,k+1}$ and $e_{3,k}$.
		\begin{align}
		&e_{3,k+1}\nonumber\\
		&=\bsx_{k+1}^\top\bsE\bsF\Big(\bsv_{k+1}+\frac{1}{\gamma_{k+1}}\bsg_{k+1}^0\Big)\nonumber\\
		&=\bsx_{k+1}^\top\bsE\bsF\Big(\bsv_{k+1}+\frac{1}{\gamma_{k}}\bsg_{k+1}^0
		+\Big(\frac{1}{\gamma_{k+1}}-\frac{1}{\gamma_{k}}\Big)\bsg_{k+1}^0\Big)\nonumber\\
		&=\bsx_{k+1}^\top\bsE\bsF\Big(\bsv_{k+1}+\frac{1}{\gamma_{k}}\bsg_{k+1}^0\Big)
		-\Delta_k\bsx_{k+1}^\top\bsE\bsF\bsg_{k+1}^0\nonumber\\
		&\le\bsx_{k+1}^\top\bsE\bsF\Big(\bsv_{k+1}+\frac{1}{\gamma_{k}}\bsg_{k+1}^0\Big) \nonumber\\
		 &\quad +\frac{1}{2}\Delta_k(\|\bsx_{k+1}\|^2_{\bsE}+\|\bsg_{k+1}^0\|_{\bsF^2}^2).\label{zerosg:v3k-1}
		\end{align}

	For the first term on the right-hand side of \eqref{zerosg:v3k-1}, we have
	\begin{align}
		&\mathbb{E}_{\mathcal{B}_k}\Big[\bsx_{k+1}^\top\bsE\bsF\Big(\bsv_{k+1}+\frac{1}{\gamma_k}\bsg_{k+1}^0\Big)\Big]\nonumber\\
		&=\mathbb{E}_{\mathcal{B}_k}\Big[(\bsx_k-\alpha_k(\beta_k\bsL\hat{\bsx}_k+\gamma_k\bsv_k+\bsg_k^0+\bsg_k^z-\bsg_k^0))^\top
		\bsE\bsF\nonumber\\
		&\quad \times\Big(\bsv_k
		+\frac{1}{\gamma_k}\bsg_{k}^0+\alpha_k\gamma_k\bsL\hat{\bsx}_k+\frac{1}{\gamma_k}(\bsg_{k+1}^0-\bsg_{k}^0)\Big)\Big]\nonumber\\
		&=(\bsx_k^\top\bsE\bsF-\alpha_k(\beta_k+\alpha_k\gamma_k^2)\hat{\bsx}_k^\top\bsE)
		\Big(\bsv_k+\frac{1}{\gamma_k}\bsg_{k}^0\Big)\nonumber\\
		&\quad+\alpha_k\gamma_k\bsx_k^\top\bsE\hat{\bsx}_k-\|\hat{\bsx}_k\|^2_{\alpha_k^2\beta_k\gamma_k\bsL}\nonumber\\
		&\quad
		+\frac{1}{\gamma_k}(\bsx_k^\top\bsE\bsF-\alpha_k\beta_k\hat{\bsx}_k^\top\bsE)\mathbb{E}_{\mathcal{B}_k}[\bsg_{k+1}^0-\bsg_{k}^0]\nonumber\\
		&\quad-\alpha_k(\gamma_k\bsv_k+\bsg_{k}^0+\bsg_k^\mu-\bsg_k^0-\bar{\bsg}_k^\mu)^\top\bsF
		\Big(\bsv_k+\frac{1}{\gamma_k}\bsg_{k}^0\Big)\nonumber\\
		&\quad 
		-\alpha_k\Big(\bsv_k+\frac{1}{\gamma_k}\bsg_{k}^0\Big)^\top\bsE\bsF\mathbb{E}_{\mathcal{B}_k}[\bsg_{k+1}^0-\bsg_{k}^0]\nonumber\\
		&\quad-\alpha_k^2\gamma_k(\bsg_k^\mu-\bsg_k^0)^\top\bsE\hat{\bsx}_k \nonumber\\
		&\quad-\mathbb{E}_{\mathcal{B}_k}\Big[\frac{\alpha_k}{\gamma_k}(\bsg_k^z-\bsg_k^0)^\top\bsE\bsF(\bsg_{k+1}^0-\bsg_{k}^0)\Big]\label{Pre_Lemma3_Reproof}\\			
		&\le(\bsx_k^\top\bsE\bsF-\alpha_k\beta_k\hat{\bsx}_k^\top\bsE)
		\Big(\bsv_k+\frac{1}{\gamma_k}\bsg_{k}^0\Big) \nonumber\\
	    &\quad +\|\hat{\bsx}_k\|^2_{\frac{\alpha_k^2\gamma_k^2}{2}\bsE}
		+\frac{\alpha_k^2\gamma_k^2}{2}\Big\|\bsv_k+\frac{1}{\gamma_k}\bsg_{k}^0\Big\|^2 \nonumber\\
		&\quad +\|\bsx_k\|^2_{\frac{\alpha_k\gamma_k}{4}\bsE} + \|\hat{\bsx}_k\|^2_{\alpha_k\gamma_k(\bsE-\alpha_k\beta_k\bsL)}\nonumber\\
		&\quad+\|\bsx_k\|^2_{\frac{\alpha_k}{2}\bsE} + \mathbb{E}_{\mathcal{B}_k}\big[\|\bsg_{k+1}^0-\bsg_{k}^0\|^2_{\frac{1}{2\alpha_k\gamma_k^2}\bsF^2}\big]
		\nonumber\\
		&\quad +\|\hat{\bsx}_k\|^2_{\frac{\alpha^2_k\beta^2_k}{2}\bsE}
		+\frac{1}{2\gamma^2_k}\mathbb{E}_{\mathcal{B}_k}[\|\bsg_{k+1}^0-\bsg_{k}^0\|^2]\nonumber\\
		&\quad -\Big\|\bsv_k+\frac{1}{\gamma_k}\bsg_{k}^0\Big\|^2_{\alpha_k\gamma_k\bsF}\nonumber\\
		&\quad+\frac{\alpha_k}{4}\|\bsg_k^\mu-\bsg_k^0\|^2
			+\Big\|\bsv_k+\frac{1}{\gamma_k}\bsg_{k}^0\Big\|^2_{\alpha_k\bsF^2}\nonumber\\
				&\quad	+ \frac{\alpha_k}{\gamma_k}\rho_2^{-1}(L)\|\bar{\bsg}_k^\mu\|^2 
				+  \Big\|\bsv_k+\frac{1}{\gamma_k}\bsg_{k}^0\Big\|^2_{\frac{1}{4}\rho_2(L)\alpha_k\gamma_k \bsF^2} \nonumber\\
		&\quad+ \frac{\alpha_k^2\gamma_k^2}{2}\Big\|\bsv_k+\frac{1}{\gamma_k}\bsg_{k}^0\Big\|^2
		+\mathbb{E}_{\mathcal{B}_k}\big[\|\bsg_{k+1}^0-\bsg_{k}^0\|^2_{\frac{1}{2\gamma_k^2}\bsF^2}\big] \nonumber\\
		&\quad +\frac{\alpha_k^2}{2}\|\bsg_k^\mu-\bsg_k^0\|^2
		 +\|\hat{\bsx}_k\|^2_{\frac{\alpha_k^2\gamma_k^2}{2}\bsE} \nonumber\\
		&\quad +\frac{\alpha_k^2}{2}\mathbb{E}_{\mathcal{B}_k}\big[\|\bsg_k^z-\bsg_k^0\|^2\big]
		+\mathbb{E}_{\mathcal{B}_k}\big[\|\bsg_{k+1}^0-\bsg_{k}^0\|^2_{\frac{1}{2\gamma_k^2}\bsF^2}\big] \label{Lemma3_Reproof}\\
		&=(\bsx_k^\top\bsE\bsF-\alpha_k\beta_k\hat{\bsx}_k^\top\bsE)\Big(\bsv_k+\frac{1}{\gamma_k}\bsg_{k}^0\Big) \nonumber\\
		&\quad +\|\bsx_k\|^2_{\frac{\alpha_k(\gamma_k+2)}{4}\bsE}
		 +\|\hat{\bsx}_k\|^2_{\alpha_k\gamma_k\bsE
			+\alpha_k^2\big((\frac{1}{2}\beta_k^2+\gamma_k^2)\bsE-\beta_k\gamma_k\bsL\big)} \nonumber\\
		&\quad -\Big\|\bsv_k+\frac{1}{\gamma_k}\bsg_{k}^0\Big\|^2_{\alpha_k(\gamma_k\bsF-\frac{1}{4}\rho_2(L)\gamma_k\bsF^2-\bsF^2)-\alpha_k^2\gamma_k^2{\bf I}_{np}}\nonumber\\
		&\quad +\mathbb{E}_{\mathcal{B}_k}\big[\|\bsg_{k+1}^0-\bsg_{k}^0\|^2_{\frac{1}{2\gamma_k^2}{\bf I}_{np}+\frac{1}{2\gamma_k^2}(2+\frac{1}{\alpha_k})\bsF^2}\big] \nonumber\\
		&\quad +\frac{\alpha_k}{4}(1+2\alpha_k)\|\bsg_k^\mu-\bsg_k^0\|^2
		 +\frac{\alpha_k^2}{2}\mathbb{E}_{\mathcal{B}_k}\big[\|\bsg_k^z-\bsg_k^0\|^2\big] \nonumber\\
		&\quad+\frac{\alpha_k}{\gamma_k}\rho_2^{-1}(L)\|\bar{\bsg}_k^\mu\|^2 \nonumber\\
		&\le(\bsx_k^\top\bsE\bsF-\alpha_k\beta_k\hat{\bsx}_k^\top\bsE)\Big(\bsv_k+\frac{1}{\gamma_k}\bsg_{k}^0\Big) \nonumber\\
		&\quad +\|\bsx_k\|^2_{\frac{\alpha_k(\gamma_k+2)}{4}\bsE}
		 +\|\hat{\bsx}_k\|^2_{\alpha_k\gamma_k\bsE
			+\alpha_k^2\big((\frac{1}{2}\beta_k^2+\gamma_k^2)\bsE-\beta_k\gamma_k\bsL\big)} \nonumber\\
		&\quad -\Big\|\bsv_k+\frac{1}{\gamma_k}\bsg_{k}^0\Big\|^2_{\alpha_k(\frac{3}{4}\gamma_k-\rho_2^{-1}(L))\bsF-\alpha_k^2\gamma_k^2\rho(L)\bsF}\nonumber\\
		&\quad + \big(\frac{1}{2\gamma_k^2}+\frac{1}{2\gamma_k^2}(2+\frac{1}{\alpha_k})\rho_2^{-2}(L)\big)\mathbb{E}_{\mathcal{B}_k}\big[\|\bsg_{k+1}^0-\bsg_{k}^0\|^2\big]\nonumber\\
		&\quad +\frac{\alpha_k}{4}(1+2\alpha_k)\|\bsg_k^\mu-\bsg_k^0\|^2
		 +\frac{\alpha_k^2}{2}\mathbb{E}_{\mathcal{B}_k}\big[\|\bsg_k^z-\bsg_k^0\|^2\big] \nonumber\\
		&\quad+\frac{\alpha_k}{\gamma_k}\rho_2^{-1}(L)\|\bar{\bsg}_k^\mu\|^2 
		\nonumber\\
		&\le\bsx_k^\top\bsE\bsF\Big(\bsv_k+\frac{1}{\gamma_k}\bsg_{k}^0\Big)
		-(1+\Delta_k)\alpha_k\beta_k\hat{\bsx}_k^\top\bsE\big(\bsv_k+\frac{1}{\gamma_k}\bsg_{k}^0\big)\nonumber\\
		&\quad +\Delta_k\alpha_k\beta_k\hat{\bsx}_k^\top\bsE\Big(\bsv_k+\frac{1}{\gamma_k}\bsg_{k}^0\Big)\nonumber\\
		&\quad+\|\bsx_k\|^2_{\alpha_k(\frac{\gamma_k+2}{4}+\frac{1}{2} \ell ^2)\bsE+3\alpha_k^2 \ell ^2\bsE}\nonumber\\
		&\quad +\|\hat{\bsx}_k\|^2_{\alpha_k\gamma_k\bsE
		 +\alpha_k^2\big((\frac{1}{2}\beta_k^2+\gamma_k^2)\bsE-\beta_k\gamma_k\bsL\big)} \nonumber\\
		&\quad -\Big\|\bsv_k+\frac{1}{\gamma_k}\bsg_{k}^0\Big\|^2_{\alpha_k(\frac{3}{4}\gamma_k-\rho_2^{-1}(L))\bsF-\alpha_k^2\gamma_k^2\rho(L)\bsF}\nonumber\\
		&\quad +(\alpha_k\epsilon_6+\alpha_k^2\epsilon_7) \ell ^2\mathbb{E}_{\mathcal{B}_k}\big[\|\bar{\bsg}_k^z\|^2\big]  +n \ell ^2\alpha_k(\frac{1}{2}+3\alpha_k)\mu_k^2 \nonumber\\
		&\quad + \alpha^2_k\mathbb{E}_{\mathcal{B}_k}[\|\bsg_k^z\|^2] +\frac{\alpha_k}{\gamma_k}\rho_2^{-1}(L)\|\bar{\bsg}_k^\mu\|^2 
		,\label{nonconvex:v3k}
	\end{align}
	where the first equality holds due to \eqref{nonconvex:kia-algo-dc-compact-x} and \eqref{nonconvex:kia-algo-dc-compact-v}; the second equality holds since \eqref{nonconvex:KL-L-eq}, \eqref{nonconvex:lemma-eq3}, $\bsE={\bf I}_{nd}-\bsH$, \eqref{zerosg:rand-grad-esti1}, and that $x_{i,k}$ and $v_{i,k}$ are independent of $\mathcal{B}_k$; the first inequality holds due to \eqref{nonconvex:lemma:inequality-arithmetic-equ} and $\rho(\bsE)=1$; the second inequality holds due to \eqref{nonconvex:lemma-eq5}; and the last  inequality holds due to \eqref{zerosg:rand-grad-esti8}, \eqref{zerosg:rand-grad-esti6}, and \eqref{zerosg:gg-rand-pd}.

	For the third term on the right-hand side of \eqref{nonconvex:v3k}, we have
		\begin{align}\label{zerosg:v3k-3}
			&\Delta_k\alpha_k\beta_k\hat{\bsx}_k^\top\bsE\Big(\bsv_k+\frac{1}{\gamma_k}\bsg_{k}^0\Big) \nonumber\\
			&\le\|\hat{\bsx}_k\|^2_{\frac{1}{2}\Delta_k\alpha_k\beta_k\bsE}
			+\Big\|\bsv_k+\frac{1}{\gamma_k}\bsg_k^0\Big\|^2_{\frac{1}{2}\rho(L)\Delta_k\alpha_k\beta_k\bsF}.
		\end{align}
		
		Then, from \eqref{zerosg:v3k-1}--\eqref{zerosg:v3k-3} and \eqref{nonconvex:lemma-eq5}, we have
		\begin{align}
			&\mathbb{E}_{\mathcal{B}_k}[e_{3,k+1}]\nonumber\\
			&\quad\le e_{3,k}
			-(1+\Delta_k)\alpha_k\beta_k\hat{\bsx}_k^\top\bsE\Big(\bsv_k+\frac{1}{\gamma_k}\bsg_{k}^0\Big)\nonumber\\
			&\quad+\|\bsx_k\|^2_{\alpha_k(\frac{\gamma_k+2}{4}+\frac{1}{2} \ell ^2)\bsE+3\alpha_k^2 \ell ^2\bsE}\nonumber\\
			&\quad +\|\hat{\bsx}_k\|^2_{\alpha_k\gamma_k\bsE
			 +\alpha_k^2\big((\frac{1}{2}\beta_k^2+\gamma_k^2)\bsE-\beta_k\gamma_k\bsL\big)+\frac{1}{2}\Delta_k\alpha_k\beta_k\bsE} \nonumber\\
			&\quad -\Big\|\bsv_k+\frac{1}{\gamma_k}\bsg_{k}^0\Big\|^2_{\big(\alpha_k(\frac{3}{4}\gamma_k-\rho_2^{-1}(L))-\alpha_k^2\gamma_k^2\rho(L) -\frac{1}{2}\rho(L)\Delta_k\alpha_k\beta_k\big)\bsF}\nonumber\\
			&\quad +(\alpha_k\epsilon_6+\alpha_k^2\epsilon_7) \ell ^2\mathbb{E}_{\mathcal{B}_k}\big[\|\bar{\bsg}_k^z\|^2\big]  +n \ell ^2\alpha_k(\frac{1}{2}+3\alpha_k)\mu_k^2 \nonumber\\
			&\quad + \alpha^2_k\mathbb{E}_{\mathcal{B}_k}[\|\bsg_k^z\|^2] +\frac{\alpha_k}{\gamma_k}\rho_2^{-1}(L)\|\bar{\bsg}_k^\mu\|^2 \nonumber\\
			&\quad+\frac{1}{2}\Delta_k\mathbb{E}_{\mathcal{B}_k}[2e_{1,k+1}+\rho_2^{-2}(L)\|\bsg_{k+1}^0\|^2].
		\label{zerosg:v3k_B_k}
		\end{align}
		Since $\mathcal{B}_k$ and $\mathcal{C}_k$ are independent and $\mathcal{A}_k = \mathcal{B}_k \cup \mathcal{C}_k$, taking the expectation with respect to $\mathcal{A}_k$ on both sides of \eqref{zerosg:v3k_B_k} yields \eqref{zerosg:v3k}.

		\subsubsection{Proof of Lemma~\ref{Lemma Optimality 2}}
		We show the relation between $e_{4,k+1}$ and $e_{4,k}$.	
		\begin{align}
			&\mathbb{E}_{\mathcal{B}_k}[e_{4,k+1}]
			=\mathbb{E}_{\mathcal{B}_k}[\tilde{f}(\bar{\bsx}_{k+1})-nf^*]\nonumber\\
			&=\mathbb{E}_{\mathcal{B}_k}[\tilde{f}(\bar{\bsx}_k)-nf^*+\tilde{f}(\bar{\bsx}_{k+1})
			-\tilde{f}(\bar{\bsx}_k)]\nonumber\\
			&\le\mathbb{E}_{\mathcal{B}_k}\Big[\tilde{f}(\bar{\bsx}_k)-nf^*
			-\alpha_k(\bar{\bsg}_{k}^z)^\top\bsg^0_k
			+\frac{1}{2}\alpha^2_k \ell \|\bar{\bsg}_{k}^z\|^2\Big]\nonumber\\
			&=e_{4,k}
			-\alpha_k(\bar{\bsg}_{k}^\mu)^\top\bsg^0_k
			+\frac{1}{2}\alpha^2_k \ell \mathbb{E}_{\mathcal{B}_k}[\|\bar{\bsg}_{k}^z\|^2]\nonumber\\
			&=e_{4,k}
			-\alpha_k(\bar{\bsg}_{k}^\mu)^\top\bar{\bsg}^0_k
			+\frac{1}{2}\alpha^2_k \ell \mathbb{E}_{\mathcal{B}_k}[\|\bar{\bsg}_{k}^z\|^2]\nonumber\\
			&=e_{4,k}
			-\frac{1}{2}\alpha_k(\bar{\bsg}_{k}^\mu)^\top(\bar{\bsg}^\mu_k+\bar{\bsg}^0_k-\bar{\bsg}^\mu_k)\nonumber\\
			&\quad-\frac{1}{2}\alpha_k(\bar{\bsg}^\mu_{k}-\bar{\bsg}^0_k+\bar{\bsg}^0_k)^\top\bar{\bsg}^0_k
			+\frac{1}{2}\alpha^2_k \ell \mathbb{E}_{\mathcal{B}_k}[\|\bar{\bsg}_{k}^z\|^2]\nonumber\\
			&\le e_{4,k}-\frac{1}{4}\alpha_k(\|\bar{\bsg}^\mu_{k}\|^2
			-\|\bar{\bsg}^0_k-\bar{\bsg}^\mu_k\|^2+\|\bar{\bsg}_{k}^0\|^2\nonumber\\
			&\quad-\|\bar{\bsg}^0_k-\bar{\bsg}^\mu_k\|^2)
			+\frac{1}{2}\alpha^2_k \ell \mathbb{E}_{\mathcal{B}_k}[\|\bar{\bsg}_{k}^z\|^2]\nonumber\\
			&= e_{4,k}-\frac{1}{4}\alpha_k\|\bar{\bsg}^\mu_{k}\|^2
			+\frac{1}{2}\alpha_k\|\bar{\bsg}^0_k-\bar{\bsg}^\mu_k\|^2\nonumber\\
			&\quad-\frac{1}{4}\alpha_k\|\bar{\bsg}_{k}^0\|^2
			+\frac{1}{2}\alpha^2_k \ell \mathbb{E}_{\mathcal{B}_k}[\|\bar{\bsg}_{k}^z\|^2]\nonumber\\
			&\le e_{4,k} - \frac{1}{4}\alpha_k\|\bar{\bsg}^\mu_{k}\|^2 + \|\bsx_k\|^2_{\alpha_k \ell ^2\bsE}\nonumber\\
			&\quad + n \ell ^2\alpha_k\mu^2_k-\frac{1}{4}\alpha_k\|\bar{\bsg}_{k}^0\|^2
			+\frac{1}{2}\alpha^2_k \ell \mathbb{E}_{\mathcal{B}_k}[\|\bar{\bsg}^z_{k}\|^2],\label{zerosg:v4k_B_k}
			\end{align}
		where the first inequality holds since that $\tilde{f}$ is smooth and \eqref{nonconvex:lemma:lipschitz}; the third equality holds since \eqref{zerosg:rand-grad-esti1} and that $x_{i,k}$ and $v_{i,k}$ are independent of $\mathcal{B}_k$; the fourth equality holds due to $(\bar{\bsg}_{k}^\mu)^\top\bsg^0_k=(\bsg_{k}^\mu)^\top\bsH\bsg^0_k=(\bsg_{k}^\mu)^\top\bsH\bsH\bsg^0_k
		=(\bar{\bsg}_{k}^\mu)^\top\bar{\bsg}^0_k$; the second inequality holds due to \eqref{nonconvex:lemma:inequality-arithmetic-equ}; and the last inequality holds due to \eqref{zerosg:rand-grad-esti9}.
		Then, taking the expectation with respect to $\mathcal{A}_k$ on both sides of \eqref{zerosg:v4k_B_k} yields \eqref{zerosg:v4k}.

		\subsubsection{Proof of Lemma~\ref{Lemma Compression}}
		We show the relation between $e_{5,k+1}$ and $e_{5,k}$.	
		
		Denote $\mathcal{C}_r(\cdot)=\mathcal{C}(\cdot)/r$, then we have
		\begin{align}
			&\mathbb{E}_{\mathcal{C}_k}[\|\bsx_{k+1}-\bsy_{k+1}\|^2]\nonumber\\
			&=\mathbb{E}_{\mathcal{C}_k}[\|\bsx_{k+1}-\bsx_{k}+\bsx_{k}-\bsy_{k}-\omega\bsq_k\|^2]\nonumber\\
			&=\mathbb{E}_{\mathcal{C}_k}[\|\bsx_{k+1}-\bsx_{k}+(1-\omega r)(\bsx_{k}-\bsy_{k})\nonumber\\
			&\quad+\omega r(\bsx_{k}-\bsy_{k}-\mathcal{C}_r(\bsx_{k}-\bsy_{k}))\|^2]\nonumber\\
			&\le(1+c_1^{-1})\mathbb{E}_{\mathcal{C}_k}[\|\bsx_{k+1}-\bsx_{k}\|^2]\nonumber\\
			&\quad+(1+c_1)\mathbb{E}_{\mathcal{C}_k}[\|(1-\omega r)(\bsx_{k}-\bsy_{k})\nonumber\\
			&\quad+\omega r(\bsx_{k}-\bsy_{k}-\mathcal{C}_r(\bsx_{k}-\bsy_{k}))\|^2]\nonumber\\
			&\le(1+c_1^{-1})\mathbb{E}_{\mathcal{C}_k}[\|\bsx_{k+1}-\bsx_{k}\|^2]\nonumber\\
			&\quad+(1+c_1)(1-\omega r)\mathbb{E}_{\mathcal{C}_k}[\|\bsx_{k}-\bsy_{k}\|^2]\nonumber\\
			&\quad+(1+c_1)\omega r\mathbb{E}_{\mathcal{C}_k}[\|\bsx_{k}-\bsy_{k}-\mathcal{C}_r(\bsx_{k}-\bsy_{k})\|^2]\nonumber\\
			&\le(1+c_1^{-1})\mathbb{E}_{\mathcal{C}_k}[\|\bsx_{k+1}-\bsx_{k}\|^2]\nonumber\\
			&\quad+(1+c_1)(1-\omega r)\|\bsx_{k}-\bsy_{k}\|^2\nonumber\\
			&\quad+(1+c_1)\omega r(1-\delta)\|\bsx_{k}-\bsy_{k}\|^2\nonumber\\
			&=(1+c_1^{-1})\mathbb{E}_{\mathcal{C}_k}[\|\bsx_{k+1}-\bsx_{k}\|^2]\nonumber\\
			&\quad+(1-c_1-2c_1^2)\|\bsx_{k}-\bsy_{k}\|^2,
			\label{nonconvex:xminush_compress_B_k}
		\end{align}
		where the first and second equalities hold due to \eqref{nonconvex:kia-algo-dc-compact-a} and \eqref{nonconvex:kia-algo-dc-compact-q}, respectively; the first inequality holds due to \eqref{nonconvex:lemma:inequality-arithmetic-equ} and $c_1>0$; the second inequality holds due to \eqref{nonconvex:lemma:inequality-arithmetic-equ} and $\omega r\in(0,1]$; 
		and the last inequality holds due to  $\bsx_k$ and $\bsy_k$ being independent of $\mathcal{C}_k$, along with \eqref{nonconvex:ass:compression_equ_scaling}.
		Then, taking the expectation with respect to $\mathcal{A}_k$ on both sides of \eqref{nonconvex:xminush_compress_B_k}, we have
		\begin{align}
			&\mathbb{E}_{\mathcal{A}_k}[\|\bsx_{k+1}-\bsy_{k+1}\|^2]\nonumber\\
			&\le(1+c_1^{-1})\mathbb{E}_{\mathcal{A}_k}[\|\bsx_{k+1}-\bsx_{k}\|^2]\nonumber\\
			&\quad+(1-c_1-2c_1^2)\|\bsx_{k}-\bsy_{k}\|^2.
			\label{nonconvex:xminush_compress}
		\end{align}

		For the first term on the right-hand side of \eqref{nonconvex:xminush_compress}, we have
		\begin{align}
			&\|\bsx_{k+1}-\bsx_k\|^2=\alpha_k^2\|\beta_k\bsL\hat{\bsx}_k+\gamma_k\bsv_k+\bsg_k^z\|^2\nonumber\\
			&=\alpha_k^2\|\beta_k\bsL(\hat{\bsx}_k-\bsx_k)+\beta_k\bsL\bsx_k+\gamma_k\bsv_k
			+\bsg_k^0+\bsg_k^z-\bsg_k^0\|^2\nonumber\\
			&\le4\alpha_k^2\big(\beta_k^2\|\hat{\bsx}_k-\bsx_k\|^2_{\bsL^2}+\beta_k^2\|\bsx_k\|^2_{\bsL^2}+\|\gamma_k\bsv_k+\bsg_k^0\|^2 \nonumber\\
			&\qquad  +\|\bsg_k^z-\bsg_k^0\|^2\big),
		\end{align}
		where the first equality holds due to \eqref{nonconvex:kia-algo-dc-compact-x}; and the first inequality holds due to \eqref{nonconvex:lemma:inequality-arithmetic-equ} and $\rho(\bsE)=1$.
		Then, taking expectation with respect to $\mathcal{A}_k$, from the independence of $\mathcal{B}_k$ and $\mathcal{C}_k$, $\mathcal{A}_k = \mathcal{B}_k \cup \mathcal{C}_k$,  \eqref{nonconvex:KL-L-eq2}, \eqref{nonconvex:lemma-eq5} and \eqref{zerosg:rand-grad-esti6}, we have
		\begin{align}\label{nonconvex:xkoneminusx}
			&\mathbb{E}_{\mathcal{A}_k}[\|\bsx_{k+1}-\bsx_{k}\|^2]
			\le 4\rho^2(L)\alpha_k^2\beta_k^2 \mathbb{E}_{\mathcal{C}_k}[\|\hat{\bsx}_k-\bsx_k\|_\bsE^2] 
			\nonumber\\
            &\quad
			+\|\bsx_k\|^2_{4\alpha_k^2(\beta_k^2\rho^2(L)+4 \ell ^2)\bsE}+ \Big\|\bsv_k+\frac{1}{\gamma_k}\bsg_k^0\Big\|^2_{4\rho(L)\alpha_k^2\gamma_k^2\bsF} 
			\nonumber\\
			&\quad+ 16n \ell ^2\alpha_k^2\mu_k^2 
			+ 8\alpha_k^2\mathbb{E}_{\mathcal{B}_k}[\|\bsg_k^z\|^2]. 	
		\end{align}	
	
		Then, from \eqref{nonconvex:xminush_compress} and \eqref{nonconvex:xkoneminusx}, we have
			\begin{align}
				&\mathbb{E}_{\mathcal{A}_k}[e_{5,k+1}]
					\le e_{5,k}-c_2\|\bsx_{k}-\bsy_{k}\|^2
				\nonumber\\
				&\quad +4(1+c_1^{-1})\rho^2(L)\alpha_k^2\beta_k^2 \mathbb{E}_{\mathcal{C}_k}[\|\bsx_k-\hat{\bsx}_k\|_\bsE^2] \nonumber\\
				&\quad + \|\bsx_k\|^2_{4(1+c_1^{-1})\alpha^2_k(\beta_k^2\rho^2(L)+4 \ell ^2)\bsE}\nonumber\\
				&\quad + 16(1+c_1^{-1})n \ell ^2\alpha_k^2\mu_k^2 + 8(1+c_1^{-1})\alpha_k^2\mathbb{E}_{\mathcal{B}_k}[\|\bsg_k^z\|^2]\nonumber\\
				&\quad + \Big\|\bsv_k+\frac{1}{\gamma_k}\bsg_k^0\Big\|^2_{4(1+c_1^{-1})\rho(L)\alpha_k^2\gamma_k^2\bsF}.
			\end{align}

		\subsubsection{Proof of Lemma~\ref{zerosg:lemma:sg:L1}}
		For the purpose of analyzing the Lyapunov function $\mathcal{L}_{1,k}$, we divide it into two parts: $\mathcal{L}_{2,k} = \sum_{i=1}^{3} e_{i,k} + e_{5,k}$ and $\mathcal{L}_{3,k} = e_{4,k}$.		 
		We then prove Lemma~\ref{zerosg:lemma:sg}, which serves as an extended version of Lemma~\ref{zerosg:lemma:sg:L1}.

		  
		\begin{lemmap}{\ref{zerosg:lemma:sg:L1}$'$}\label{zerosg:lemma:sg}
				Suppose  Assumptions~\ref{nonconvex:ass:compression}--\ref{zerosg:ass:fig} hold, $\{\gamma_k\}$ is non-decreasing, $\beta_k/\gamma_k=\epsilon_1$, $\alpha_k\gamma_k=\epsilon_2$, $\epsilon_1>\kappa_1$, $\epsilon_2>0$, and $\gamma_k\ge\varepsilon_{0}$, 
				then the Lyapunov function satisfies:
				\begin{subequations}
				\begin{align}
				&\mathbb{E}_{\mathcal{A}_k}[\mathcal{L}_{1,k+1}]  \nonumber\\
				&\quad \le \mathcal{L}_{1,k} -\|\bsx_k\|^2_{(2 a_1-\varepsilon_5\Delta_k-b_{1,k})\bsE}
				-\Big\|\bsv_k+\frac{1}{\gamma_k}\bsg_{k}^0\Big\|^2_{b_{2,k}\bsF}\nonumber\\
				&\quad - \alpha_k\Big(\frac{1}{4}-(1+\eta_2^2)(b_{3,k}+8p(1+\eta_1^2)b_{4,k})\alpha_k\Big)\|\bar{\bsg}^0_{k}\|^2 \nonumber\\
				&\quad +2pn\sigma^2_1b_{4,k}\alpha_k^2 + n\sigma^2_2(b_{3,k}+4p(1+\eta_1^2)b_{4,k})\alpha_k^2 \nonumber\\
				&\quad +b_{5,k}\alpha_k\mu_k^2	- (c_2-\delta_0\varepsilon_2-\delta_0b_{6,k}\Delta_k)\|\bsx_{k}-\bsy_{k}\|^2, 
				\nonumber\\	
					\label{zerosg:sgproof-vkLya}
					\\
				&\mathbb{E}_{\mathcal{A}_k}[\mathcal{L}_{2,k+1}]  \nonumber\\
				&\le  \mathcal{L}_{2,k} -\|\bsx_k\|^2_{(2 a_1-\varepsilon_5\Delta_k-b_{1,k})\bsE}
				-\Big\|\bsv_k+\frac{1}{\gamma_k}\bsg_{k}^0\Big\|^2_{b_{2,k}\bsF}\nonumber\\
				&\quad  + \Big(  b_{7,k} + (1+\eta_2^2)\big(b_{3,k}+8p(1+\eta_1^2)b_{4,k}\big)\Big)\alpha_k^2   \|\bar{\bsg}^0_{k}\|^2 \nonumber\\
				&\quad + 2pn\sigma^2_1b_{4,k}\alpha_k^2 + n\sigma^2_2\big(b_{3,k}+4p(1+\eta_1^2)b_{4,k}\big)\alpha_k^2 \nonumber\\
				&\quad 	+ b_{5,k}^{\prime}\alpha_k\mu_k^2 
				 - (c_2-\delta_0\varepsilon_2-\delta_0b_{6,k}\Delta_k)\|\bsx_{k}-\bsy_{k}\|^2,
				\label{zerosg:sgproof-vkLya-bounded}\\
				&\mathbb{E}_{\mathcal{A}_k}[\mathcal{L}_{3,k+1}]\le \mathcal{L}_{3,k}+\|\bsx_k\|^2_{2\alpha_k \ell ^2\bsE}
					-\frac{3}{16}\alpha_k\|\bar{\bsg}_{k}^0\|^2 \nonumber\\
					&\quad +p\alpha_k^2 a_7+(n+p) \ell ^2\alpha_k\mu^2_k, \label{Lyapunov 3}
				\end{align}
				where 
				\begin{align*}
					&b_{1,k}=8p(1+\eta_1^2)(1+\eta_2^2)\rho_2^{-2}(L) \ell ^4\frac{\alpha_k}{\gamma_k^2}\\
					&\quad +8p(1+\eta_1^2)(1+\eta_2^2)\Big(2\rho_2^{-2}(L)+\\
					&\qquad (\epsilon_1+1)\rho_2^{-1}(L)+2\Big) \ell ^4\frac{\alpha_k^2}{\gamma_k^2}\\
					&\quad +16p(1+\eta_1^2)(1+\eta_2^2)(\epsilon_1+1)^2\rho_2^{-1}(L) \ell ^4\frac{\alpha_k}{\gamma_k^3}\\
					&\quad+\Big(\frac{1}{2}+ \ell ^2\Big)\alpha_k\Delta_k + \Big(5+32p(1+\eta_1^2)(1+\eta_2^2)\\
					&\qquad +24p(1+\eta_1^2)(1+\eta_2^2)\varepsilon_8 \ell ^2\Big) \ell ^2\alpha_k^2\Delta_k\\
					&\quad+ 24p(1+\eta_1^2)(1+\eta_2^2)\varepsilon_1 \ell ^4\alpha_k^2\Delta_k^2\\
					&\quad+ 8p(1+\eta_1^2)(1+\eta_2^2)\Big((\epsilon_1+1)\rho_2^{-1}(L)+1\Big) \ell ^4\frac{\alpha_k^2\Delta_k}{\gamma_k^2}\\
					&\quad+ 16p(1+\eta_1^2)(1+\eta_2^2)(\epsilon_1+1)^2\rho_2^{-1}(L) \ell ^4\frac{\alpha_k\Delta_k}{\gamma_k^3}
					,\\	
				&b_{2,k}=2 a_2 - \frac{1}{2}\Delta_k(\rho(L)\epsilon_1\epsilon_2+3\rho(L)\epsilon_2^2+\epsilon_1+\epsilon_2+1),\\
				&b_{3,k}=\frac{3}{2}\varepsilon_8\frac{\Delta_k}{\alpha_k^2}
					+\frac{3}{2}\varepsilon_1\frac{\Delta_k^2}{\alpha_k^2},\\
				&b_{4,k}=6 + 16(1+c_1^{-1})+  \ell  + \frac{\rho_2^{-2}(L) \ell ^2}{\epsilon_2}\frac{1}{\gamma_k}\\
					&\quad +(2\rho_2^{-2}(L)+\frac{2(\epsilon_1+1)^2+\epsilon_1\epsilon_2+\epsilon_2}{\epsilon_2}\rho_2^{-1}(L)+2) \ell ^2\frac{1}{\gamma_k^2}\\
					&\quad +(4+3\varepsilon_8 \ell ^2)\Delta_k +3\varepsilon_1 \ell ^2\Delta_k ^2\\
					&\quad +(\frac{2(\epsilon_1+1)^2+\epsilon_1\epsilon_2+\epsilon_2}{\epsilon_2}\rho_2^{-1}(L)+1) \ell ^2\frac{\Delta_k}{\gamma_k^2},\\
				&b_{5,k}= n \ell ^2\Big(\frac{5}{2}+\big(8+16(1+c_1^{-1})+\frac{1}{4}p^2b_{4,k}\big)\alpha_k\\
					&\quad +\Delta_k+5\alpha_k\Delta_k\Big),\\
				&b_{6,k}=\epsilon_1^2\epsilon_2^2+2\epsilon_1\epsilon_2^2+\epsilon_2^2+\epsilon_1\epsilon_2\\
					&\quad +(\epsilon_1\epsilon_2^2+\epsilon_2^2+\frac{1}{2}\epsilon_1\epsilon_2+\epsilon_2)\rho(L) + 3\epsilon_1^2\epsilon_2^2\rho^2(L),\\
				&b_{5,k}^{\prime} = b_{5,k} + 2n \ell^2 b_{7,k} \frac{\epsilon_2}{\varepsilon_0},\\
				& a_2=\frac{1}{16}\epsilon_2-\big(\frac{5}{4}+2(1+c_1^{-1})\big)\rho(L)\epsilon_2^2,\\
				& a_7=2(\sigma^2_1+2(1+\eta_1^2)\sigma^2_2) \ell.
				\end{align*}
			\end{subequations}
			\end{lemmap}

		\begin{proof}

		We first denote the notations to be used next: 		
		\begin{align*}
			&\bsM_{1,k}=\frac{\beta_k}{2}\bsL-\frac{9\gamma_k}{4}\bsE-(1+\frac{5}{2} \ell ^2)\bsE,\\	
			&\bsM^0_{2,k}= 3\beta_k^2\bsL^2 - (\beta_k\gamma_k-\gamma_k^2)\bsL + \Big(4(1+c_1^{-1})\rho^2(L)\beta_k^2\\
			 	&\quad + 2\beta_k^2 +2\beta_k\gamma_k +3\gamma_k^2 + 8 \ell ^2+16(1+c_1^{-1}) \ell ^2\Big)\bsE,\\
			&\bsM_{2,k}=\bsM^0_{2,k}
				\\
				&\quad +8p(1+\eta_1^2)(1+\eta_2^2)\big(6+ 16(1+c_1^{-1})+ \ell \big) \ell ^2\bsE,\\
			&\bsM_{3}=3\epsilon_1^2\epsilon_2^2\bsL^2 + (\epsilon_1\epsilon_2^2+\epsilon_2^2-\frac{1}{2}\epsilon_1\epsilon_2)\bsL \\
			 	&\quad + (\epsilon_1^2\epsilon_2^2+2\epsilon_1\epsilon_2^2+\epsilon_2^2+\epsilon_1\epsilon_2+\frac{1}{2})\bsE,\\	
			&\bsM_{4,k}=\frac{\rho(L)}{2}(\beta_k+2\gamma_k){\bf I}_{np} + 2\gamma_k\bsE,\\
			&\bsM_{5,k}=3\beta_k^2\bsL^2-(\beta_k\gamma_k-\gamma_k^2)\bsL + \Big(4(1+c_1^{-1})\rho^2(L)\beta_k^2\\
				&\quad + 2\beta_k^2 +2\beta_k\gamma_k +3\gamma_k^2 
				\Big)\bsE,\\
			&\bsM_{6}=\frac{\rho(L)}{2}(\epsilon_1\epsilon_2+2\epsilon_2){\bf I}_{np} + 3\epsilon_1^2\epsilon_2^2\bsL^2 +(\epsilon_1\epsilon_2^2+\epsilon_2^2)\bsL \\
			 	&\quad + (\epsilon_1^2\epsilon_2^2+2\epsilon_1\epsilon_2^2+\epsilon_2^2+\epsilon_1\epsilon_2)\bsE,\\	
			&\bsM_{2,k}^{\prime} = \bsM_{2,k} + 2 b_{7,k} \ell^2 \bsE,\\
			&b^0_{2,k}=\alpha_k(\frac{1}{4}\gamma_k-\rho_2^{-1}(L))-\big(\frac{5}{2}+4(1+c_1^{-1})\big)\rho(L)\epsilon_2^2\\
				&\quad -\frac{1}{2}\Delta_k(\rho(L)\epsilon_1\epsilon_2+3\rho(L)\epsilon_2^2+\epsilon_1+\epsilon_2+1),\\
			&b^0_{4,k}=3 + 8(1+c_1^{-1}) +\frac{1}{2} \ell +\frac{\rho_2^{-2}(L) \ell ^2}{2\epsilon_2}\frac{1}{\gamma_k}\\
				&\quad +(\rho_2^{-2}(L)+\frac{2(\epsilon_1+1)^2+\epsilon_1\epsilon_2+\epsilon_2}{2\epsilon_2}\rho_2^{-1}(L)+1) \ell ^2\frac{1}{\gamma_k^2}\\
				&\quad +2\Delta_k + (\frac{2(\epsilon_1+1)^2+\epsilon_1\epsilon_2+\epsilon_2}{2\epsilon_2}\rho_2^{-1}(L)+\frac{1}{2}) \ell ^2\frac{\Delta_k}{\gamma_k^2}
				.\\
			\end{align*}

		Then, we start to show the relation between $\mathcal{L}_{i,k+1}$ and $\mathcal{L}_{i,k}$ for $i=1,2,3$.

		\noindent {\bf (i)} This step is to show the relation between $\mathcal{L}_{1,k+1}$ and $\mathcal{L}_{1,k}$. We have		
		\begin{align}
			&\mathbb{E}_{\mathcal{A}}[\mathcal{L}_{1,k+1}]\nonumber\\
			&\le  \mathcal{L}_{1,k}+\frac{1}{2}\Delta_k\|\bsx_k\|^2_{\bsE} \nonumber\\
			&\quad - (1+\Delta_k)\|\bsx_k\|^2_{\frac{\alpha_k\beta_k}{2}\bsL-\frac{\alpha_k}{2}\bsE-\alpha_k(1+5\alpha_k) \ell ^2\bsE}\nonumber\\
			&\quad+ (1+\Delta_k)\mathbb{E}_{\mathcal{C}_k}\Big[\|\hat{\bsx}_k-\bsx_k+\bsx_k\|^2_{\frac{3}{2}\alpha_k^2\beta_k^2\bsL^2}\Big] \nonumber\\
			&\quad + (1+\Delta_k)n \ell ^2\alpha_k(1+5\alpha_k)\mu^2_k\nonumber\\
			&\quad + (1+\Delta_k)\frac{\alpha_k}{2}(\beta_k+2\gamma_k)\rho(L)\mathbb{E}_{\mathcal{C}_k}[\|\bsx_k-\hat{\bsx}_k\|^2] \nonumber\\
			&\quad + 2(1+\Delta_k)\alpha^2_k\mathbb{E}_{\mathcal{B}_k}[\|\bsg_k^z\|^2]\nonumber\\
			&\quad + (1+\Delta_k)\Big\|\bsv_k+\frac{1}{\gamma_k}\bsg_k^0\Big\|^2_{\frac{6\alpha_k^2\gamma_k^2\rho(L)+\alpha_k\gamma_k}{4}\bsF} 
			\nonumber\\ 
			&\quad +\frac{1}{2}\varepsilon_1(\Delta_k+\Delta_k^2)\mathbb{E}_{\mathcal{A}_k}[\|\bsg_{k+1}^0\|^2] \nonumber\\
			&\quad+\mathbb{E}_{\mathcal{C}_k}\Big[\|\hat{\bsx}_k-\bsx_k+\bsx_k\|^2_{(1+\Delta_k)\frac{1}{2}\alpha_k^2(\beta_k\gamma_k+\gamma_k^2)\bsL }\Big] \nonumber\\
			&\quad+\mathbb{E}_{\mathcal{C}_k}\Big[\|\hat{\bsx}_k-\bsx_k+\bsx_k\|^2_{(1+\Delta_k)\frac{1}{2}\alpha_k^2(\beta_k+\gamma_k)^2\bsE}\Big] \nonumber\\
			&\quad +\Big\|\bsv_k+\frac{1}{\gamma_k}\bsg_{k}^0\Big\|^2_{\big(\alpha_k\frac{\gamma_k}{4}+\Delta_k(\frac{\beta_k+\gamma_k}{2\gamma_k}+\frac{\alpha_k\gamma_k}{4})\big)\bsF}\nonumber\\
			&\quad +(1+\Delta_k)(\epsilon_4\alpha_k+\epsilon_5\alpha_k^2) \ell ^2\mathbb{E}_{\mathcal{B}_k}\big[\|\bar{\bsg}_k^z\|^2\big] 
			\nonumber\\ 
			&\quad+\|\bsx_k\|^2_{\alpha_k(\frac{\gamma_k+2}{4}+\frac{1}{2} \ell ^2)\bsE+3\alpha_k^2 \ell ^2\bsE}\nonumber\\
			&\quad +\mathbb{E}_{\mathcal{C}_k}\Big[\|\hat{\bsx}_k-\bsx_k+\bsx_k\|^2_{\alpha_k\gamma_k\bsE
				+\alpha_k^2\big((\frac{1}{2}\beta_k^2+\gamma_k^2)\bsE-\beta_k\gamma_k\bsL\big)}\Big] \nonumber\\
			&\quad +\mathbb{E}_{\mathcal{C}_k}\Big[\|\hat{\bsx}_k-\bsx_k+\bsx_k\|^2_{\frac{1}{2}\Delta_k\alpha_k\beta_k\bsE} \Big]\nonumber\\
			&\quad -\Big\|\bsv_k+\frac{1}{\gamma_k}\bsg_{k}^0\Big\|^2_{\big(\alpha_k(\frac{3}{4}\gamma_k-\rho_2^{-1}(L))-\alpha_k^2\gamma_k^2\rho(L) -\frac{1}{2}\rho(L)\Delta_k\alpha_k\beta_k\big)\bsF}\nonumber\\
			&\quad +(\alpha_k\epsilon_6+\alpha_k^2\epsilon_7) \ell ^2\mathbb{E}_{\mathcal{B}_k}\big[\|\bar{\bsg}_k^z\|^2\big]  +n \ell ^2\alpha_k(\frac{1}{2}+3\alpha_k)\mu_k^2 \nonumber\\
			&\quad + \alpha^2_k\mathbb{E}_{\mathcal{B}_k}[\|\bsg_k^z\|^2] +\frac{\alpha_k}{\gamma_k}\rho_2^{-1}(L)\|\bar{\bsg}_k^\mu\|^2 \nonumber\\
			&\quad+\frac{1}{2}\rho_2^{-2}(L)\Delta_k\mathbb{E}_{\mathcal{A}_k}[\|\bsg_{k+1}^0\|^2]
			\nonumber\\ 
			&\quad- \frac{1}{4}\alpha_k\|\bar{\bsg}^\mu_{k}\|^2 + \|\bsx_k\|^2_{\alpha_k \ell ^2\bsE}\nonumber\\
			&\quad + n \ell ^2\alpha_k\mu^2_k-\frac{1}{4}\alpha_k\|\bar{\bsg}_{k}^0\|^2
			+\frac{1}{2}\alpha^2_k \ell \mathbb{E}_{\mathcal{B}_k}[\|\bar{\bsg}^z_{k}\|^2]
			\nonumber\\ 
			&\quad - c_2\|\bsx_{k}-\bsy_{k}\|^2 \nonumber\\
			&\quad +4(1+c_1^{-1})\rho^2(L)\alpha_k^2\beta_k^2 \mathbb{E}_{\mathcal{C}_k}[\|\bsx_k-\hat{\bsx}_k\|_\bsE^2] \nonumber\\
			&\quad + \|\bsx_k\|^2_{4(1+c_1^{-1})\alpha^2_k(\beta_k^2\rho^2(L)+4 \ell ^2)\bsE}\nonumber\\
			&\quad + 16(1+c_1^{-1})n \ell ^2\alpha_k^2\mu_k^2 + 8(1+c_1^{-1})\alpha_k^2\mathbb{E}_{\mathcal{B}_k}[\|\bsg_k^z\|^2]\nonumber\\
			&\quad + \Big\|\bsv_k+\frac{1}{\gamma_k}\bsg_k^0\Big\|^2_{4(1+c_1^{-1})\rho(L)\alpha_k^2\gamma_k^2\bsF}
			\nonumber\\ 
			&\le \mathcal{L}_{1,k} - \|\bsx_k\|^2_{\alpha_k\bsM_{1,k}-\alpha_k^2\bsM^0_{2,k}-\Delta_k\bsM_{3}}\nonumber\\
			&\quad -\|\bsx_k\|^2_{((\frac{1}{2}+ \ell ^2)\alpha_k\Delta_k+5 \ell ^2\alpha_k^2\Delta_k)\bsE} \nonumber\\
			&\quad + \mathbb{E}_{\mathcal{C}_k}\big[\|\bsx_k-\hat{\bsx}_k\|^2_{\alpha_k\bsM_{4,k}+\alpha_k^2\bsM_{5,k}
			 + \Delta_k\bsM_6}\big] \nonumber\\
			&\quad - \Big\|\bsv_k+\frac{1}{\gamma_k}\bsg_{k}^0\Big\|^2_{b^0_{2,k}\bsF}
			  -\frac{1}{4}\alpha_k\|\bar{\bsg}_{k}^0\|^2+b^0_{4,k}\alpha^2_k 
			  \mathbb{E}_{\mathcal{B}_k}[\|\bsg_k^z\|^2] \nonumber\\
			&\quad + \frac{1}{2}(\varepsilon_8\Delta_k+\varepsilon_1\Delta_k^2)
			\mathbb{E}_{\mathcal{A}_k}[\|\bsg_{k+1}^0\|^2] - c_2\|\bsx_{k}-\bsy_{k}\|^2 \nonumber\\
			&\quad + n \ell ^2\Big(\frac{5}{2}+\big(8+16(1+c_1^{-1})\big)\alpha_k+\Delta_k+5\alpha_k\Delta_k\Big)\alpha_k\mu_k^2 \nonumber\\
				&\quad -\alpha_k\Big(\frac{1}{4}-\frac{\rho_2^{-1}(L)}{\gamma_k}\Big) \|\bar{\bsg}_{k}^{\mu}\|^2 			\label{zerosg:vkLya-1}\\
			&\le  \mathcal{L}_{1,k} - \|\bsx_k\|^2_{\alpha_k\bsM_{1,k}-\alpha_k^2\bsM_{2,k}-\Delta_k\bsM_{3}
			-b_{1,k}\bsE} \nonumber\\
			&\quad + \mathbb{E}_{\mathcal{C}_k}\big[\|\bsx_k-\hat{\bsx}_k\|^2_{\alpha_k\bsM_{4,k}+\alpha_k^2\bsM_{5,k} + \Delta_k\bsM_6}\big] \nonumber\\
			&\quad - \Big\|\bsv_k+\frac{1}{\gamma_k}\bsg_{k}^0\Big\|^2_{b^0_{2,k}\bsF} - c_2\|\bsx_{k}-\bsy_{k}\|^2 \nonumber\\
			&\quad -\alpha_k\Big(\frac{1}{4}-(1+\eta_2^2)(b_{3,k}+8p(1+\eta_1^2)b_{4,k})\alpha_k\Big)\|\bar{\bsg}^0_{k}\|^2 \nonumber\\
			&\quad + 2pn\sigma^2_1b_{4,k}\alpha_k^2+n\sigma^2_2(b_{3,k}+4p(1+\eta_1^2)b_{4,k})\alpha_k^2 \nonumber\\
				&\quad +b_{5,k}\alpha_k\mu_k^2  -\alpha_k\Big(\frac{1}{4}-\frac{\rho_2^{-1}(L)}{\gamma_k}\Big) \|\bar{\bsg}_{k}^{\mu}\|^2\label{zerosg:vkLya}\\
			&\le  \mathcal{L}_{1,k} - \|\bsx_k\|^2_{\alpha_k\bsM_{1,k}-\alpha_k^2\bsM_{2,k}-\Delta_k\bsM_{3}-b_{1,k}\bsE} \nonumber\\
			&\quad - \Big\|\bsv_k+\frac{1}{\gamma_k}\bsg_{k}^0\Big\|^2_{b^0_{2,k}\bsF} \nonumber\\
			&\quad -\alpha_k\Big(\frac{1}{4}-(1+\eta_2^2)(b_{3,k}+8p(1+\eta_1^2)b_{4,k})\alpha_k\Big)\|\bar{\bsg}^0_{k}\|^2 \nonumber\\
			&\quad + 2pn\sigma^2_1b_{4,k}\alpha_k^2+n\sigma^2_2(b_{3,k}+4p(1+\eta_1^2)b_{4,k})\alpha_k^2 \nonumber\\
			&\quad +b_{5,k}\alpha_k\mu_k^2 - (c_2-\delta_0\varepsilon_2-\delta_0b_{6,k}\Delta_k)\|\bsx_{k}-\bsy_{k}\|^2 \nonumber\\		
			&\quad -\alpha_k\Big(\frac{1}{4}-\frac{\rho_2^{-1}(L)}{\gamma_k}\Big) \|\bar{\bsg}_{k}^{\mu}\|^2   \label{nonconvex:vkLya_xhat}
		\end{align}
		where the first inequality holds due to \eqref{nonconvex:v1k}, \eqref{zerosg:v2k}, \eqref{zerosg:v3k}, \eqref{zerosg:v4k} and \eqref{nonconvex:xminush};
		  the second inequality holds due to $\|\hat{\bsx}_k\|^2=\|\hat{\bsx}_k-\bsx_k+\bsx_k\|^2\le2\|\hat{\bsx}_k-\bsx_k\|^2+2\|\bsx_k\|^2$, the independence between $\bsx_k$ and $\mathcal{C}_k$, $\|\bar{\bsg}^z_{k}\|^2\le\|\bsg^z_{k}\|^2$, $\beta_k=\epsilon_1\gamma_k$ and  $\alpha_k=\frac{\epsilon_2}{\gamma_k}$; 
		   the tird inequality hods due to \eqref{zerosg:rand-grad-esti2} and \eqref{zerosg:rand-grad-esti4};
		    the last inequality hods due to \eqref{nonconvex:xminusxhat}, the independence of $\bsx_k$ and $\bsy_k$ from $\mathcal{C}_k$, \eqref{nonconvex:KL-L-eq2}, $\rho(\bsE)=1$, $\gamma_k\ge\varepsilon_{0}\ge1$ and $\alpha_k\le\epsilon_2$.

		Next, we scale and bound some coefficients to prove the Lemma~\ref{zerosg:lemma:sg}.

		From \eqref{nonconvex:KL-L-eq2}, $\beta_k=\epsilon_1\gamma_k$, $\epsilon_1>\kappa_1\ge\frac{13}{2\rho_2(L)}$, $\gamma_k\ge\varepsilon_{0}\ge1+\frac{5}{2} \ell ^2$, and $\alpha_k=\frac{\epsilon_2}{\gamma_k}$, we have
		\begin{align}\label{zerosg:m1-rand-pd}
		\alpha_k\bsM_{1,k}\ge\varepsilon_3\epsilon_2\bsE.
		\end{align}
		From \eqref{nonconvex:KL-L-eq2}, $\beta_k=\epsilon_1\gamma_k$, $\gamma_k\ge\varepsilon_{0}\ge\big(8+16(1+c_1^{-1})+8p(1+\eta_1^2)(1+\eta_2^2)(6+ 16(1+c_1^{-1})+ \ell )\big)^\frac{1}{2} \ell $, and $\alpha_k=\frac{\epsilon_2}{\gamma_k}$, we have
		\begin{align}\label{zerosg:m2-rand-pd}
		\alpha_k^2\bsM_{2,k}\le\varepsilon_4\epsilon_2^2\bsE.
		\end{align}
		From \eqref{nonconvex:KL-L-eq2}, $\beta_k=\epsilon_1\gamma_k$,  and $\alpha_k=\frac{\epsilon_2}{\gamma_k}$, we have
		\begin{align}\label{zerosg:m3-rand-pd}
		\bsM_{3}\le\varepsilon_5\bsE.
		\end{align}
		From $\gamma_k\ge\varepsilon_{0}\ge 12\rho_2^{-1}(L)$ and $\alpha_k=\frac{\epsilon_2}{\gamma_k}$, we have
		\begin{align}\label{zerosg:vkLya-b1}
		b^0_{2,k}\ge&b_{2,k}.
		\end{align}
		From $\gamma_k\ge\varepsilon_{0}\ge 12\rho_2^{-1}(L)$, we have
		\begin{align}\label{zerosg:vkLya-b2}
			\frac{1}{4}-\frac{\rho_2^{-1}(L)}{\gamma_k}  \ge 0
		\end{align}

		Finally, from \eqref{zerosg:vkLya}--\eqref{zerosg:vkLya-b2}, we know that \eqref{zerosg:sgproof-vkLya} holds.

		\noindent {\bf (ii)} This step is to show the relation between $\mathcal{L}_{2,k+1}$ and $\mathcal{L}_{2,k}$.		
		
		Similar to the derivation of~\eqref{zerosg:sgproof-vkLya}, 
				we obtain
		\begin{align}
			&\mathbb{E}_{\mathcal{A}_k}[\mathcal{L}_{2,k+1}]  \nonumber\\
			&\le  \mathcal{L}_{2,k} - \|\bsx_k\|^2_{\alpha_k\bsM_{1,k}-\alpha_k^2\bsM_{2,k}-\Delta_k\bsM_{3}-b_{1,k}\bsE} \nonumber\\
			&\quad - \Big\|\bsv_k+\frac{1}{\gamma_k}\bsg_{k}^0\Big\|^2_{b^0_{2,k}\bsF} 
				+ \frac{\alpha_k}{\gamma_k}\rho_2^{-1}(L)\|\bar{\bsg}_k^\mu - \bar{\bsg}_k^0 + \bar{\bsg}_k^0\|^2 \nonumber\\
			&\quad  + (1+\eta_2^2)(b_{3,k}+8p(1+\eta_1^2)b_{4,k})\alpha_k^2 \|\bar{\bsg}^0_{k}\|^2 \nonumber\\
			&\quad + 2pn\sigma^2_1b_{4,k}\alpha_k^2+n\sigma^2_2(b_{3,k}+4p(1+\eta_1^2)b_{4,k})\alpha_k^2 \nonumber\\
			&\quad +b_{5,k}\alpha_k\mu_k^2 - (c_2-\delta_0\varepsilon_2-\delta_0b_{6,k}\Delta_k)\|\bsx_{k}-\bsy_{k}\|^2 \nonumber\\
			&\le  \mathcal{L}_{2,k} - \|\bsx_k\|^2_{\alpha_k\bsM_{1,k}-\alpha_k^2\bsM_{2,k}^{\prime}-\Delta_k\bsM_{3}-b_{1,k}\bsE} \nonumber\\
			&\quad - \Big\|\bsv_k+\frac{1}{\gamma_k}\bsg_{k}^0\Big\|^2_{b^0_{2,k}\bsF} 
				+ b_{5,k}\alpha_k\mu_k^2 + 2n \ell^2 b_{7,k} \alpha_k^2 \mu_k^2 \nonumber\\
			&\quad  + \Big(  b_{7,k} + (1+\eta_2^2)\big(b_{3,k}+8p(1+\eta_1^2)b_{4,k}\big)\Big)\alpha_k^2   \|\bar{\bsg}^0_{k}\|^2 \nonumber\\
			&\quad + 2pn\sigma^2_1b_{4,k}\alpha_k^2 + n\sigma^2_2\big(b_{3,k}+4p(1+\eta_1^2)b_{4,k}\big)\alpha_k^2 \nonumber\\
			&\quad - (c_2-\delta_0\varepsilon_2-\delta_0b_{6,k}\Delta_k)\|\bsx_{k}-\bsy_{k}\|^2 \nonumber\\
			&\le  \mathcal{L}_{2,k} -\|\bsx_k\|^2_{(2 a_1-\varepsilon_5\Delta_k-b_{1,k})\bsE}
				-\Big\|\bsv_k+\frac{1}{\gamma_k}\bsg_{k}^0\Big\|^2_{b_{2,k}\bsF}\nonumber\\
			&\quad  + \Big(  b_{7,k} + (1+\eta_2^2)\big(b_{3,k}+8p(1+\eta_1^2)b_{4,k}\big)\Big)\alpha_k^2   \|\bar{\bsg}^0_{k}\|^2 \nonumber\\
			&\quad + 2pn\sigma^2_1b_{4,k}\alpha_k^2 + n\sigma^2_2\big(b_{3,k}+4p(1+\eta_1^2)b_{4,k}\big)\alpha_k^2 \nonumber\\
			&\quad 	+ b_{5,k}^{\prime}\alpha_k\mu_k^2 
			 - (c_2-\delta_0\varepsilon_2-\delta_0b_{6,k}\Delta_k)\|\bsx_{k}-\bsy_{k}\|^2,
		\end{align}
		where the first inequality holds in the same manner as~\eqref{nonconvex:vkLya_xhat};
		the second inequality hold due to $\|\bar{\bsg}_k^\mu - \bar{\bsg}_k^0 + \bar{\bsg}_k^0\|^2 \le 2\|  \bar{\bsg}_k^0\|^2 +2 \|\bar{\bsg}_k^\mu - \bar{\bsg}_k^0\| $, \eqref{zerosg:rand-grad-esti9}, and $\gamma_k=\epsilon_2/\alpha_k$;
		the last inequality holds in the same manner as~\eqref{zerosg:m1-rand-pd}--\eqref{zerosg:vkLya-b2} and $\alpha_k = \frac{\epsilon_2}{\gamma_k} \le \frac{\epsilon_2}{\varepsilon_0}$.

		\noindent {\bf (iii)} This step is to show the relation between $\mathcal{L}_{3,k+1}$ and $\mathcal{L}_{3,k}$.		
		From \eqref{zerosg:v4k}, \eqref{zerosg:rand-grad-esti5}, \eqref{zerosg:rand-grad-esti2}, and $\mathcal{L}_{3,k}=e_{4,k}$, we have
		\begin{align}
		&\mathbb{E}_{\mathcal{A}_k}[\mathcal{L}_{3,k+1}] \nonumber\\
		&\quad \le \mathcal{L}_{3,k}-\frac{1}{4}\alpha_k\|\bar{\bsg}^\mu_{k}\|^2
		+\|\bsx_k\|^2_{\alpha_k \ell ^2\bsE}+n \ell ^2\alpha_k\mu^2_k\nonumber\\
		&\quad-\frac{1}{4}\alpha_k\|\bar{\bsg}_{k}^0\|^2
		+\frac{1}{2}\alpha_k^2 \ell \Big(\frac{16p(1+\eta_1^2)(1+\eta_2^2)}{n}\|\bar{\bsg}_{k}^0\|^2 \nonumber\\
		&\quad 	+\frac{16p(1+\eta_1^2)(1+\eta_2^2)}{n} \ell ^2\|\bsx_{k}\|^2_{\bsE}\nonumber\\
		&\quad+4p\sigma^2_1+8p(1+\eta_1^2)\sigma^2_2+\frac{1}{2}p^2 \ell ^2\mu_k^2
		+\|\bar{\bsg}^\mu_{k}\|^2\Big).\label{zerosg:v4kspeed-diminishing}
		\end{align}

		From $\alpha_k=\frac{\epsilon_2}{\gamma_k}$, and $\gamma_k\ge \varepsilon_0\ge128p(1+\eta_1^2)(1+\eta_2^2)\epsilon_2 \ell $, we have
		\begin{subequations}
		\begin{align}
		&\frac{8p(1+\eta_1^2)(1+\eta_2^2)}{n}\alpha_k^2 \ell  \nonumber\\
			&\quad \le\frac{8p(1+\eta_1^2)(1+\eta_2^2)\epsilon_2}{\gamma_k}\alpha_k \ell 
			\ge\frac{1}{16}\alpha_k,\label{zerosg:v4kspeed-diminishing-1.1}\\
		&\frac{8p(1+\eta_1^2)(1+\eta_2^2)}{n}\alpha_k^2 \ell ^3
			\ge\alpha_k \ell ^2,\label{zerosg:v4kspeed-diminishing-1.2}\\
		&\frac{1}{2}\alpha_k^2 \ell \ge\frac{1}{16}\alpha_k,\label{zerosg:v4kspeed-diminishing-1.3}\\
		&\frac{1}{4}p^2\alpha_k^2 \ell ^3\ge p \ell ^2\alpha_k.\label{zerosg:v4kspeed-diminishing-1.4}
		\end{align}
		\end{subequations}
	
		Then, from \eqref{zerosg:v4kspeed-diminishing}--\eqref{zerosg:v4kspeed-diminishing-1.4}, we have \eqref{Lyapunov 3}.

	\end{proof}

\subsection{Proof of Theorem~\ref{zerosg:thm-sg-smT}}\label{zerosg:proof-thm-sg-smT}
In addition to the notations defined in Appendix~\ref{zerosg:proof-thm-random-pd-sm},
we also denote the following notations.
\begin{align*}
&\tilde{\kappa}_0(\epsilon_1,\epsilon_2)=\max\Big\{\varepsilon_{0},
	~\Big(\frac{p(1+\eta_1^2)(1+\eta_2^2)\tilde{\varepsilon}_9}{ a_1}\Big)^{\frac{1}{3}},\\
	&\quad ~64p(1+\eta_1^2)(1+\eta_2^2)\epsilon_2\tilde{\varepsilon}_{12}\Big\},\\
&\kappa_2(\epsilon_1)=\min\Big\{\frac{\varepsilon_3}{\varepsilon_4},~\frac{\rho^{-1}(L)}{20+32(1+c_1^{-1})}, \\
	&\quad ~\frac{\sqrt{\varepsilon_{13}^2+4\varepsilon_{14}c_2}-\varepsilon_{13}}{2\varepsilon_{14}},~8\Big\},\\
&\varepsilon_6=\max\Big\{\frac{1+\epsilon_1\rho_2(L)}{2},
~\frac{1+\epsilon_1}{2}+\frac{1}{2\epsilon_1\rho_2^2(L)}\Big\},\\
&\varepsilon_7=\frac{\epsilon_1\rho_2(L)-1}{2\epsilon_1\rho_2(L)},\\			
&\tilde{\varepsilon}_9
	=8\big(\rho_2^{-2}(L)+ 2(\epsilon_1+1)^2\rho_2^{-1}(L)\big) \ell ^4\epsilon_2 \\
	&\quad +8\big(2\rho_2^{-2}(L)+(\epsilon_1+1)\rho_2^{-1}(L)+2\big) \ell ^4\epsilon_2^2,\\
&\tilde{\varepsilon}_{12}=6+16(1+c_1^{-1})+ \ell +\frac{\rho_2^{-2}(L)}{\epsilon_2} \ell ^2\\
	&\quad +\big(2\rho_2^{-2}(L)+\frac{2(\epsilon_1+1)^2+2\epsilon_1\epsilon_2+\epsilon_2}{\epsilon_2}\rho_2^{-1}(L)+2\big) \ell ^2,\\
&\varepsilon_{13}=\frac{1}{2}\rho(L)\delta_0\epsilon_1+2\delta_0+\rho(L)\delta_0,\\
&\varepsilon_{14}=\big(3\rho^2(L)+4(1+c_1^{-1})\rho^2(L)+2\big)\delta_0\epsilon_1^2\\
	& \quad +\big(2-\rho_2(L)\big)\delta_0\epsilon_1+\big(3+\rho(L) \ell ^2\big)\delta_0,\\			
& a_2=\frac{1}{16}\epsilon_2-\big(\frac{5}{4}+2(1+c_1^{-1})\big)\rho(L)\epsilon_2^2,\\
& a_3=\frac{1}{2}c_2-\frac{1}{2}\delta_0\varepsilon_2 \\
&\tilde{ a}_4=2(\sigma^2_1+2(1+\eta_1^2)\sigma^2_2)\tilde{\varepsilon}_{12},\\
&\tilde{ a}_5= \ell ^2\Big(\frac{1}{256(1+\eta_1^2)(1+\eta_2^2)}+	\frac{(24+16c_1^{-1})\epsilon_2+\frac{5}{2}}{p}\Big),\\
&\tilde{ a}_6=8(1+\eta_1^2)(1+\eta_2^2)\tilde{\varepsilon}_{12} + b_{7,k}/p,\\
&\tilde{ a}_5^{\prime}= \ell ^2\Big(\frac{1}{256(1+\eta_1^2)(1+\eta_2^2)}+	\frac{(24+16c_1^{-1})\epsilon_2+\frac{5}{2}}{p}\\
	&\quad + 2b_{7,k}\frac{\epsilon_2}{\varepsilon_0p}\Big),\\
&\tilde{d}_1=\frac{1}{\varepsilon_6}\min\{a_1, 2a_2, 2a_3\}.
\end{align*}

To prove Theorem~\ref{zerosg:thm-sg-smT}, the following lemma (results of Lemma \ref{zerosg:lemma:sg} under fixed algorithm parameters) is used.
\begin{lemma}\label{zerosg:lemma:sg2-T}
Suppose Assumptions~\ref{nonconvex:ass:compression}--\ref{zerosg:ass:fig} hold, $\beta_k=\beta=\epsilon_1\gamma$, $\gamma_k=\gamma$, and $\alpha_k=\alpha=\frac{\epsilon_2}{\gamma}$, where
$\epsilon_1>\kappa_1$, and $\epsilon_2\in(0,\kappa_2(\epsilon_1))$, and $\gamma\ge\tilde{\kappa}_0(\epsilon_1,\epsilon_2)$ are constants. Let $\{\bsx_k\}$ be the sequence generated by Algorithm~\ref{nonconvex:algorithm-pdgd}, then
\begin{subequations}
\begin{align}
&\mathbb{E}_{\mathcal{A}_k}[\mathcal{L}_{1,k+1}]\nonumber\\
&\quad\le   \mathcal{L}_{1,k}- a_1\|\bsx_k\|^2_{\bsE}
-2 a_2\Big\|\bsv_k+\frac{1}{\gamma}\bsg_{k}^0\Big\|^2_{\bsF} -\frac{1}{8}\alpha\|\bar{\bsg}^0_{k}\|^2\nonumber\\
&\quad-2 a_3\|\bsx_{k}-\bsy_{k}\|^2+pn\tilde{ a}_4\alpha^2+pn\tilde{ a}_5\alpha\mu_k^2,
\label{zerosg:sgproof-vkLya2T}\\
&\mathbb{E}_{\mathcal{A}_k}[\mathcal{L}_{2,k+1}]\nonumber\\
&\quad \le   \mathcal{L}_{2,k}- a_1\|\bsx_k\|^2_{\bsE}-2 a_2\Big\|\bsv_k+\frac{1}{\gamma}\bsg_{k}^0\Big\|^2_{\bsF} + p\tilde{ a}_6\alpha^2\|\bar{\bsg}^0_{k}\|^2\nonumber\\
&\quad -2 a_3\|\bsx_{k}-\bsy_{k}\|^2+pn\tilde{ a}_4\alpha^2+pn\tilde{ a}_5^{\prime}\alpha\mu_k^2,
\label{zerosg:sgproof-vkLya2T-bounded}\\
&\mathbb{E}_{\mathcal{A}_k}[\mathcal{L}_{3,k+1}]
\le  \mathcal{L}_{3,k}+\|\bsx_k\|^2_{2\alpha  \ell ^2\bsE}-\frac{1}{8}\alpha\|\bar{\bsg}_{k}^0\|^2\nonumber\\
&\quad +p a_7\alpha^2
+(n+p) \ell ^2\alpha\mu^2_k.\label{zerosg:v4kspeed}
\end{align}
\end{subequations}
\end{lemma}
\begin{proof}
 {\bf (i)}Substituting $\beta_k=\beta=\epsilon_1\gamma$, $\gamma_k=\gamma$, $\alpha_k=\alpha=\frac{\epsilon_2}{\gamma}$, and $\Delta_k=0$ into \eqref{nonconvex:v1k}, \eqref{zerosg:v2k}, \eqref{zerosg:v3k}, \eqref{zerosg:v4k} and \eqref{nonconvex:xminush}, similar to the way to get \eqref{nonconvex:vkLya_xhat}, we have
\begin{align}
&\mathbb{E}_{\mathcal{A}_k}[\mathcal{L}_{1,k+1}]\nonumber\\
&\quad \le \mathcal{L}_{1,k} - \|\bsx_k\|^2_{\alpha\tilde{\bsM}_{1}-\alpha^2\tilde{\bsM}_{2}-\tilde{b}_{1}\bsE} \nonumber\\
&\quad - \Big\|\bsv_k+\frac{1}{\gamma}\bsg_{k}^0\Big\|^2_{\tilde{b}^0_{2}\bsF} \nonumber\\
&\quad -\alpha\Big(\frac{1}{4}-8p(1+\eta_1^2)(1+\eta_2^2)\tilde{b}_{4}\alpha\Big)\|\bar{\bsg}^0_{k}\|^2 \nonumber\\
&\quad +2pn(\sigma^2_1+2(1+\eta_1^2)\sigma^2_2)\tilde{b}_{4}\alpha^2\nonumber\\
&\quad +\tilde{b}_{5}\alpha\mu_k^2 - 2 a_3\|\bsx_{k}-\bsy_{k}\|^2,
\label{zerosg:vkLyaT}
\end{align}
where
\begin{align*}
	&\tilde{\bsM}_{1}=\frac{\beta}{2}\bsL-\frac{9\gamma}{4}\bsE-(1+\frac{5}{2} \ell ^2)\bsE,\\	
	&\tilde{\bsM}_{2}^0= 3\beta^2\bsL^2 - (\beta\gamma-\gamma^2)\bsL + \Big(4(1+c_1^{-1})\rho^2(L)\beta^2\\
		 &\quad + 2\beta^2 +2\beta\gamma +3\gamma^2 + 8 \ell ^2+16(1+c_1^{-1}) \ell ^2\Big)\bsE,\\
	&\tilde{\bsM}_{2}=\tilde{\bsM}_{2}^0
		\\
		&\quad +8p(1+\eta_1^2)(1+\eta_2^2)\big(6+ 16(1+c_1^{-1})+ \ell \big) \ell ^2\bsE,\\
	&\tilde{b}_{1}=8p(1+\eta_1^2)(1+\eta_2^2)\rho_2^{-2}(L) \ell ^4\frac{\alpha}{\gamma^2}\\
		&\quad +8p(1+\eta_1^2)(1+\eta_2^2)\Big(2\rho_2^{-2}(L)+\\
		&\qquad (\epsilon_1+1)\rho_2^{-1}(L)+2\Big) \ell ^4\frac{\alpha^2}{\gamma^2}\\
		&\quad +16p(1+\eta_1^2)(1+\eta_2^2)(\epsilon_1+1)^2\rho_2^{-1}(L) \ell ^4\frac{\alpha}{\gamma^3}
		,\\	
	&\tilde{b}^0_{2}=\alpha(\frac{1}{4}\gamma-\rho_2^{-1}(L))-\big(\frac{5}{2}+4(1+c_1^{-1})\big)\rho(L)\epsilon_2^2,\\
	&\tilde{b}_{4}=6 + 16(1+c_1^{-1})+  \ell  + \frac{\rho_2^{-2}(L) \ell ^2}{\epsilon_2}\frac{1}{\gamma}\\
		&\quad +(2\rho_2^{-2}(L)+\frac{2(\epsilon_1+1)^2+\epsilon_1\epsilon_2+\epsilon_2}{\epsilon_2}\rho_2^{-1}(L)+2) \ell ^2\frac{1}{\gamma^2},\\
	&\tilde{b}_{5}= n \ell ^2\Big(\frac{5}{2}+\big(8+16(1+c_1^{-1})+\frac{1}{4}p^2\tilde{b}_{4}\big)\alpha\Big).\\
\end{align*}

From \eqref{zerosg:vkLyaT}, similar to the way to get \eqref{zerosg:sgproof-vkLya}, we have
\begin{align}\label{zerosg:sgproof-vkLyaT}
&\mathbb{E}_{\mathcal{A}_k}[\mathcal{L}_{1,k+1}]\nonumber\\
&\quad \le  \mathcal{L}_{1,k}-\|\bsx_k\|^2_{(2 a_1-\tilde{b}_{1})\bsE}
-\|\bsv_k+\frac{1}{\gamma_k}\bsg_{k}^0\|^2_{\tilde{b}_2^0\bsF}\nonumber\\
&\quad-\alpha\Big(\frac{1}{4}-8p(1+\eta_1^2)(1+\eta_2^2)\tilde{b}_{4}\alpha\Big)\|\bar{\bsg}^0_{k}\|^2\nonumber\\
&\quad+2pn(\sigma^2_1+2(1+\eta_1^2)\sigma^2_2)\tilde{b}_{4}\alpha^2+\tilde{b}_{5}\alpha\mu_k^2,\nonumber\\
&\quad -2 a_3\|\bsx_{k}-\bsy_{k}\|^2.
\end{align}

From $\gamma\ge\tilde{\kappa}_0(\epsilon_1,\epsilon_2)\ge\varepsilon_0\ge1$ and $\gamma\ge\tilde{\kappa}_0(\epsilon_1,\epsilon_2)\ge\Big(\frac{p(1+\eta_1^2)(1+\eta_2^2)\tilde{\varepsilon}_9}{ a_1}\Big)^{\frac{1}{3}}$, we have
\begin{align}
	2 a_1-\tilde{b}_{1}\ge2 a_1-\frac{p(1+\eta_1^2)(1+\eta_2^2)\tilde{\varepsilon}_9}{\gamma^3}\ge a_1.
\end{align}

From $\gamma\ge\tilde{\kappa}_0(\epsilon_1,\epsilon_2)\ge\varepsilon_0\ge12\rho_2^{-1}(L)$, we have 
\begin{align}
	\tilde{b}^0_{2}\ge2 a_2.
\end{align}

From $\gamma\ge1$, we have
\begin{align}
   \tilde{b}_4\le\tilde{\varepsilon}_{12}.
\end{align}

From $\gamma\ge\tilde{\kappa}_0(\epsilon_1,\epsilon_2)\ge~64p(1+\eta_1^2)(1+\eta_2^2)\epsilon_2\tilde{\varepsilon}_{12}$, we have
\begin{align}\label{zerosg:b3kb4keta:1}
	&\frac{1}{4}-8p(1+\eta_1^2)(1+\eta_2^2)\tilde{b}_4\alpha\nonumber\\
	&\ge\frac{1}{4}-\frac{8p(1+\eta_1^2)(1+\eta_2^2)\epsilon_2\tilde{\varepsilon}_{12}}{\gamma}
	\ge\frac{1}{8},\\
	&\tilde{b}_4\alpha\le\frac{1}{64p(1+\eta_1^2)(1+\eta_2^2)}.\label{mid}
	\end{align}

	From \eqref{mid} and $\gamma\ge1$, we have 
\begin{align}
	\tilde{b}_5\le pn\tilde{a}_5.
\end{align}

From $\epsilon_1>\kappa_1\ge\frac{13}{2\rho_2(L)}$, we have
$\varepsilon_3>0$.
From $\rho(L)\epsilon_1\ge\rho_2(L)\epsilon_1\ge\frac{13}{2}$, we have $\varepsilon_4\ge7\rho^2(L)\epsilon_1^2-\rho(L)\epsilon_1>0$.
From $\varepsilon_3>0$, $\varepsilon_4>0$, and $\epsilon_2 < \kappa_2(\epsilon_1) \le\min\{\frac{\varepsilon_3}{\varepsilon_4},~\frac{\rho^{-1}(L)}{20+32(1+c_1^{-1})}\}$, we have
\begin{align}
 a_1>0~\text{and}~
 a_2>0.\label{zerosg:kappa4-6}
\end{align}

From $\epsilon_2 < \kappa_2(\epsilon_1) \le\frac{\sqrt{\varepsilon_{13}^2+4\varepsilon_{14}c_2}-\varepsilon_{13}}{2\varepsilon_{14}}$, we have
\begin{align}\label{zerosg:c1}
	 a_3=\frac{1}{2}(c_2-\delta_0\varepsilon_2)=\frac{1}{2}(c_2-\varepsilon_{13}\epsilon_2-\varepsilon_{14}\epsilon_2^2)>0.
\end{align}

From \eqref{zerosg:sgproof-vkLyaT}--\eqref{zerosg:kappa4-6}, we have \eqref{zerosg:sgproof-vkLya2T}.

\noindent {\bf (ii)}
Similarly, we obtain \eqref{zerosg:sgproof-vkLya2T-bounded}.

\noindent {\bf (iii)}
Noting that $\alpha_k=\alpha$, $\gamma\ge\tilde{\kappa}_0(\epsilon_1,\epsilon_2)\ge64p(1+\eta_1^2)(1+\eta_2^2)\epsilon_2\tilde{\varepsilon}_{12}
\ge64p(1+\eta_1^2)(1+\eta_2^2)\epsilon_2 \ell $, and $\alpha=\frac{\epsilon_2}{\gamma}$, similar to the way to get \eqref{zerosg:v4kspeed-diminishing-2}, we have \eqref{zerosg:v4kspeed}.
\end{proof}

We are now ready to prove Theorem~\ref{zerosg:thm-sg-smT}.

From $\gamma_k=\gamma=\frac{\epsilon_2\sqrt{pT}}{\sqrt{n}}$ and $T\ge  \frac{n(\tilde{\kappa}_0(\epsilon_1,\epsilon_2))^2}{p\epsilon_2^2}$, we have $\gamma\ge\tilde{\kappa}_0(\epsilon_1,\epsilon_2)$. Thus, all conditions needed in Lemma~\ref{zerosg:lemma:sg2-T} are satisfied. So \eqref{zerosg:sgproof-vkLya2T}--\eqref{zerosg:v4kspeed} hold.

Since $\alpha_k=\alpha=\frac{\sqrt{n}}{\sqrt{pT}}$ and $\mu_{i,k}\le\frac{\kappa_\mu\sqrt{p\alpha}}{\sqrt{n+p}}$ as stated in \eqref{zerosg:step:eta2-sm}, taking expectation with respect to $\mathcal{F}_T$, summing \eqref{zerosg:sgproof-vkLya2T} over $ k\in[0,T-1]$ yields
\begin{align}
&\frac{1}{nT}\sum_{k=0}^{T-1}\mathbb{E}[\|\bsx_k\|^2_{\bsE}]
\le\frac{1}{ a_1}\Big(\frac{{\mathcal{L}}_{1,0}}{nT}
+\frac{n(\tilde{ a}_4+\tilde{ a}_5\kappa_\mu^2)}{T}\Big).
\label{zerosg:thm-sg-sm-equ3.1p}
\end{align}
Similarly, from \eqref{zerosg:v4kspeed} and \eqref{zerosg:step:eta2-sm},  we have
\begin{align}\label{zerosg:thm-sg-sm-equ3p}
&\frac{1}{T}\sum_{k=0}^{T}\mathbb{E}[\|\nabla f(\bar{x}_k)\|^2]=\frac{1}{nT}\sum_{k=0}^{T}\mathbb{E}[\|\bar{\bsg}_{k}^0\|^2]\nonumber\\
&\le 8\Big(\frac{{\mathcal{L}}_{3,0}}{nT\alpha}
+\frac{2 \ell ^2}{nT}\sum_{k=0}^{T}\mathbb{E}[\|\bsx_k\|^2_{\bsE}]+\frac{p\alpha( a_7+\ell^2k_\mu^2)}{n}\Big).
\end{align}
Noting that $\alpha=\frac{\sqrt{n}}{\sqrt{pT}}$, from \eqref{zerosg:thm-sg-sm-equ3p} and \eqref{zerosg:thm-sg-sm-equ3.1p}, we have
\begin{align*}
&\frac{1}{T}\sum_{k=0}^{T-1}\mathbb{E}[\|\nabla f(\bar{x}_k)\|^2]\nonumber\\
&=8(f(\bar{x}_0)-f^*+ a_7+ \ell ^2\kappa_\mu^2)\frac{\sqrt{p}}{\sqrt{nT}}
+\mathcal{O}\Big(\frac{n}{T}\Big),
\end{align*}
which gives \eqref{zerosg:coro-sg-sm-equ3}.

Taking expectation with respect to $\mathcal{F}_T$, summing \eqref{zerosg:v4kspeed} over $ k\in[0,T]$, and using \eqref{zerosg:step:eta2-sm}  yield
\begin{align}\label{zerosg:thm-sg-sm-equ4p}
&n(\mathbb{E}[f(\bar{x}_{T})]-f^*)=\mathbb{E}[{\mathcal{L}}_{3,T}]\nonumber\\
&\le {\mathcal{L}}_{3,0}+\frac{2 \ell ^2\sqrt{n}}{\sqrt{pT}} \sum_{k=0}^{T-1}\mathbb{E}[\|\bsx_k\|^2_{\bsE}]+n a_7
+n \ell ^2\kappa_\mu^2.
\end{align}

Noting that ${\mathcal{L}}_{3,0}=\mathcal{O}(n)$ and $\frac{\sqrt{n}n}{\sqrt{pT}}\le1$ due to $T\ge \frac{n^3}{p}$, from \eqref{zerosg:thm-sg-sm-equ3.1p} and \eqref{zerosg:thm-sg-sm-equ4p}, we have \eqref{zerosg:coro-sg-sm-equ4}. Then, from \eqref{zerosg:rand-grad-smooth}, we know that there exists a constant $\tilde{c}_g>0$, such that
\begin{align}\label{zerosg:coro-sg-sm-gbark0}
\mathbb{E}[\|\bar{\bsg}^0_k\|^2]\le n\tilde{c}_g,~\forall k\in\mathbb{N}_0.
\end{align}

Denote	$\hat{\mathcal{L}}_{2,k}=\|\bsx_k\|^2_{\bsE}+\Big\|\bsv_k	+\frac{1}{\gamma}\bsg_k^0\Big\|^2_{\bsF}+\|\bsx_{k}-\bsy_{k}\|^2$.
We have
\begin{align}
	\mathcal{L}_{2,k} &=\frac{1}{2}\|\bsx_{k}\|^2_{\bsE}+\frac{1}{2}\Big\|\bsv_k+\frac{1}{\gamma_k}\bsg_k^0\Big\|^2_{\frac{\beta_k+\gamma_k}{\gamma_k}\bsF} + \|\bsx_{k}-\bsy_{k}\|^2\nonumber\\
	&\quad +\bsx_k^\top\bsE\bsF\Big(\bsv_k
		+\frac{1}{\gamma_k}\bsg_k^0\Big) \nonumber\\
	&\ge\frac{1}{2}\|\bsx_{k}\|^2_{\bsE}
		+\frac{1}{2}\Big(1+\frac{\beta_k}{\gamma_k}\Big)\Big\|\bsv_k+\frac{1}{\gamma_k}\bsg_k^0\Big\|^2_{\bsF}\nonumber\\
		&\quad-\frac{\gamma_k}{2\beta_k\rho_2(L)}\|\bsx_{k}\|^2_{\bsE}
		-\frac{\beta_k}{2\gamma_k}\Big\|\bsv_k+\frac{1}{\gamma_k}\bsg_k^0\Big\|^2_{\bsF} \nonumber\\
		&\quad +\|\bsx_{k}-\bsy_{k}\|^2\nonumber\\
	&\ge\varepsilon_7\Big(\|\bsx_{k}\|^2_{\bsE}
		+\Big\|\bsv_k+\frac{1}{\gamma_k}\bsg_k^0\Big\|^2_{\bsF}\Big)\nonumber\\
		&\quad+\|\bsx_{k}-\bsy_{k}\|^2\label{zerosg:vkLya3.2}\\
	&\ge\varepsilon_7\hat{\mathcal{L}}_{2,k}\ge0,\label{zerosg:vkLya3}
\end{align}
where the first inequality holds due to \eqref{nonconvex:lemma-eq5} and the Cauchy--Schwarz inequality; and the last inequality holds due to $0<\varepsilon_7<\frac{1}{2}$. Similarly, we have
\begin{align}\label{zerosg:vkLya3.1}
	\mathcal{L}_{2,k}\le\varepsilon_6\hat{\mathcal{L}}_{2,k}.
\end{align}

Denote $\psi_{2,k}=\mathbb{E}[\mathcal{L}_{2,k}]$. From \eqref{zerosg:sgproof-vkLya2T-bounded}, \eqref{zerosg:coro-sg-sm-gbark0}, \eqref{zerosg:vkLya3}--\eqref{zerosg:vkLya3.1}
, and \eqref{zerosg:step:eta2-sm}, we have
\begin{align}\label{zerosg:vkLya4-bound-tilde}
&\psi_{2,k+1}\le(1-\tilde{d}_1)\psi_{2,k}
+\frac{n^2(\tilde{ a}_4+\tilde{ a}_5^{\prime}\kappa_\mu^2+\tilde{ a}_6\tilde{c}_g)}{T},\nonumber\\
&\qquad\qquad\qquad\qquad\qquad\qquad\qquad ~\forall 0\le k\le T.
\end{align}

From $\epsilon_1>\kappa_1 \ge\frac{13}{2\rho_2(L)}$, we have $\varepsilon_6>1$. From $\epsilon_2\in(0,\kappa_2(\epsilon_1))$, we have $ a_2 = \frac{1}{16}\epsilon_2-\big(\frac{5}{4}+2(1+c_1^{-1})\big)\rho(L)\epsilon_2^2\le\frac{1}{16}\epsilon_2\le\frac{1}{2}$.
Thus,
\begin{align}\label{zerosg:vkLya2-a1-bounded-thm2}
0<\tilde{d}_1\le\frac{2 a_2}{\varepsilon_6}\le1.
\end{align}

From \eqref{zerosg:vkLya4-bound-tilde} and \eqref{zerosg:vkLya2-a1-bounded-thm2}, similar to the way to get \eqref{zerosg:serise:lemma:sequence-equ6}, we have $\psi_{2,T}=\mathcal{O}(\frac{n^2}{T})$. Then, from \eqref{zerosg:vkLya3.1}, we have \eqref{zerosg:coro-sg-sm-equ3.1}.

\subsection{Proof of Theorem~\ref{zerosg:thm-sg-diminishingt}}\label{zerosg:proof-thm-sg-diminishingt}
In addition to the notations defined in Appendices~\ref{zerosg:proof-thm-random-pd-sm} and \ref{zerosg:proof-thm-sg-smT}, we also denote the following notations.
\begin{align*}
&\kappa_0(\epsilon_1,\epsilon_2)=\max\Big\{\varepsilon_{0},~\frac{2\varepsilon_5}{ a_1},
	(\frac{2p(1+\eta_1^2)(1+\eta_2^2)\varepsilon_9}{ a_1})^{\frac{1}{2}},\frac{\varepsilon_{10}}{2 a_2},\\
	&
	~\frac{24(1+\eta_2^2)\varepsilon_1}{\epsilon_2},
	~\frac{\delta_0b_{6,k}}{ a_3}, ~128p(1+\eta_1^2)(1+\eta_2^2)\epsilon_2\varepsilon_{12},p\varepsilon_8\Big\},\\		
&\kappa_3(\epsilon_1,\epsilon_2)=\frac{24(1+\eta_2^2)\varepsilon_8}{\epsilon_2},\\
&\hat{\kappa}_3(\epsilon_1,\epsilon_2,\epsilon_3)\\
	&\quad=\max\Big\{\frac{\kappa_0(\epsilon_1,\epsilon_2)}{\epsilon_3},
	~\kappa_3(\epsilon_1,\epsilon_2),~\frac{2}{3 a_3},~\frac{2}{3 a_1},
	~\frac{2}{3 a_2}\Big\},\\
&\varepsilon_9
=8\big(\rho_2^{-2}(L)+4(\epsilon_1+1)^2\rho_2^{-1}(L)\big) \ell ^4\epsilon_2 \\
&\quad+\frac{(1+2 \ell ^2)\epsilon_2}{2p(1+\eta_1^2)(1+\eta_2^2)}
+\big(\frac{5}{p(1+\eta_1^2)(1+\eta_2^2)}+32\big) \ell ^2\epsilon_2^2\\
&\quad + 8\big(2\rho_2^{-2}(L)+2(\epsilon_1+1)\rho_2^{-1}(L)+3\varepsilon_8+3\varepsilon_1+3\big) \ell ^4\epsilon_2^2,\\
&\varepsilon_{10}=\rho(L)\epsilon_1\epsilon_2+3\rho(L)\epsilon_2^2+\epsilon_1+\epsilon_2+1,\\
&\varepsilon_{11}=\frac{3\epsilon_3}{2\epsilon_2^2}(\varepsilon_8+\varepsilon_1),\\
&\varepsilon_{12}=10+16(1+c_1^{-1})+ \ell \\
	&\quad +\frac{\big(4(\epsilon_1+1)^2+2\epsilon_1\epsilon_2+2\epsilon_2\big)\rho_2^{-1}(L)+\rho_2^{-2}(L)}{\epsilon_2} \ell ^2\\
	&\quad +(3+2\rho_2^{-2}(L)+3\varepsilon_8+3\varepsilon_1) \ell ^2,\\	
&\varepsilon_{15}=\frac{1}{\varepsilon_6}\min\Big\{\frac{ a_1\epsilon_3m}{\epsilon_2},\frac{ a_2\epsilon_3m}{\epsilon_2},
	~\frac{ a_3\epsilon_3m}{\epsilon_2},~\frac{\nu}{8}\Big\},\\
&\varepsilon_{16}
=\frac{4\epsilon_2^2( a_7+ \ell ^2\kappa_\mu^2)}{\epsilon_3^2(\frac{3\nu\epsilon_2}{8\epsilon_3}-1)},\\
&d_2=\frac{\epsilon_2\varepsilon_{15}}{\epsilon_3},\\
&d_3=pn( a_4+ a_5\kappa_\mu^2)\frac{\epsilon_2^2}{\kappa^2_0},\\
&d_4=pn( a_4+ a_5^{\prime}\kappa_\mu^2
+ a_6\hat{c}_g)\frac{\epsilon_2^2}{\epsilon_3^2},\\
&\hat{d}_2=\frac{2}{3\varepsilon_6}.
\end{align*}

To prove Theorem~\ref{zerosg:thm-sg-diminishingt}, the following lemma (results of Lemma \ref{zerosg:lemma:sg} under time-varying algorithm parameters) is used.

\begin{lemma}\label{zerosg:lemma:sg2}
	Suppose  Assumptions~\ref{nonconvex:ass:compression}--\ref{zerosg:ass:fig} hold. Suppose $\beta_k=\epsilon_1\gamma_k$, $\gamma_k=\epsilon_3(k+m)$, and $\alpha_k=\frac{\epsilon_2}{\gamma_k}$, $\epsilon_3\ge \frac{\kappa_0(\epsilon_1,\epsilon_2)}{m}$, $\epsilon_1>\kappa_1$, $\epsilon_2\in(0,\kappa_2(\epsilon_1))$, and $m\ge \kappa_3(\epsilon_1,\epsilon_2)$. Let $\{\bsx_k\}$ be the sequence generated by Algorithm~\ref{nonconvex:algorithm-pdgd}, then
	\begin{subequations}
	\begin{align}
	&\mathbb{E}_{\mathcal{A}_k}[\mathcal{L}_{1,k+1}] \nonumber\\
	&\quad \le  \mathcal{L}_{1,k}- a_1\|\bsx_k\|^2_{\bsE}
	- a_2\Big\|\bsv_k+\frac{1}{\gamma_k}\bsg_{k}^0\Big\|^2_{\bsF} 
	 - \frac{1}{16}\alpha_k\|\bar{\bsg}^0_{k}\|^2\nonumber\\
	&\quad- a_3\|\bsx_{k}-\bsy_{k}\|^2 
	+pn a_4\alpha_k^2+pn a_5\alpha_k\mu_k^2,
		\label{zerosg:sgproof-vkLya2}\\
	&\mathbb{E}_{\mathcal{A}_k}[\mathcal{L}_{2,k+1}] \nonumber\\
	& \quad \le \mathcal{L}_{2,k}- a_1\|\bsx_k\|^2_{\bsE}
	- a_2\Big\|\bsv_k+\frac{1}{\gamma_k}\bsg_{k}^0\Big\|^2_{\bsF} 
	 +p a_6\alpha_k^2\|\bar{\bsg}^0_{k}\|^2 \nonumber\\
	&\quad- a_3\|\bsx_{k}-\bsy_{k}\|^2 + pn a_4\alpha_k^2+pn a_5^{\prime}\alpha_k\mu_k^2,
		\label{zerosg:sgproof-vkLya2-bounded}\\
	&\mathbb{E}_{\mathcal{A}_k}[{\mathcal{L}}_{3,k+1}]\le {\mathcal{L}}_{3,k}+\|\bsx_k\|^2_{2\alpha_k \ell ^2\bsE}
	-\frac{3}{16}\alpha_k\|\bar{\bsg}_{k}^0\|^2 \nonumber\\
	&\quad +p\alpha_k^2 a_7+(n+p) \ell ^2\alpha_k\mu^2_k,
		\label{zerosg:v4kspeed-diminishing-2}
	\end{align}
	\end{subequations}
\end{lemma}
where
\begin{align*}
	& a_4=2\varepsilon_{12}\sigma^2_1+\frac{\varepsilon_{11}\sigma^2_2}{p}
	+4(1+\eta_1^2)\varepsilon_{12}\sigma^2_2,\\
	& a_5= \ell ^2\Big(\frac{1}{512(1+\eta_1^2)(1+\eta_2^2)}
		+\frac{(29+16c_1^{-1})\epsilon_2+\frac{7}{2}}{p}\Big),\\
	& a_6=\frac{b_{7,k}}{p} + \frac{(1+\eta_2^2)\varepsilon_{11}}{p}
		+8(1+\eta_1^2)(1+\eta_2^2)\varepsilon_{12},\\
	& a_5= \ell ^2\Big(\frac{1}{512(1+\eta_1^2)(1+\eta_2^2)}
		+\frac{(29+16c_1^{-1})\epsilon_2+\frac{7}{2}}{p} \\
		&\quad +2b_{7,k}\frac{\epsilon_2}{\varepsilon_0p} \Big).
\end{align*}

\begin{proof}
Noting that $\epsilon_1>\kappa_1$ and $\gamma_k=\epsilon_3(k+m)\ge\epsilon_3m\ge \kappa_0(\epsilon_1,\epsilon_2)\ge\varepsilon_{0}$, we know that all conditions needed in Lemma~\ref{zerosg:lemma:sg} are satisfied, so \eqref{zerosg:sgproof-vkLya}--\eqref{zerosg:sgproof-vkLya-bounded} hold, and \eqref{Lyapunov 3} i.e. \eqref{zerosg:v4kspeed-diminishing-2} hold.

Next, we prove \eqref{zerosg:sgproof-vkLya2} by scaling and bounding some coefficients.



Recalling $\gamma_k=\epsilon_3(k+m)$,
\begin{align}\label{zerosg:omegak}
\Delta_k&=\frac{1}{\gamma_{k}}-\frac{1}{\gamma_{k+1}}
=\frac{1}{\epsilon_3}(\frac{1}{(k+m)}-\frac{1}{(k+m+1)})
\nonumber\\
&\le\frac{1}{\epsilon_3(k+m)(k+m+1)}
\le\frac{\epsilon_3}{\gamma_k^2}.
\end{align}
From $\gamma_k=\epsilon_3(k+m)\ge\epsilon_3m\ge \kappa_0(\epsilon_1,\epsilon_2)\ge\varepsilon_{0}\ge1$, we have $\Delta_k\le1$.	

From $\epsilon_1>\kappa_1 \ge\rho_2(L)$, we have $\varepsilon_8=(1+\epsilon_1)\rho_2^{-1}(L)+\rho_2^{-2}(L)\ge1$. From  $\epsilon_2<\kappa_2(\epsilon_1)\le\frac{1}{20+32(1+c_1^{-1})}\le1$ and $\varepsilon_8\ge1$, we have $m\ge \kappa_3(\epsilon_1,\epsilon_2)=\frac{24(1+\eta_2^2)\varepsilon_8}{\epsilon_2}\ge1$.

From \eqref{zerosg:omegak}, $\alpha_k=\frac{\epsilon_2}{\gamma_k}$, $\gamma_k\ge1$, $\Delta_k\le1$, $\epsilon_3\ge \frac{\kappa_0(\epsilon_1,\epsilon_2)}{m}
\ge(\frac{2p(1+\eta_1^2)(1+\eta_2^2)\varepsilon_9}{ a_1})^{\frac{1}{2}}\times\frac{1}{m}$, and $m\ge1$, we have
\begin{align}
b_{1,k}\le\frac{p(1+\eta_1^2)(1+\eta_2^2)\varepsilon_9}{\epsilon_3^2m^2}\le\frac{ a_1}{2} .\label{zerosg:b1k}
\end{align}

From \eqref{zerosg:kappa4-6}-\eqref{zerosg:b1k}, $\epsilon_3\ge \frac{\kappa_0(\epsilon_1,\epsilon_2)}{m}\ge\frac{2\varepsilon_5}{ a_1m}$, and \eqref{zerosg:kappa4-6}, whether $\varepsilon_5>0$ or not, we have
\begin{align}
2 a_1-\varepsilon_5\Delta_k-b_{1,k}
\ge2 a_1-\frac{\varepsilon_5}{\epsilon_3m}-\frac{ a_1}{2}
\ge a_1>0.\label{zerosg:varepsilon4}
\end{align}

From \eqref{zerosg:omegak}, $\alpha_k=\frac{\epsilon_2}{\gamma_k}$, $\epsilon_3\ge \frac{\kappa_0(\epsilon_1,\epsilon_2)}{m}\ge\frac{\varepsilon_{10}}{2 a_2m}\ge\frac{\varepsilon_{10}}{2 a_2m}$, and \eqref{zerosg:kappa4-6}, we have
\begin{align}
b_{2,k}\ge2 a_2-\frac{\varepsilon_{10}}{2\epsilon_3m}
\ge a_2>0.\label{zerosg:b2k}
\end{align}

From \eqref{zerosg:kappa4-6}, \eqref{zerosg:omegak}, $m\ge\kappa_3(\epsilon_1,\epsilon_2)\ge(\frac{24(1+\eta_2^2)\varepsilon_8}{\epsilon_2})^\frac{1}{\theta}$, $\epsilon_3\ge \frac{\kappa_0(\epsilon_1,\epsilon_2)}{m}\ge\frac{24(1+\eta_2^2)\varepsilon_1}{\epsilon_2m}$ , and $\alpha_k=\frac{\epsilon_2}{\gamma_k}$ ,we have

\begin{align}\label{zerosg:b3keta}
	b_{3,k}\alpha_k\le\frac{3\varepsilon_8}{2\epsilon_2m}
+\frac{3\varepsilon_1}{2\epsilon_2\epsilon_3m}\le\frac{1}{8(1+\eta_2^2)}.
\end{align}

From $\eqref{zerosg:omegak}, \gamma_k\ge1$ and $\Delta_k\le1$, we have
\begin{subequations}
\begin{align}
b_{3,k}&\le \varepsilon_{11},\label{zerosg:b3k}\\
b_{4,k}&\le \varepsilon_{12}.\label{zerosg:b4k}
\end{align}
\end{subequations}

From $\epsilon_3\ge \frac{\kappa_0(\epsilon_1,\epsilon_2)}{m}\ge\frac{128p(1+\eta_1^2)(1+\eta_2^2)\epsilon_2\varepsilon_{12}}{m}$, \eqref{zerosg:b3keta} and \eqref{zerosg:b4k}, we have
	\begin{align}\label{zerosg:b3kb4keta}
	&\frac{1}{4}-(1+\eta_2^2)(b_{3,k}+8p(1+\eta_1^2)b_{4,k})\alpha_k\nonumber\\
	&\ge\frac{1}{8}-8p(1+\eta_1^2)(1+\eta_2^2)b_{4,k}\alpha_k \nonumber\\
	&\ge\frac{1}{8}-\frac{8p(1+\eta_1^2)(1+\eta_2^2)\epsilon_2\varepsilon_{12}}{\epsilon_3m}
	\ge\frac{1}{16},\\
	&b_{4,k}\alpha_k\le\frac{1}{128p(1+\eta_1^2)(1+\eta_2^2)}.\label{mid 2}
	\end{align}

From \eqref{mid 2}, $\alpha_k=\frac{\epsilon_2}{\gamma_k}$, $\gamma_k\ge1$, and $\Delta_k\le1$, we have
\begin{align}
b_{5,k}\le pn a_5.\label{zerosg:b5k}
\end{align}


From \eqref{zerosg:omegak} and $\epsilon_3\ge \frac{\kappa_0(\epsilon_1,\epsilon_2)}{m}\ge\frac{\delta_0b_{6,k}}{ a_3m}$, we have
\begin{align}\label{zerosg:c2}
	\delta_0b_{6,k}\Delta_k\le \frac{\delta_0b_{6,k}}{\epsilon_3m}\le  a_3.
\end{align}

From \eqref{zerosg:c1} and \eqref{zerosg:c2}, we have
\begin{align}\label{zerosg:c}
	c_2-\delta_0\varepsilon_2-\delta_0b_{6,k}\Delta_k\ge  a_3>0.
\end{align}

Then, from \eqref{zerosg:sgproof-vkLya}, \eqref{zerosg:varepsilon4}, \eqref{zerosg:b2k}, \eqref{zerosg:b3k}--\eqref{zerosg:b5k} and \eqref{zerosg:c}, we know that \eqref{zerosg:sgproof-vkLya2} holds.
\hspace*{\fill}

 From \eqref{zerosg:sgproof-vkLya-bounded}, \eqref{zerosg:varepsilon4}, \eqref{zerosg:b2k}, \eqref{zerosg:b3k}, \eqref{zerosg:b4k}, \eqref{zerosg:b5k} and \eqref{zerosg:c}, we have \eqref{zerosg:sgproof-vkLya2-bounded}.

\end{proof}

Now it is ready to prove Theorem~\ref{zerosg:thm-sg-diminishingt}.

From $m\ge\hat{\kappa}_3(\epsilon_1,\epsilon_2,\epsilon_3)\ge\frac{\kappa_0(\epsilon_1,\epsilon_2)}{\epsilon_3}$, we have
\begin{align*}
\epsilon_3\ge\frac{\kappa_0(\epsilon_1,\epsilon_2)}{m}.
\end{align*}
Thus, all conditions needed in Lemma~\ref{zerosg:lemma:sg2} are satisfied, so \eqref{zerosg:sgproof-vkLya2}--\eqref{zerosg:v4kspeed-diminishing-2} still hold.

From \eqref{nonconvex:equ:plc}, we have that
\begin{align}\label{nonconvex:gg3}
\|\bar{\bsg}^0_k\|^2=n\|\nabla f(\bar{x}_k)\|^2\ge2\nu n(f(\bar{x}_k)-f^*)
=2\nu e_{4,k}.
\end{align}

From \eqref{zerosg:sgproof-vkLya2}, \eqref{nonconvex:gg3}, \eqref{zerosg:vkLya3}, and \eqref{zerosg:vkLya3.1}, we have
\begin{align}\label{zerosg:vkLya2-pl}
&\mathbb{E}_{\mathcal{A}_k}[\mathcal{L}_{1,k+1}]\nonumber\\
&\le \mathcal{L}_{1,k}- a_1\|\bsx_k\|^2_{\bsE}
 - a_2\Big\|\bsv_k+\frac{1}{\gamma_k}\bsg_{k}^0\Big\|^2_{\bsF} 
 -\frac{\alpha_k\nu}{8}e_{4,k}\nonumber\\
&\quad - a_3\|\bsx_{k}-\bsy_{k}\|^2 + pn a_4\alpha_k^2+pn a_5\alpha_k\mu_k^2\nonumber\\
&\le \mathcal{L}_{1,k}-\frac{\alpha_k}{\varepsilon_6}\min\Big\{\frac{ a_1}{\alpha_k}, \frac{ a_2}{\alpha_k},
	\frac{a_3}{\alpha_k}, \frac{\nu}{8}\Big\}\mathcal{L}_{1,k}\nonumber\\
 &\quad +pn a_4\alpha_k^2+pn a_5\alpha_k\mu_k^2\nonumber\\
&\le \mathcal{L}_{1,k}-\alpha_k\varepsilon_{15}\mathcal{L}_{1,k}+pn a_4\alpha_k^2
+pn a_5\alpha_k\mu_k^2,~\forall k\in\mathbb{N}_0.
\end{align}

Denote $\hat{\mathcal{L}}_{1,k}=\|\bsx_k\|^2_{\bsE}+\Big\|\bsv_k+\frac{1}{\gamma}\bsg_k^0\Big\|^2_{\bsF}+n(f(\bar{x}_k)-f^*)+\|\bsx_{k}-\bsy_{k}\|^2$.
From \eqref{zerosg:vkLya3} and \eqref{zerosg:vkLya3.1}, we have
\begin{align}
	&\mathcal{L}_{1,k}\ge\varepsilon_7\hat{\mathcal{L}}_{1,k}\ge0,\label{ge7}\\
	&\mathcal{L}_{1,k}\le\varepsilon_6\hat{\mathcal{L}}_{1,k}.\label{le6}
\end{align}

Denote $\psi_{1,k}=\mathbb{E}[\mathcal{L}_{1,k}]$, $r_{1,k}=\alpha_k\varepsilon_{15}$, and $r_{2,k}=pn a_4\alpha_k^2
+pn a_5\alpha_k\mu_k^2$. From \eqref{zerosg:vkLya2-pl}, we have
\begin{align}
\psi_{1,k+1}
\le (1-r_{1,k})\psi_{1,k}+r_{2,k},~\forall k\in\mathbb{N}_0.
\label{zerosg:vkLya2-pl-z}
\end{align}

From \eqref{zerosg:step:eta1t1}, we have
\begin{align}
r_{1,k}&=\alpha_k\varepsilon_{15}
=\frac{d_2}{(k+m)},\label{zerosg:vkLya2-pl-r1}\\
r_{2,k}&=pn a_4\alpha_k^2
+pn a_5\alpha_k\mu_k^2
\le\frac{d_3}{(k+m)}.\label{zerosg:vkLya2-pl-r2}
\end{align}

From \eqref{zerosg:vkLya2-a1-bounded-thm2}, we have
\begin{align}\label{zerosg:vkLya2-pl-r1.1}
r_{1,k}\le\frac{ a_2}{\varepsilon_6}\le\frac{1}{2}\le1.
\end{align}

From \eqref{zerosg:kappa4-6} and \eqref{zerosg:b2k}, we know that
\begin{align}\label{zerosg:vkLya2-pl-a1a2}
d_2>0~\text{and}~d_3>0.
\end{align}

Due to $m\ge\hat{\kappa}_3(\epsilon_1,\epsilon_2,\epsilon_3)\ge \frac{2}{3 a_3}$,
\begin{align}\label{zerosg:vkLya2-pl-a1-1}
\frac{ a_3m}{\varepsilon_6}\ge\frac{2}{3\varepsilon_6}.
\end{align}

From $m\ge\hat{\kappa}_3(\epsilon_1,\epsilon_2,\epsilon_3)\ge \frac{2}{3 a_1}$, we have
\begin{align}\label{zerosg:vkLya2-pl-a1-1d}
\frac{ a_1m}{\varepsilon_6}\ge\frac{2}{3\varepsilon_6}.
\end{align}

From $m\ge\hat{\kappa}_3(\epsilon_1,\epsilon_2,\epsilon_3)\ge \frac{2}{3 a_2}$, we have
\begin{align}\label{zerosg:vkLya2-pl-a1-2}
\frac{ a_2m}{\varepsilon_6}\ge\frac{2}{3\varepsilon_6}.
\end{align}

Because $\epsilon_3\in[\frac{3\hat{\kappa}_0\nu\epsilon_2}{16},\frac{3\nu\epsilon_2}{16})$,
\begin{align}\label{zerosg:vkLya2-pl-a1-3}
\frac{16}{3\nu}<\frac{\epsilon_2}{\epsilon_3}\le\frac{16}{3\hat{\kappa}_0\nu}.
\end{align}
Thus,
\begin{align}\label{zerosg:vkLya2-pl-a1-2.1}
\frac{\nu\epsilon_2}{8\varepsilon_6\epsilon_3}>\frac{2}{3\varepsilon_6}.
\end{align}

Hence, from \eqref{zerosg:vkLya2-pl-a1-1}, \eqref{zerosg:vkLya2-pl-a1-2}, \eqref{zerosg:vkLya2-pl-a1-2.1}, and $\varepsilon_6>1$ due to $\epsilon_1>1$, we have
\begin{align}\label{zerosg:vkLya2-pl-a1}
d_2>\hat{d}_2>0~\text{and}~\hat{d}_2<\frac{2}{3}.
\end{align}

Then from \eqref{zerosg:vkLya2-pl-z}--\eqref{zerosg:vkLya2-pl-a1a2}, \eqref{zerosg:vkLya2-pl-a1}, and \eqref{zerosg:serise:lemma:sequence-equ5},
\begin{align}\label{zerosg:vkLya2-pl-theta0t}
\psi_{1,k}\le\phi_1(k,m,\hat{d}_2,d_3,2,\psi_{1,0}),~\forall k\in\mathbb{N}_+,
\end{align}
where the function $\phi_1$ is defined in \eqref{zerosg:serise:lemma:sequence-equ5-phi3}.

From \eqref{zerosg:vkLya2-pl-a1}, we have $\phi_1(k,m,\hat{d}_2,d_3,2,z_0)
=\mathcal{O}(\frac{nm^{\hat{d}_2}}{(k+m)^{\hat{d}_2}}
+\frac{pn}{(k+m)^{\hat{d}_2}m^{1-\hat{d}_2}})$. Hence, from \eqref{zerosg:rand-grad-smooth}, \eqref{zerosg:vkLya2-pl-theta0t}, we get
\begin{align}\label{zerosg:gbark0-pl-speed-1}
\mathbb{E}[\|\bar{\bsg}^0_k\|^2]
=\mathcal{O}\Big(\frac{nm^{\hat{d}_2}}{(k+m)^{\hat{d}_2}}
+\frac{pn}{(k+m)^{\hat{d}_2}m^{1-\hat{d}_2}}\Big),~\forall k\in\mathbb{N}_+.
\end{align}

Noting that $m>\hat{\kappa}_3(\epsilon_1,\epsilon_2,\epsilon_3)\ge \frac{\kappa_0(\epsilon_1,\epsilon_2)}{\epsilon_3}
\ge\frac{p\varepsilon_8}{\epsilon_3}$, i.e. $m=\mathcal{O}(p)$, from \eqref{zerosg:gbark0-pl-speed-1}, we know that there exists a constant $\hat{c}_g>0$, such that
\begin{align}\label{zerosg:gbark0-pl-speed}
\mathbb{E}[\|\bar{\bsg}^0_k\|^2]\le n\hat{c}_g,~\forall k\in\mathbb{N}_0.
\end{align}

From \eqref{zerosg:sgproof-vkLya2-bounded}, \eqref{zerosg:gbark0-pl-speed}, \eqref{zerosg:vkLya3}, and \eqref{zerosg:step:eta1t1}, we have
\begin{align}\label{zerosg:vkLya4-bound-pl-1}
\psi_{2,k+1}\le(1-\tilde{d}_1)\psi_{2,k}+\frac{d_4}{(t+m)^{2}}.
\end{align}

Using \eqref{zerosg:serise:lemma:sequence-equ6}, from \eqref{zerosg:vkLya2-a1-bounded-thm2} and \eqref{zerosg:vkLya4-bound-pl-1}, we have
\begin{align}\label{zerosg:lemma:sequence-equ6-bounded-pl-speed}
\psi_{2,k}&\le \phi_2(k,m,\tilde{d}_1,d_4,2,\psi_{2,0}),~\forall k\in\mathbb{N}_+,
\end{align}
where the function $\phi_2$ is defined in \eqref{zerosg:serise:lemma:sequence-equ6-phi4}.
From \eqref{zerosg:lemma:sequence-equ6-bounded-pl-speed}, \eqref{ge7}, \eqref{zerosg:serise:lemma:sequence-equ6-phi4}, and \eqref{zerosg:vkLya2-pl-a1-3}, we have
\begin{align}\label{zerosg:vkLya4-bound-brevez-speed}
&\mathbb{E}[\|\bsx_k\|^2_{\bsE}]\le\frac{1}{\varepsilon_7}\psi_{2,k}\nonumber\\
&\le\frac{1}{\varepsilon_7}\phi_2(k,m,\tilde{d}_1,d_4,2,\psi_{2,0})
=\mathcal{O}\Big(\frac{pn}{(k+m)^{2}}\Big),
\end{align}
which yields \eqref{zerosg:thm-sg-diminishing-equ2.1bounded}.

From $\epsilon_3<\frac{3\nu\epsilon_2}{16}$, we have
\begin{align}\label{zerosg:vkLya2-pl-a1-3-bounded}
\frac{3\nu\epsilon_2}{8\epsilon_3}>2.
\end{align}

Same to the way to prove \eqref{zerosg:serise:lemma:sequence-equ5}, from \eqref{zerosg:v4kspeed-diminishing-2}, \eqref{zerosg:vkLya4-bound-brevez-speed}, \eqref{zerosg:vkLya2-pl-a1-3-bounded}, and \eqref{zerosg:vkLya2-pl-a1-3}, we have \begin{align}\label{zerosg:v4kspeed-diminishing-6}
&\mathbb{E}[f(\bar{x}_{T})-f^*]\nonumber\\
&\le\mathcal{O}\Big(\frac{m^{\frac{3\nu\epsilon_2}{8\epsilon_3}}}{(T+m)^{\frac{3\nu\epsilon_2}{8\epsilon_3}}}\Big) +\mathcal{O}\Big(\frac{p}{n(T+m)^{2}}\Big)+\frac{\varepsilon_{16}p}{n(T+m)}\nonumber\\
&+\mathcal{O}\Big(\frac{m^{\frac{3\nu\epsilon_2}{8\epsilon_3}}}{(T+m)^{\frac{3\nu\epsilon_2}{8\epsilon_3}}}\Big) +\mathcal{O}\Big(\frac{p}{(T+m)^{3}}\Big) +\mathcal{O}\Big(\frac{p}{(T+m)^{2}}\Big).
\end{align}
From \eqref{zerosg:vkLya2-pl-a1-3}, we have
\begin{align}\label{zerosg:v4kspeed-diminishing-7}
&\varepsilon_{16}
=\frac{4\epsilon_2^2( a_7+ \ell ^2\kappa_\mu^2)}{\epsilon_3^2(\frac{3\nu\epsilon_2}{8\epsilon_3}-1)}
\le\frac{1024( a_7+ \ell ^2\kappa_\mu^2)}{9\hat{\kappa}_0(2-\hat{\kappa}_0)\nu^2}.
\end{align}
Thus, from \eqref{zerosg:v4kspeed-diminishing-6}, \eqref{zerosg:v4kspeed-diminishing-7}, and $m=\mathcal{O}(p)$, we have \eqref{zerosg:thm-sg-diminishing-equ2bounded}.

\subsection{Proof of Theorem~\ref{zerosg:thm-random-pd-fixed}}\label{zerosg:proof-thm-random-pd-fixed}

In addition to the notations defined in Appendix~\ref{zerosg:proof-thm-sg-smT},
we also denote the following notations.
\begin{align*}
&\varepsilon=\frac{1}{2}+\frac{1}{2}\max\{1-\tilde{\varepsilon}_{15},~\tilde{\varepsilon}^2\},\\
&\tilde{\varepsilon}_{15}=\frac{1}{4\varepsilon_6}\min\{4a_1, 8a_2, 8a_3, \alpha\nu\}.
\end{align*}

All conditions needed in Lemma~\ref{zerosg:lemma:sg2-T} are satisfied, so \eqref{zerosg:sgproof-vkLya2T} still holds.

If Assumption~\ref{nonconvex:ass:fil} also holds, then \eqref{nonconvex:gg3} holds.
From  \eqref{zerosg:sgproof-vkLya2T}, \eqref{nonconvex:gg3},  and \eqref{zerosg:vkLya3.1}, for any $k\in\mathbb{N}_0$, we have
\begin{align}\label{zerosg:vkLya2-pl-fixed}
&\mathbb{E}[\mathcal{L}_{1,k+1}]\nonumber\\
&\le \mathbb{E}\Big[\mathcal{L}_{1,k}- a_1\|\bsx_k\|^2_{\bsE}
-2 a_2\|\bsv_k+\frac{1}{\gamma}\bsg_{k}^0\|^2_{\bsE}\nonumber\\
&\quad -\frac{\alpha\nu }{4}e_{4,k}-2 a_3\|\bsx_{k}-\bsy_{k}\|^2\nonumber\\
&\quad
+2pn(\sigma^2_1+2(1+\eta_1^2)\sigma^2_2)\tilde{\varepsilon}_{12}\alpha^2+pn\tilde{ a}_5\alpha\mu_k^2\nonumber\\
&\le \mathbb{E}[1-\tilde{\varepsilon}_{15}\mathcal{L}_{1,k}]\nonumber\\
&\quad+2pn(\sigma^2_1+2(1+\eta_1^2)\sigma^2_2)\tilde{\varepsilon}_{12}\alpha^2
+pn\tilde{ a}_5\alpha\mu_k^2.
\end{align}

From \eqref{zerosg:vkLya2-a1-bounded-thm2}, we have
\begin{align}\label{zerosg:vkLya2-pl-r1.1-fixed}
0<\tilde{\varepsilon}_{15}\le\frac{2 a_2}{\varepsilon_6}\le1.
\end{align}

From \eqref{zerosg:vkLya2-pl-fixed}, \eqref{zerosg:vkLya3},  \eqref{zerosg:vkLya2-pl-r1.1-fixed}, and $\mu_{i,k}\in(0,\kappa_\mu\tilde{\varepsilon}^{k}]$, we have
\begin{align}\label{zerosg:vkLya2-pl-fixed2}
&~\mathbb{E}[\mathcal{L}_{1,k+1}]\nonumber\\
&\le (1-\tilde{\varepsilon}_{15})^{k+1}{\mathcal{L}}_{1,0}\nonumber\\
&+2pn(\sigma^2_1+2(1+\eta_1^2)\sigma^2_2)\tilde{\varepsilon}_{12}\alpha^2
\sum_{\tau=0}^{k}(1-\tilde{\varepsilon}_{15})^\tau\nonumber\\
&+pn\tilde{ a}_5\kappa_\mu^2\alpha
\sum_{\tau=0}^{k}(1-\tilde{\varepsilon}_{15})^\tau
\tilde{\varepsilon}^{2(k-\tau)},~\forall k\in\mathbb{N}_0.
\end{align}

From \eqref{zerosg:vkLya3}, \eqref{zerosg:vkLya2-pl-fixed2}, \eqref{zerosg:lemma:sumgeo-equ}, and $\varepsilon>\max\{1-\tilde{\varepsilon}_{15},~\tilde{\varepsilon}^2\}$,  we have
\begin{align}\label{zerosg:vkLya2-pl-fixed3}
&\mathbb{E}[\|\bsx_k\|^2_{\bsE}+n(f(\bar{x}_k)-f^*)]
\le\frac{1}{\varepsilon_7}\mathbb{E}[\mathcal{L}_{1,k}]\nonumber\\
&\le \frac{n}{\varepsilon_7}\Big(\frac{{\mathcal{L}}_{1,0}}{n}
+\frac{p\tilde{ a}_5\kappa_\mu^2\alpha}{\varepsilon-\tilde{\varepsilon}^2}\Big)\varepsilon^{k}\nonumber\\
&+\frac{2n\tilde{\varepsilon}_{12}\alpha}{\varepsilon_7\tilde{\varepsilon}_{15}}
(\sigma^2_1+2(1+\eta_1^2)\sigma^2_2)p\alpha,~\forall k\in\mathbb{N}_0,
\end{align}
which gives \eqref{zerosg:thm-sg-fixed-equ1-coro1} under the assumption that $\sigma_1 = \sigma_2 = 0$.

\end{document}